\documentclass{article}
\usepackage{amsmath,amsthm,mathtools,amsfonts,amssymb,stmaryrd,xy,dsfont,mathabx}
\usepackage[a4paper]{geometry}
\usepackage[utf8]{inputenc}
\theoremstyle{definition}
\newtheorem{remark}{Remark}[section]
\newtheorem{remarks}[remark]{Remarks}

 \theoremstyle{plain}
\newtheorem{definition}[remark]{Definition}
\newtheorem{theorem}[remark]{Theorem}
\newtheorem{proposition}[remark]{Proposition}
\newtheorem{corollary}[remark]{Corollary}
\newtheorem{lemma}[remark]{Lemma}

\newtagform{empty}{}{}
\newcommand{\complex}{\mathbb{C}}
\newcommand{\reals}{\mathbb{R}}

\DeclareMathOperator{\Id}{id}

\DeclareMathOperator{\Ad}{Ad}

\DeclareMathOperator{\fix}{Fix}
\DeclareMathOperator{\cofix}{Cofix}

\DeclareMathOperator{\Mor}{Mor}

\newcommand{\smalldiagram}{}

\newcommand{\hDelta}{\widehat{\Delta}}
\newcommand{\hepsilon}{\widehat{\epsilon}}
\newcommand{\halpha}{\widehat{\alpha}}

\newcommand{\mutimes}{\underset{\mu}{\otimes}}
\newcommand{\nutimes}{\underset{\nu}{\otimes}}
\newcommand{\stensor}[2]{{_{#1}\otimes_{#2}}} 
\newcommand{\sfsource}{\stensor{\hbeta}{\alpha}}
\newcommand{\sfrange}{\stensor{\alpha}{\beta}}

\newcommand{\lnspan}{\big[} \newcommand{\rnspan}{\big]}

\newcommand{\frakA}{\mathfrak{A}}

\newcommand{\frakB}{\mathfrak{B}}

\newcommand{\frakC}{\mathfrak{C}}

\newcommand{\frakH}{\mathfrak{H}}
\newcommand{\frakK}{\mathfrak{K}}
\newcommand{\frakL}{\mathfrak{L}}

\newcommand{\cbasel}[2]{(\mathfrak{#1}, \mathfrak{#2}, \mathfrak{#1}^{\dag})}

\newcommand{\lt}{\smalltriangleleft}
\newcommand{\rt}{\smalltriangleright}

\newcommand{\hbeta}{\widehat{\beta}}
\newcommand{\hA}{\widehat{A}} \newcommand{\ha}{\widehat{a}}
 
\newcommand{\hR}{\widehat{R}}

\newcommand{\mycong}{\xrightarrow{\cong}}

\newcommand{\rtensor}[3]{ {_{#1}\! \underset{#2}{\otimes}\! {}_{#3}}}

\newcommand{\htensor}[2]{\rtensor{#1}{\frakH}{#2}}

\newcommand{\rtensorab}{\htensor{\alpha}{\beta}}

\newcommand{\rtensorh}{\underset{\frakH}{\otimes}}

\newcommand{\fibre}[3]{ {_{#1}\! \underset{#2}{\ast}\! {}_{#3}}}

\newcommand{\rfibre}[2]{ {_{#1}\ast_{#2}}}

\newcommand{\fsource}{\htensor{\hbeta}{\alpha}}
\newcommand{\frange}{\rtensorab}

\newcommand{\tl}{\ensuremath \olessthan}
\newcommand{\tr}{\ensuremath \ogreaterthan}

\newcommand{\Hsource}{H \fsource H}
\newcommand{\Hrange}{H \frange H}

\newcommand{\HsrH}{H \rtensor{\tilde s}{\tilde \mu}{\tilde r} H}

\newcommand{\checkHsource}{H \stensor{\halpha}{\hbeta} H}
\newcommand{\checkHrange}{H \stensor{\hbeta}{\alpha} H}
\newcommand{\hatHsource}{H \stensor{\alpha}{\beta} H}
\newcommand{\hatHrange}{H \stensor{\beta}{\halpha} H}

\newcommand{\checkV}{\widecheck{V}}
\newcommand{\hatV}{\widehat{V}}

\newcommand{\sHsource}{H \sfsource H}
\newcommand{\sHrange}{H \sfrange H}

\newcommand{\gHsource}{K \stensor{\sigma}{\rho} K}
\newcommand{\gHrange}{K \stensor{\rho}{\rho} K}

\newcommand{\fibreab}{\fibre{\alpha}{}{\beta}}
\newcommand{\AfibreA}{A \fibreab A}

\newcommand{\kalpha}[1]{|\alpha\rangle_{\leg{#1}}}
\newcommand{\balpha}[1]{\langle\alpha|_{\leg{#1}}}
\newcommand{\kbeta}[1]{|\beta{}\rangle_{\leg{#1}}}
\newcommand{\bbeta}[1]{\langle\beta|_{\leg{#1}}}
\newcommand{\khbeta}[1]{|\hbeta{}\rangle_{\leg{#1}}}
\newcommand{\bhbeta}[1]{\langle\hbeta|_{\leg{#1}}}
\newcommand{\khalpha}[1]{|\halpha{}\rangle_{\leg{#1}}}

\newcommand{\kgamma}[1]{|\gamma\rangle_{\leg{#1}}}
\newcommand{\bgamma}[1]{\langle\gamma|_{\leg{#1}}}

\DeclareMathOperator{\Dom}{Dom}

\newcommand{\leg}[1]{#1}

\newcommand{\Fix}{\mathrm{Fix}}
\newcommand{\Cofix}{\mathrm{Cofix}}

\newcommand{\pmu}{(H,\hat\beta,\alpha,\beta,V)}

\newcommand{\sHone}{H \sfsource H \sfsource H}
\newcommand{\sHtwo}{H \sfrange H \sfsource  H}
\newcommand{\sHthree}{H \sfrange H \sfrange H}
\newcommand{\sHfour}{(\sHsource) \stensor{\alpha \lt \alpha}{\beta} H}
\newcommand{\sHfourlt}{H \sfsource H \rtensor{\beta}{}{\alpha} H}
\newcommand{\sHfourrt}{\big(\sHrange\big) \rtensor{\hbeta \lt
    \beta}{}{\alpha} H}
\newcommand{\sHfive}{H \stensor{\hbeta}{\alpha \rt \alpha} (\sHrange)}

\newcommand{\gtensor}{\stensor{\gamma^{op}}{\gamma}}

\newcommand{\ckac}{(H,\alpha,\halpha,\beta,\hbeta,U,V)}

\title{Compact quantum groupoids in the setting of $C^{*}$-algebras}

\author{Thomas Timmermann\\[1ex]
  \texttt{timmermt@math.uni-muenster.de}\\ NWF-I Mathematik,
  Universität Regensburg\\ 93040 Regensburg}

\date{\today}

\usepackage[dvips]{hyperref}
\begin{document}
 \xyrequire{matrix}
\xyrequire{arrow}
\xyrequire{curve}

\maketitle 

\abstract{We propose a definition of compact quantum groupoids in
  the setting of $C^{*}$-algebras, associate to such a quantum
  groupoid a regular $C^{*}$-pseudo-multiplicative unitary, and use
  this unitary to construct a dual Hopf $C^{*}$-bimodule and to pass
  to a measurable quantum groupoid in the sense of Enock and
  Lesieur. Moreover, we discuss examples related to compact and to
  \'etale groupoids and study principal compact $C^{*}$-quantum
  groupoids. }

\tableofcontents

\section{Introduction}

\paragraph{Overview}
In the setting of von Neumann algebras, measurable quantum groupoids
were studied by Lesieur and Enock
\cite{enock:1,enock:2,enock:9,lesieur}, building Vallin's notions of 
Hopf-von Neumann bimodules and  pseudo-multiplicative unitaries
\cite{vallin:1,vallin:2} and Haagerup's operator-valued
weights.  

In this article, we propose a definition of compact quantum
groupoids in the setting of $C^{*}$-algebras, building on the
notions of Hopf-$C^{*}$-bimodules and $C^{*}$-pseudo-multiplicative
unitaries introduced in \cite{timmer:cpmu,timmer:ckac}. To each
compact $C^{*}$-quantum groupoid, we associate a regular
$C^{*}$-pseudo-multiplicative unitary, a von-Neumann algebraic
completion, and a dual Hopf $C^{*}$-bimodule. Moreover, we extend this
$C^{*}$-pseudo-multiplicative unitary to a weak $C^{*}$-pseudo-Kac
system; hence, the results of \cite{timmer:ckac} can be applied to
the study of coactions of compact $C^{*}$-quantum groupoids on
$C^{*}$-algebras.

To illustrate the general theory, we discuss several examples of compact
$C^{*}$-quantum groupoids: the $C^{*}$-algebra of continuous
functions on a compact groupoid, the reduced $C^{*}$-algebra of
an \'etale groupoid with compact base, and
principal  compact $C^{*}$-quantum groupoids.

Let us mention that many constructions and results seem to extend to a more
general notion of quantum $C^{*}$-groupoids where the Haar weights
are still assumed to be bounded but where the $C^{*}$-algebra of
units need no longer be unital and where the KMS-state on this
$C^{*}$-algebra is replaced by a proper KMS-weight. 

\paragraph{Plan}
Let us outline the contents and organization of this article in some
more detail. 

In the first part of this article (Sections 2,3,4), we introduce the
definition of a compact $C^{*}$-quantum groupoid.  Recall that a
measured compact groupoid consists of a base space $G^{0}$, a total
space $G$, range and source maps $r,s \colon G \to G^{0}$, a
multiplication $G {_{s}\times_{r}} G \to G$, a left and a right Haar
system, and a quasi-invariant measure on $G^{0}$.  Roughly, the
corresponding ingredients of a compact $C^{*}$-quantum groupoid are
unital $C^{*}$-algebras $B$ and $A$, representations $r,s \colon
B^{(op)} \to A$, a comultiplication $\Delta \colon A \to A \ast A$,
a left and a right Haar weight $\phi,\psi \colon A \to B^{(op)}$,
and a KMS-state on $B$, subject to several axioms.  These
ingredients are introduced in several steps. In Section 2, we focus
on the tuple $(B,A,r,\phi,s,\psi)$, which can be considered as a
compact $C^{*}$-quantum graph, and review some related
GNS-constructions. In Section 3, we recall from
\cite{timmer:cpmu,timmer:ckac} the definition of the fiber product
$A \ast A$ and of the underlying relative tensor product of Hilbert
modules over $C^{*}$-algebras.  Finally, in Section 4, we give the
definition of a compact $C^{*}$-quantum groupoid and obtain first
properties like uniqueness of the Haar weights and the existence of
an invariant state on the basis.

In the second part of this article (Sections 5,6,7), we associate to
every compact $C^{*}$-quantum groupoid a fundamental unitary, a
von-Neumann-algebraic completion, and a dual Hopf
$C^{*}$-bimodule. This fundamental unitary satisfies a pentagon
equation, generalizes the multiplicative unitaries of Baaj and
Skandalis \cite{baaj:2}, and can be considered as a
pseudo-multiplicative unitary in the sense of Vallin \cite{vallin:2}
equipped with additional data. The unitary and the completion will
be constructed in Section 5.  In Section 6, we study a
particular feature of this unitary --- the existence of fixed or
cofixed elements --- and show that for a general
$C^{*}$-pseudo-multiplicative unitary, such (co)fixed elements yield
invariant conditional expectations and bounded counits on the
associated Hopf $C^{*}$-bimodules. In Section 7, we return to
compact $C^{*}$-quantum groupoids and discuss their duals.

The last part of this article (Sections 8,9) is devoted to examples
of compact $C^{*}$-quantum groupoids which are obtained from compact
and from \'etale groupoids one side and from center-valued traces on
$C^{*}$-algebras on the other side. For these examples, we give a
detailed description of the ingredients, the associated fundamental
unitaries, and the dual objects.

\paragraph{Preliminaries}
Let us fix some general notation and concepts used in this article.

Given a subset $Y$ of a normed space $X$, we denote by $[Y]
\subseteq X$ the closed linear span of $Y$.  Given a $C^{*}$-algebra
$A$ and a $C^{*}$-subalgebra $B \subseteq M(A)$, we denote by $A
\cap B'$ the relative commutant $\{a \in A \mid ab=ba$ for all $b
\in B\}$. Given a Hilbert space $H$ and a subset $X \subseteq {\cal
  L}(H)$, we denote by $X'$ the commutant of $X$.  All sesquilinear
maps like inner products of Hilbert spaces are assumed to be
conjugate-linear in the first component and linear in the second
one.

We shall make extensive use of Hilbert $C^{*}$-modules. A standard
reference is \cite{lance}.

Let $A$ and $B$ be $C^{*}$-algebras.  Given Hilbert $C^{*}$-modules $E$
and $F$ over $B$, we denote the space of all adjointable operators $E\to
F$ by ${\cal L}_{B}(E,F)$.  Let $E$ and $F$ be $C^{*}$-modules over $A$
and $B$, respectively, and let $\pi \colon A \to {\cal L}_{B}(F)$ be a
$*$-homomorphism. Then one can form the internal tensor product $E
\otimes_{\pi} F$, which is a Hilbert $C^{*}$-module over $B$
\cite[Chapter 4]{lance}. This Hilbert $C^{*}$-module is the closed
linear span of elements $\eta \otimes_{A} \xi$, where $\eta \in E$ and
$\xi \in F$ are arbitrary, and $\langle \eta \otimes_{\pi} \xi|\eta'
\otimes_{\pi} \xi'\rangle = \langle
\xi|\pi(\langle\eta|\eta'\rangle)\xi'\rangle$ and $(\eta \otimes_{\pi}
\xi)b=\eta \otimes_{\pi} \xi b$ for all $\eta,\eta' \in E$, $\xi,\xi'
\in F$, and $b \in B$.  We denote the internal tensor product by
``$\tr$''; thus, for example, $E \tr_{\pi} F=E \otimes_{\pi} F$. If the
representation $\pi$ or both $\pi$ and $A$ are understood, we write
``$\tr_{A}$'' or ``$\tr$'', respectively, instead of $"\tr_{\pi}$''.

Given $A$, $B$, $E$, $F$ and $\pi$ as above, we define a {\em
  flipped internal tensor product} $F {_{\pi}\tl} E$ as follows. We
equip the algebraic tensor product $F \odot E$ with the structure
maps $\langle \xi \odot \eta | \xi' \odot \eta'\rangle := \langle
\xi| \pi(\langle \eta|\eta'\rangle) \xi'\rangle$, $(\xi \odot \eta)
b := \xi b \odot \eta$, and by factoring out the null-space of the
semi-norm $\zeta\mapsto \| \langle \zeta|\zeta\rangle\|^{1/2}$ and
taking completion, we obtain a Hilbert $C^{*}$-$B$-module $F
{_{\pi}\tl} E$.  This is the closed linear span of elements $\xi
{_{\pi}\tl} \eta$, where $\eta \in E$ and $\xi \in F$ are arbitrary,
and $\langle \xi {_{\pi}\tl} \eta|\xi' {_{\pi}\tl} \eta'\rangle =
\langle \xi|\pi(\langle\eta|\eta'\rangle)\xi'\rangle$ and $(\xi
{_{\pi}\tl} \eta)b=\xi b {_{\pi}\tl} \eta$ for all $\eta,\eta' \in
E$, $\xi,\xi' \in F$, and $b\in B$. As above, we write
``${_{A}\tl}$'' or simply ``$\tl$'' instead of ``${_{\pi}\tl}$'' if
the representation $\pi$ or both $\pi$ and $A$ are understood,
respectively.  Evidently, the usual and the flipped internal tensor
product are related by a unitary map $\Sigma \colon F \tr E \mycong
E \tl F$, $\eta \tr \xi \mapsto \xi \tl \eta$.

\section{Compact $C^{*}$-quantum graphs}
\label{section:graph}

The first basic ingredient in the definition of a compact
$C^{*}$-quantum groupoids are compact $C^{*}$-quantum
graphs. Roughly, such a compact $C^{*}$-quantum graph consists of a
$C^{*}$-algebra $B$ (of units) with a faithful KMS-state, a $C^{*}$-algebra
$A$ (of arrows), and two compatible module structures consisting of
representations $B,B^{(op)} \to A$ and  conditional expectations $A
\to B,B^{(op)}$. Thinking of (the underlying graph of) a groupoid,
the representations correspond to the range and the source map, and
the conditional expectations to the left and the right Haar weight.

Throughout the following sections, we will use several GNS- and
Rieffel-constructions for compact $C^{*}$-quantum graphs.  We first
recall the GNS-construction for KMS-states and present the
Rieffel-construction for a single module structure, before we turn
to compact $C^{*}$-quantum graphs.  To prepare for the
definition of the unitary antipode of a compact
$C^{*}$-quantum groupoid,  we finally discuss
coinvolutions on compact $C^{*}$-quantum graphs.

\paragraph{KMS-states on $C^{*}$-algebras and associated
  GNS-constructions}

We shall use the theory of KMS-states on $C^{*}$-algebras, see
\cite[\S 5]{bratteli:2}, \cite[\S 8.12]{pedersen}, and adopt the
following conventions. Let $\mu$ be a faithful KMS-state on a
$C^{*}$-algebra $B$. We denote by $\sigma^{\mu}$ the modular
automorphism group, by $H_{\mu}$ the GNS-space, by $\Lambda_{\mu}
\colon B \to H_{\mu}$ the GNS-map, by $\zeta_{\mu}
=\Lambda_{\mu}(1_{B})$ the cyclic vector, and by
$J_{\mu} \colon H_{\mu} \to H_{\mu}$ the modular conjugation
associated to $\mu$. We shall frequently use the formula
\begin{align} \label{eq:graph-modular-conjugation}
  J_{\mu}\Lambda_{\mu}(b) &=
  \Lambda_{\mu}(\sigma^{\mu}_{i/2}(b)^{*}) \quad \text{for all } b
  \in \Dom(\sigma^{\mu}_{i/2}).
\end{align}
We omit explicit mentioning of the GNS-representation $\pi_{\mu}
\colon B \to {\cal L}(H_{\mu})$ and identify $B$ with
$\pi_{\mu}(B)$; thus, $b\Lambda_{\mu}(x) = \pi_{\mu}(b)
\Lambda_{\mu}(x) = \Lambda_{\mu}(bx) = bx\zeta_{\mu}$ for all $b,x
\in B$.

We denote by $B^{op}$ the opposite $C^{*}$-algebra of $B$, which
coincides with $B$ as a Banach space with involution but has the
reversed multiplication, and by $\mu^{op} \colon B^{op} \to
\complex$ the opposite state of $\mu$, given by by
$\mu^{op}(b^{op}):=\mu(b)$ for all $b \in B$.  Using formula
\eqref{eq:graph-modular-conjugation}, one easily verifies that
$\mu^{op}$ is a KMS-state, that the modular automorphism group
$\sigma^{\mu^{op}}$ is given by
$\sigma^{\mu^{op}}_{t}(b^{op})=\sigma^{\mu}_{-t}(b)^{op}$ for all $b
\in B$, $t \in \reals$, and that one can always choose the GNS-space
and GNS-map for $\mu^{op}$ such that $H_{\mu^{op}}=H_{\mu}$ and
$\Lambda_{\mu^{op}}(b^{op})=J_{\mu}\Lambda_{\mu}(b^{*})$ for all $b
\in B$. Then $\zeta_{\mu^{op}}=\zeta_{\mu}$, $J_{\mu^{op}}=J_{\mu}$,
 $\pi_{\mu^{op}}(b)=J_{\mu}\pi_{\mu}(b)^{*}J_{\mu}$ for all $b
\in B$, and
\begin{align*}
  \Lambda_{\mu^{op}}(b^{op}) &=
  \Lambda_{\mu}(\sigma^{\mu}_{-i/2}(b)), &
  b^{op} \Lambda_{\mu}(x) &= \Lambda_{\mu}(x\sigma^{\mu}_{-i/2}(b))
\quad \text{for all } b \in \Dom(\sigma^{\mu}_{-i/2}), \, x \in B.
\end{align*}

\paragraph{Module structures and associated Rieffel constructions}
We shall use the following kind of module structures on
$C^{*}$-algebras relative to KMS-states:
\begin{definition} \label{definition:graph-module} Let $\mu$ be
  faithful KMS-state on a unital $C^{*}$-algebra $B$.  A {\em
    $\mu$-module structure} on a unital $C^{*}$-algebra $A$ is pair
  $(r,\phi)$ consisting of a unital embedding $r \colon B \to A$ and
  a completely positive map $\phi \colon A \to B$ such that $r \circ
  \phi \colon A \to r(B)$ is a unital conditional expectation,
  $\nu:=\mu \circ \phi$ is a KMS-state, and 
  $\sigma^{\nu}_{t}(r(B)) \subseteq r(B)$ for all $t \in \reals$.
\end{definition}
Given a module structure as above, we can form a
GNS-Rieffel-construction as follows:
\begin{lemma} \label{lemma:graph-rieffel} Let $\mu$ be a
  faithful KMS-state on a unital $C^{*}$-algebra $B$, let
  $(r,\phi)$ be a $\mu$-module structure on a unital
  $C^{*}$-algebra $A$, and put $\nu:=\mu \circ \phi$.
\begin{enumerate}
 \item $\sigma^{\nu}_{t} \circ r =r \circ \sigma^{\mu}_{t}$ for all
    $t \in \reals$.
  \item There exists a unique isometry $\zeta \colon H_{\mu}
    \hookrightarrow H_{\nu}$ such that $\zeta\Lambda_{\mu}(b) =
    \Lambda_{\nu}(r(b))$ for all $b \in B$.
  \item $\zeta J_{\mu}=J_{\nu} \zeta$, $\zeta b= r(b)\zeta$,
    $\zeta^{*}\Lambda_{\nu}(a)=\Lambda_{\mu}(\phi(a))$,
    $\zeta^{*}a = \phi(a)\zeta^{*}$ for all $b \in B$, $a \in
    A$.
  \item There exists a $\mu^{op}$-module structure
    $(r^{op},\phi^{op})$ on $A^{op}$ such that
    $r^{op}(b^{op})=r(b)^{op}$ and
    $\phi^{op}(a^{op})=\phi(a)^{op}$ for all $b \in B$, $a \in
    A$. For all $b \in B$, $\zeta \Lambda_{\mu^{op}}(b^{op}) =
    \Lambda_{\nu^{op}}(r^{op}(b^{op}))$.
    \end{enumerate}  
\end{lemma}
\begin{proof}
  i) This follows easily from the uniqueness of the modular
  automorphism group of a faithful  KMS-state.

  ii) Straightforward.

  iii) $\zeta J_{\mu} = J_{\nu} \zeta$ because
  $\Dom(\sigma^{\mu}_{i/2})$ is dense in $B$ and $\zeta J_{\mu}
  \Lambda_{\mu}(b) = \zeta \Lambda_{\mu}(\sigma^{\mu}_{i/2}(b)^{*}) =
  \Lambda_{\nu}(r(\sigma^{\mu}_{i/2}(b)^{*})) =
  \Lambda_{\nu}(\sigma^{\nu}_{i/2}(r(b))^{*}) =
  J_{\nu}\zeta\Lambda_{\mu}(b)$ for all $b \in \Dom(\sigma^{\mu}_{i/2})$
  by i). The proof of the remaining assertions is routine.

  iv) Straightforward.
\end{proof}

\paragraph{Compact $C^{*}$-quantum graphs} The definition of a
compact $C^{*}$-quantum graph involves the following simple variant
of a Radon-Nikodym derivative for KMS-states:
\begin{lemma} \label{lemma:graph-rn}
  Let $A$ be a unital $C^{*}$-algebra with a KMS-state $\nu$ and a
  positive invertible element $\delta$ that satisfies $\nu(\delta) = 1$
  and $\sigma^{\nu}_{t}(\delta)=\delta$ for all $t \in \reals$.
  \begin{enumerate}
  \item The state $\nu_{\delta}$ on $A$ given by $\nu_{\delta}(a) =
    \nu(\delta^{1/2}a\delta^{1/2})$ for all $a \in A$ is a faithful
    KMS-state and $\sigma_{t}^{\nu_{\delta}} = \Ad_{\delta^{it}} \circ
    \sigma_{t}^{\nu} = \sigma_{t}^{\nu} \circ \Ad_{\delta^{it}}$ for all
    $t \in \reals$.
  \item The map $\Lambda_{\nu_{\delta}} \colon A \to H_{\nu}$, $a
    \mapsto \Lambda_{\nu}(a\delta^{1/2})$, is a GNS-map for
    $\nu_{\delta}$, and the associated modular conjugation
    $J_{\nu_{\delta}}$ is equal to $J_{\nu}$.
  \item If $\tilde \delta \in A$ is another positive invertible element
    satisfying $\nu(\tilde \delta)=1$, $\sigma^{\nu}_{t}(\tilde
    \delta)=\delta$ for all $t \in \reals$, and $\nu_{\tilde
      \delta}=\nu_{\delta}$, then $\delta=\tilde \delta$.  \qed
  \end{enumerate}
\end{lemma}
\begin{definition} \label{definition:graph-graph} 
A {\em compact
 $C^{*}$-quantum graph} is a tuple $(B,\mu,A,r,\phi,s,\psi,\delta)$, where
  \begin{enumerate}
  \item $B$ is a unital $C^{*}$-algebra with a faithful KMS-state
    $\mu$,
  \item $A$ is a unital $C^{*}$-algebra,
  \item $(r,\phi)$ and $(s,\psi)$ are a $\mu$-module
    structure and a $\mu^{op}$-module structure on $A$,
    respectively, such that $r(B)$ and $s(B^{op})$ commute,
  \item $\delta \in A \cap r(B)' \cap s(B^{op})'$ is a positive,
    invertible, $\sigma^{\nu}$-invariant element such that
    $\nu(\delta)=1$ and $\mu^{op} \circ \psi = (\mu \circ
    \phi)_{\delta}$.
  \end{enumerate}
  Given a compact  $C^{*}$-quantum graph $(B,\mu,A,r,\phi,s,\psi)$, we
  put $\nu:=\mu \circ \phi$, $\nu^{-1}:=\mu^{op} \circ \psi$, and
  denote by $\zeta_{\phi}, \zeta_{\psi} \colon H_{\mu} \to H_{\nu}$ the
  isometries associated to $(r,\phi)$, $(s,\psi)$ as in Lemma
  \ref{lemma:graph-rieffel}.
\end{definition}
For every compact $C^{*}$-quantum graph $(B,\mu,A,r,\phi,s,\psi)$,
we have $\nu \circ r = \mu \circ \phi \circ r = \mu$ and $ \nu^{-1}
\circ s = \mu^{op} \circ \psi \circ s = \mu^{op}$. The compositions
$\nu \circ s$ and $\nu^{-1} \circ r$ are related to $\mu^{op}$ and
$\mu$, respectively, as follows.
\begin{lemma} \label{lemma:graph} Let $(B,\mu,A,r,\phi,s,\psi,\delta)$
  be a compact $C^{*}$-quantum graph.
  \begin{enumerate}
  \item $\phi(\delta) \in B$ and $\psi(\delta^{-1}) \in B^{op}$ are
    positive, invertible, central, invariant with respect to
    $\sigma^{\mu}$ and $\sigma^{\mu^{op}}$, respectively, and
    satisfy $\mu(\phi(\delta))=1=\mu^{op}(\psi(\delta^{-1}))$.
  \item $\nu^{-1} \circ r = \mu_{\phi(\delta)}$ and $\nu \circ s =
    \mu^{op}_{\psi(\delta^{-1})}$.
  \end{enumerate}
\end{lemma}
\begin{proof}
  i) We only prove the assertions concerning $\phi(\delta)$.  Since
  $\delta$ is positive, there exists an $\epsilon > 0$ such that $\delta
  > \epsilon 1_{A}$, and since $\phi$ is positive, we can conclude
  $\phi(\delta) > \epsilon \phi(1_{A})=\epsilon 1_{B}$. Therefore,
  $\phi(\delta)$ is positive and invertible.  It is central because $b
  \phi(\delta) = \phi(r(b)\delta) = \phi(\delta r(b))=\phi(\delta)b$ for
  all $b \in B$, and invariant under $\sigma^{\mu}$ because 
  $\sigma^{\mu}_{t}(\phi(\delta)) = \phi(\sigma^{\nu}_{t}(\delta)) =
\phi(\delta)$ for all $t \in \reals$.

ii) The first equation holds because for all $b \in B$,
$\nu^{-1}(r(b)) = \mu(\phi(\delta^{1/2}r(b)\delta^{1/2})) = \mu(b
\phi(\delta)) = \mu(\phi(\delta)^{1/2}b\phi(\delta)^{1/2})$.  The
second equation follows similarly.
\end{proof}

Let $(B,\mu,A,r,\phi,s,\psi,\delta)$ be a compact $C^{*}$-quantum
graph. Then for all $b,c \in B$, $\psi(r(b))c^{op} =
\psi(r(b)s(c^{op})) = \psi(s(c^{op})r(b)) = c^{op}\psi(r(b))$ and
similarly $\phi(s(b^{op}))c=c\phi(s(b^{op}))$, so that we can define
completely positive maps
\begin{align} \label{eq:graph-tau}
  \tau &:=\psi \circ r \colon B \to Z(B^{op}) & & \text{and} &
  \tau^{\dag} &:=\phi \circ s \colon B^{op} \to Z(B).
\end{align}
We identify $Z(B)$ and $Z(B^{op})$ with $B \cap B^{op} \subseteq
{\cal L}(H_{\mu})$ in the natural way. 
 \paragraph{Coinvolutions on compact $C^{*}$-quantum graphs}
 The following concept will be used to define the unitary antipode
 of a compact $C^{*}$-quantum groupoid:
 \begin{definition} \label{definition:graph} A {\em coinvolution} for a
   compact $C^{*}$-quantum graph $(B,\mu,A,r,\phi,s,\psi,\delta)$ is an
   antiautomorphism $R \colon A \to A$ satisfying $R \circ R = \Id_{A}$
   and $R(r(b)) = s(b^{op})$, $\phi(R(a)) = \psi(a)^{op}$ for all $b \in
   B$, $a \in A$.
\end{definition}
\begin{lemma} \label{lemma:graph-coinvolution} Let $R$ be a
  coinvolution for a compact $C^{*}$-quantum graph
  $(B,\mu,A,r,\phi,s,\psi,\delta)$.
  \begin{enumerate}
  \item $\sigma^{\nu}_{t} \circ R = R \circ \sigma^{\nu^{-1}}_{-t}$ for
    all $t \in \reals$, and $R(\delta)=\delta^{-1}$. In particular,
    $\phi(\delta)=\psi(\delta^{-1})$.
  \item $\tau(b)=\tau^{\dag}(b^{op})$ for all $b \in B$.
  \item There exists a unique antiunitary $I \colon H_{\nu} \to
    H_{\nu}$, $\Lambda_{\nu^{-1}}(a) \mapsto \Lambda_{\nu}(R(a)^{*})$,
    and
    \begin{align*} 
      I\Lambda_{\nu}(a) &= \Lambda_{\nu}(R(a\delta^{1/2})^{*}), &
      Ia^{*}I &= R(a), & I^{2}&= \Id_{H}, & I \zeta_{\psi} J_{\mu} &=
      \zeta_{\phi}, & IJ_{\nu} &= J_{\nu}I.
    \end{align*}
  \end{enumerate}
  \end{lemma}
  \begin{proof}
    i) The first equation follows from the fact that $R$ is an
    antiautomorphism and that $\nu \circ R = \nu^{-1}$. To prove the
    second equation, put $\delta':=R(\delta^{-1})$. Then
    \begin{gather*}
      \nu(\delta') = \nu^{-1}(\delta^{-1}) = \nu(1) = 1, \\
      \sigma^{\nu}_{t}(\delta') =
      R(\sigma^{\nu^{-1}}_{-t}(\delta^{-1})) =
      R\big(\sigma^{\nu}_{-t} \circ
      \Ad_{\delta^{it}}(\delta^{-1})\big) =
      R(\sigma^{\nu}_{-t}(\delta^{-1})) = R(\delta^{-1}) =
      \delta', \\
      \nu_{\delta}(a) = \nu^{-1}(a) = \nu(R(a)) =
      \nu^{-1}(\delta^{-1/2}R(a)\delta^{-1/2}) =
      \nu(\delta'{}^{1/2}a\delta'{}^{1/2}) = \nu_{\delta'}(a)
    \end{gather*}
    for all $a \in A$, and by Lemma \ref{lemma:graph-rn} iii),
    $\delta=\delta'$.

    ii) $(\phi \circ s)(b^{op}) = (\phi\circ R \circ R \circ
    s)(b^{op}) = (\psi \circ r)(b)^{op}$ for all $b \in B$.

    iii) The formula for $I$ defines an antiunitary because for all $a
    \in A$,
    \begin{align*}
      \big\|\Lambda_{\nu}(R(a)^{*})\big\|^{2} &= (\mu \circ
      \phi)\big(R(a)R(a)^{*}\big) = (\mu \circ \phi \circ R)(a^{*}a)
      = (\mu^{op} \circ \psi)(a^{*}a) =
      \big\|\Lambda_{\nu^{-1}}(a)\big\|^{2}.
    \end{align*}
    The first two equations given in ii)
    follow immediately. Next, $I^{2}=\Id_{H}$ because 
    \begin{align*}
      I^{2}\Lambda_{\nu}(a) &=
      \Lambda_{\nu}\big(R(R(a\delta^{1/2})^{*}\delta^{1/2})^{*}\big)
      = \Lambda_{\nu}\big(a\delta^{1/2}\delta^{-1/2}\big) =
      \Lambda_{\nu}(a)
    \end{align*}
    for all $a \in A$, and $I\zeta_{\psi} J_{\mu} = \zeta_{\phi}$
    because
  \begin{align*}
    I\zeta_{\psi} J_{\mu}\Lambda_{\mu}(b^{*}) =
    I\zeta_{\psi}\Lambda_{\mu^{op}}(b^{op}) =
    I\Lambda_{\nu^{-1}}(s(b^{op})) = \Lambda_{\nu}(r(b)^{*})
    = \zeta_{\phi}\Lambda_{\mu}(b^{*})
  \end{align*}
  for all $b \in B$.  The relation $\sigma^{\nu^{-1}}_{t} = R \circ
  \sigma^{\nu}_{-t} \circ R$ ($t \in \reals$) implies that for all
  $a \in \Dom(\sigma^{\nu^{-1}}_{i/2})$,
    \begin{align*}
      J_{\nu}I \Lambda_{\nu^{-1}}(a)  =
      \Lambda_{\nu}\big(\sigma^{\nu}_{i/2}(R(a)^{*})^{*}\big) &=
      \Lambda_{\nu}\big(R(\sigma^{\nu^{-1}}_{i/2}(a)^{*})^{*}\big) 
      = IJ_{\nu^{-1}}\Lambda_{\nu^{-1}}(a).
    \end{align*}
    Since $J_{\nu^{-1}}=J_{\nu}$, we can conclude
    $J_{\nu}I=IJ_{\nu}$.
\end{proof}

\section{The relative tensor product and the fiber product}

\label{section:modules}

Fundamental to the following development is the general language of
$C^{*}$-modules and $C^{*}$-algebras over KMS-states, the relative
tensor product of such $C^{*}$-modules, and the fiber product of
such $C^{*}$-algebras: The fiber product is
needed to define the target of the comultiplication of a compact
$C^{*}$-quantum groupoid, and the relative tensor product is needed
to define this fiber product and the domain and the range of the
fundamental unitary.

We proceed as follows. First, we introduce the language of
$C^{*}$-modules and $C^{*}$-algebras over KMS-states.  Next, we
describe the $C^{*}$-module structures that arise from a compact
$C^{*}$-quantum graph and which are needed later. Finally, we present the
relative tensor product and the fiber product. Except for the second
paragraph, the  reference is \cite{timmer:cpmu}

\paragraph{$C^{*}$-modules and $C^{*}$-algebras over KMS-states}
We adapt the framework of $C^{*}$-modules and $C^{*}$-algebras over
 $C^{*}$-bases introduced in \cite{timmer:cpmu} to our
present needs, replacing $C^{*}$-bases by KMS-states as follows.  A
{\em $C^{*}$-base} is a triple $\cbasel{B}{H}$ consisting of a
Hilbert space $\frakH$ and two commuting nondegenerate
$C^{*}$-algebras $\frakB,\frakB^{\dag} \subseteq {\cal L}(\frakH)$.
We restrict ourselves to $C^{*}$-bases of the form $(H_{\mu},B,B^{op})$,
where $H_{\mu}$ is the GNS-space of a faithful KMS-state $\mu$ on a
$C^{*}$-algebra $B$, and where $B$ and $B^{op}$ act on
$H_{\mu}=H_{\mu^{op}}$ via the GNS-representations, and obtain the
following notions of $C^{*}$-modules over
$\mu$.
\begin{definition}
  Let $\mu$ be a faithful KMS-state on a $C^{*}$-algebra $B$.  A
  {\em $C^{*}$-$\mu$-module} is a pair $(H,\alpha)$, briefly written
  $H_{\alpha}$, where $H$ is a Hilbert space and $\alpha \subseteq
  {\cal L}(H_{\mu},H)$ is a closed subspace satisfying $[\alpha
  H_{\mu}]=H$, $[\alpha B]= \alpha$, and $[\alpha^{*}\alpha] = B
  \subseteq {\cal L}(H_{\mu})$. 
  A {\em morphism} between $C^{*}$-$\mu$-modules $H_{\alpha}$ and
  $K_{\beta}$ is an operator $T \in {\cal L}(H,K)$ satisfying
  $T\alpha \subseteq \beta$ and $T^{*}\beta \subseteq \alpha$. We
  denote the set of such morphisms by ${\cal
    L}(H_{\alpha},K_{\beta})$.
\end{definition}
\begin{lemma} \label{lemma:modules}
    Let $\mu$ be a faithful KMS-state on a $C^{*}$-algebra $B$ and
    let $H_{\alpha}$ be a $C^{*}$-$\mu$-module.
    \begin{enumerate}
    \item $\alpha$ is a right Hilbert $C^{*}$-$B$-module with inner
      product given by $\langle \xi|\xi'\rangle=\xi^{*}\xi$ for all
      $\xi,\xi' \in \alpha$.
    \item There exist  isomorphisms $
      \alpha \tr H_{\mu} \to H$, $\xi \tr \zeta \mapsto \xi \zeta$,
      and $H_{\mu} \tl \alpha \to H$, $\zeta \tl
      \xi \mapsto \xi\zeta$.
    \item There exists a nondegenerate  representation $\rho_{\alpha} \colon
      B^{op} \to {\cal L}(H)$ such that $\rho_{\alpha}(b^{op})(\xi
      \zeta)= \xi b^{op} \zeta$ for all $b \in B$, $\xi \in \alpha$,
      $\zeta \in H_{\mu}$.
    \item Let $K_{\beta}$ be a $C^{*}$-$\mu$-module and $T \in {\cal
        L}(H_{\alpha},K_{\beta})$. Then left multiplication by $T$
      defines an operator in ${\cal
        L}_{B}(\alpha,\beta)$, again denoted by $T$, and
      $T\rho_{\alpha}(b^{op}) = \rho_{\beta}(b^{op})T$ for all $b
      \in B$.  \qed
    \end{enumerate}
\end{lemma}
\begin{definition}
  Let $\mu_{1},\ldots,\mu_{n}$ be faithful KMS-states on
  $C^{*}$-algebras $B_{1}, \ldots, B_{n}$. A {\em
    $C^{*}$-$(\mu_{1},\ldots,\mu_{n})$-module} is a tuple
  $(H,\alpha_{1},\ldots,\alpha_{n})$, where $H$ is a Hilbert space
  and $(H,\alpha_{i})$ is a $C^{*}$-$\mu_{i}$-module for each
  $i=1,\ldots,n$ such that
  $[\rho_{\alpha_{i}}(B_{i}^{op})\alpha_{j}]=\alpha_{j}$ whenever
  $i\neq j$. The set of morphisms of
  $C^{*}$-$(\mu_{1},\ldots,\mu_{n})$-modules
  $(H,\alpha_{1},\ldots,\alpha_{n})$ and
  $(K,\beta_{1},\ldots,\beta_{n})$ is ${\cal
    L}((H,\alpha_{1},\ldots,\alpha_{n}),
  (K,\beta_{1},\ldots,\beta_{n})) := \bigcap_{i=1}^{n} {\cal
    L}(H_{\alpha_{i}},K_{\beta_{i}}) \subseteq {\cal L}(H,K)$.
\end{definition}
\begin{remark}
  Let $\mu_{1},\ldots,\mu_{n}$ be faithful KMS-states on
  $C^{*}$-algebras $B_{1}, \ldots, B_{n}$ and let
  $(H,\alpha_{1},\ldots,\alpha_{n})$ be a
  $C^{*}$-$(\mu_{1},\ldots,\mu_{n})$-module. Then
  $\rho_{\alpha_{i}}(B^{op}_{i}) \subseteq {\cal L}(H_{\alpha_{j}})$
  whenever $i\neq j$; in particular,
  $[\rho_{\alpha_{i}}(B^{op}_{i}),\rho_{\alpha_{j}}(B^{op}_{j})]=0$
  whenever $i \neq j$.
\end{remark}

Next, we define $C^{*}$-algebras over KMS-states. 
\begin{definition}
  Let $\mu_{1},\ldots,\mu_{n}$ be faithful KMS-states on
  $C^{*}$-algebras $B_{1},\ldots,B_{n}$.  A {\em
    $C^{*}$-$(\mu_{1},\ldots,\mu_{n})$-algebra} consists of a  $C^{*}$-$(\mu_{1},\ldots,\mu_{n})$-module
  $(H,\alpha_{1},\ldots,\alpha_{n})$ and a nondegenerate
  $C^{*}$-algebra $A \subseteq {\cal L}(H)$ such that
  $[\rho_{\alpha_{i}}(B_{i}^{op})A] \subseteq A$ for each
  $i=1,\ldots,n$. In the cases $n=1,2$, we abbreviate
  $A^{\alpha}_{H}:=(H_{\alpha},A)$,
  $A^{\alpha,\beta}_{H}:=((H,\alpha,\beta),A)$.
  A {\em morphism} between
  $C^{*}$-$(\mu_{1},\ldots,\mu_{n})$-algebras
  $((H,\alpha_{1},\ldots,\alpha_{n}),A)$ and
  $((K,\gamma_{1},\ldots,\gamma_{n}),C)$ is a nondegenerate
  $*$-homomorphism $\phi \colon A \to M(C)$ such that for each
  $i=1,\ldots,n$, we have $[I_{\phi,i}\alpha_{i}] = \gamma_{i}$,
  where $I_{\phi,i}:=\{T \in {\cal L}(H_{\alpha_{i}},
  K_{\gamma_{i}}) \mid Ta = \phi(a)T$ for all $a \in A\}$. We denote
  the set of all such morphisms by
  $\Mor(((H,\alpha_{1},\ldots,\alpha_{n}),A),
  ((K,\gamma_{1},\ldots,\gamma_{n}),C))$.
\end{definition}
\begin{remark}
  If $\phi$ is a morphism between $C^{*}$-$\mu$-algebras $A_{H}^{\alpha}$ and
  $C_{K}^{\gamma}$, then $\phi(\rho_{\alpha}(b^{op}))=\rho_{\gamma}(b^{op})$
  for all $b \in B$; see \cite[Lemma 2.2]{timmer:ckac}.
\end{remark}

\paragraph{The $C^{*}$-module of a compact $C^{*}$-quantum graph} 
To proceed from compact $C^{*}$-quantum graphs to compact
$C^{*}$-quantum groupoids, we need several $C^{*}$-module structures
arising from the GNS-Rieffel-construction in Lemma
\ref{lemma:graph-rieffel}. 
\begin{lemma} \label{lemma:module} Let $\mu$ be a faithful
  KMS-state on a unital $C^{*}$-algebra $B$, let $(r,\phi)$ be a
  $C^{*}$-$\mu$-module structure on a unital $C^{*}$-algebra $A$,
  and put $\nu:=\mu \circ \phi$,  $H:=H_{\nu}$, $\halpha:=[A\zeta]$, 
  $\beta:=[A^{op}\zeta]$.
  \begin{enumerate}
  \item ${_{\halpha}H_{\beta}}$ is a
    $C^{*}$-$(\mu,\mu^{op})$-module and $\rho_{\halpha} = r^{op}$,
 $ \rho_{\beta} = r$.
  \item $A^{\beta}_{H}$ is a $C^{*}$-$\mu^{op}$-algebra.
  \item $a^{op} \zeta = \sigma^{\nu}_{-i/2}(a)\zeta$ 
    for all $a \in \Dom(\sigma^{\nu}_{-i/2}) \cap
    r(B)'$.
  \item $A + (A \cap r(B)')^{op} \subseteq {\cal L}(H_{\halpha})$
    and $A^{op} + (A \cap r(B)') \subseteq {\cal L}(H_{\beta})$.
  \end{enumerate}
   \end{lemma}
   \begin{proof}
     i) Lemma \ref{lemma:graph-rieffel} immediately implies that
     $H_{\halpha}$ is a $C^{*}$-$\mu$-module and that $H_{\beta}$ is a
     $C^{*}$-$\mu^{op}$-module.  The equations for $\rho_{\halpha}$
     and $\rho_{\beta}$ follow from the fact that by Lemma
     \ref{lemma:graph-rieffel}, $\rho_{\halpha}(b^{op}) a\zeta =
     a\zeta b^{op} = a r(b)^{op} \zeta = r(b)^{op} a\zeta$ and
     $\rho_{\beta}(b) a^{op} \zeta = a^{op}\zeta b = a^{op}r(b)\zeta
     = r(b) a^{op}\zeta$ for all $b \in B$, $a \in A$.  In
     particular, $[\rho_{\halpha}(B^{op})\beta] =
     [r(B)^{op}A^{op}\zeta] =\beta$ and
     $[\rho_{\beta}(B)\halpha]=[r(B)A\zeta]=\halpha$, whence
     ${_{\halpha} H_{\beta}}$ is a $C^{*}$-$(\mu,\mu^{op})$-module.

     ii) By i), $[\rho_{\beta}(B)A] =A$.

     iii) For all $a \in  r(B)' \cap \Dom(\sigma^{\nu}_{-i/2})$ and
     $b\in B$,
     \begin{align*}
       a^{op}   \zeta \Lambda_{\mu}(b) &=
    \Lambda_{\nu}\big(r(b)\sigma^{\nu}_{-i/2}(a)\big) 
    =\Lambda_{\nu}\big(\sigma^{\nu}_{-i/2}(a)r(b)\big)
    = \sigma^{\nu}_{-i/2}(a) \zeta \Lambda_{\mu}(b).
  \end{align*}

  iii) We only prove the first inclusion, the second one follows
  similarly. Clearly, $[A\halpha]=\halpha$. Since $\sigma_{t}^{\nu}
  (r(B)) \subseteq r(B)$ for all $t \in \reals$, the subspace $C:=
  \Dom(\sigma^{\nu}_{-i/2}) \cap r(B)'$ is dense in $A \cap r(B)'$,
  and by iii), $[(A \cap r(B)')^{op}\halpha] = [CA\zeta] = [AC\zeta]
  \subseteq [A\zeta]=\halpha$.
   \end{proof}
\begin{proposition} \label{proposition:module-graph}
Let $(B,\mu,A,r,\phi,s,\psi,\delta)$ be a compact
$C^{*}$-quantum graph. Put $\nu:=\mu \circ \phi$, $\nu^{-1}:=\mu^{op}
\circ \psi = \nu_{\delta}$ and
\begin{align}
  H &:= H_{\nu}, & \halpha &:=[A\zeta_{\phi}], & \beta
  &:=[A^{op}\zeta_{\phi}], & \hbeta &:=[A\zeta_{\psi}], & \alpha &:=
  [A^{op}\zeta_{\psi}]. 
\end{align}
\begin{enumerate}
\item $(H,\halpha,\beta,\hbeta,\alpha)$ is a
  $C^{*}$-$(\mu,\mu^{op},\mu^{op},\mu)$-module and $\rho_{\halpha}=
  r^{op}$, $\rho_{\beta}=r$, $\rho_{\hbeta}=s^{op}$,
  $\rho_{\alpha}=s$.
\item $A_{H}^{\alpha,\beta}$ is a $C^{*}$-$(\mu,\mu^{op})$-algebra.
\item Let $R$ be a coinvolution for $(B,\mu,A,r,\phi,s,\psi,\delta)$ and
  let $I \colon H_{\nu} \to H_{\nu}$ be given by $\Lambda_{\nu^{-1}}(a)
  \mapsto \Lambda_{\nu}(R(a)^{*})$. Then $I\zeta_{\phi}J_{\mu} =
  \zeta_{\psi}$, $I \zeta_{\psi}J_{\mu} = \zeta_{\phi}$ and $I\hbeta
  J_{\mu} =\halpha$, $I\beta J_{\mu} = \alpha$.
\end{enumerate}
\end{proposition}
\begin{proof}
  i), ii)  Immediate from Lemma \ref{lemma:module}.

  iii)  We have $I \zeta_{\psi} J_{\mu} = \zeta_{\phi}$ because
  for all $b \in B$,
  \begin{align*}
    I\zeta_{\psi}J_{\mu}\Lambda_{\mu}(b^{*})  = I \zeta_{\psi}
    \Lambda_{\mu^{op}}(b^{op}) = I \Lambda_{\nu^{-1}}(s(b^{op}))
    = \Lambda_{\nu}(R(s(b^{op}))^{*}) =
    \Lambda_{\nu}(r(b^{*})) = \zeta_{\phi}\Lambda_{\mu}(b^{*}).
  \end{align*}
  The remaining assertions follow easily.
\end{proof}

\paragraph{The relative tensor product of $C^{*}$-modules}
The relative tensor product of $C^{*}$-modules over KMS-states is a
symmetrized version of the internal tensor product of Hilbert
$C^{*}$-modules and a $C^{*}$-algebraic analogue of the relative
tensor product of Hilbert spaces over a von Neumann algebra. We
briefly summarize the definition and the main properties; for
details, see \cite[Section 2.2]{timmer:cpmu}.

Let $\mu$ be a faithful KMS-state on a $C^{*}$-algebra $B$, let 
$H_{\beta}$ be a $C^{*}$-$\mu$-module, and let $K_{\gamma}$ be a
$C^{*}$-$\mu^{op}$-module. The {\em relative tensor product} of
$H_{\beta}$ and $K_{\gamma}$ is the Hilbert space $H
\stensor{\beta}{\gamma} K:=\beta \tr H_{\mu} \tl \gamma$.   
It is spanned by elements $\xi \tr \zeta \tl \eta$, where $\xi \in
\beta$, $\zeta \in H_{\mu}$, $\eta \in \gamma$, and the inner product is
given by $\langle \xi \tr \zeta \tl \eta|\xi' \tr \zeta' \tl
\eta'\rangle = \langle \zeta | \xi^{*}\xi' \eta^{*}\eta'
\zeta'\rangle = \langle \zeta|\eta^{*}\eta' \xi^{*}\xi'
\zeta'\rangle$ for all $\xi,\xi' \in \beta$, $\zeta,\zeta' \in
H_{\mu}$, $\eta,\eta' \in \gamma$.

Obviously, there exists {\em flip} isomorphism
\begin{align} \label{eq:modules-flip-modules}
  \Sigma \colon H \stensor{\beta}{\gamma} K \to K
  \stensor{\gamma}{\beta} H, \quad \xi \tr \zeta \tl \eta \mapsto
  \eta \tr \zeta \tl \xi.
\end{align}

The isomorphisms $\beta \tr H_{\mu} \cong H$, $\xi \tr \zeta \equiv
\xi \zeta$, and $H_{\mu} \tl \gamma \cong K$, $\zeta \tl \zeta \equiv
\eta \zeta$, (see Lemma \ref{lemma:modules}) induce the following
isomorphisms which we use without further notice:
\begin{gather*}
  H {_{\rho_{\beta}} \tl} \gamma \cong H \stensor{\beta}{\gamma} K
  \cong \beta \tr_{\rho_{\gamma}} K, \qquad \xi \zeta \tl \eta
  \equiv \xi \tr \zeta \tl \eta \equiv \xi \tr \eta\zeta \quad (\xi
  \in \beta,\, \zeta \in H_{\mu},\, \eta \in \gamma).
\end{gather*}

Using these isomorphisms, we define the following tensor products of
operators:
\begin{align*}
  S \stensor{\beta}{\gamma}T:= S \tr T &\in {\cal L}(\beta
  \tr_{\rho_{\gamma}} K) = {\cal L}(H \stensor{\beta}{\gamma} K)
  \quad \text{ for all } S \in {\cal L}(H_{\alpha}), \ T \in
  \rho_{\gamma}(B)' \subseteq {\cal
    L}(K), \\
  S \stensor{\beta}{\gamma}T:= S \tl T &\in {\cal L}(H
  {_{\rho_{\beta}}\tl} \gamma) = {\cal L}(H \stensor{\beta}{\gamma}
  K) \quad \text{ for all } S \in \rho_{\beta}(B^{op})' \subseteq
  {\cal L}(H), \ T \in {\cal L}(K_{\beta}).
\end{align*}
Note that $S \tr T = S \tr \Id \tl T = S \tl
T$ for all $S \in {\cal L}(H_{\beta})$, $T \in {\cal
  L}(K_{\gamma})$.

For each $\xi \in \beta$, $\eta \in \gamma$, there exist bounded
linear operators
\begin{align*}
  |\xi\rangle_{\leg{1}} \colon K &\to H \stensor{\beta}{\gamma} K, \ \omega \mapsto \xi
  \tr \omega, & \langle \xi|_{\leg{1}}:=|\xi\rangle_{\leg{1}}^{*}\colon
  \xi' \tr \omega &\mapsto
  \rho_{\gamma}(\xi^{*}\xi')\omega, \\
  |\eta\rangle_{\leg{2}} \colon H &\to H \stensor{\beta}{\gamma} K, \ \omega \mapsto \omega
  \tl \eta, & \langle\eta|_{\leg{2}} := |\eta\rangle_{\leg{2}}^{*}
  \colon \omega \tl\eta &\mapsto \rho_{\beta}(
  \eta^{*}\eta')\omega.
\end{align*}
We put $\kbeta{1} := \big\{ |\xi\rangle_{\leg{1}} \,\big|\, \xi \in
\beta\big\}$ and similarly define $\bbeta{1}$, $\kgamma{2}$,
$\bgamma{2}$. 

Assume that $\frakH=(H,\alpha_{1},\ldots,\alpha_{m},\beta)$ is a
$C^{*}$-$(\sigma_{1},\ldots,\sigma_{m},\mu)$-module and that
$\frakK=(K,\gamma,\delta_{1},\ldots,\delta_{n})$ is a
$C^{*}$-$(\mu^{op},\tau_{1},\ldots,\tau_{n})$-module, where
$\sigma_{1},\ldots,\sigma_{m}, \tau_{1},\ldots, \tau_{n}$ are
faithful KMS-states on $C^{*}$-algebras
$A_{1},\ldots,A_{m},C_{1},\ldots,C_{n}$. For $i=1,\ldots,m$ and
$j=1,\ldots, n$,  put
\begin{align*}
  \alpha_{i} \lt \gamma &:= [\kgamma{2}\alpha_{i}] \subseteq {\cal
    L}(H_{\sigma_{i}}, H \stensor{\beta}{\gamma} K), &
  \beta \rt \delta_{j} &:= [\kbeta{1}\delta_{j}] \subseteq {\cal
    L}(H_{\tau_{j}}, H \stensor{\beta}{\gamma} K).
\end{align*}
Then the tuple $ \frakH \mutimes \frakK:= (H \stensor{\beta}{\gamma}
K, \alpha_{1} \lt \gamma, \ldots, \alpha_{m} \lt \gamma, \beta \rt
\delta_{1},\ldots, \beta \rt \delta_{n})$ is a
$C^{*}$-$(\sigma_{1},\ldots,\sigma_{m}$, $
\tau_{1},\ldots,\tau_{n})$-module, called the {\em relative tensor
  product} of $\frakH$ and $\frakK$. For all $i=1,\ldots,m$, $a \in
A_{i}$ and $j=1,\ldots,n$, $c \in C_{j}$,
\begin{align*}
  \rho_{(\alpha_{i} \lt \gamma)}(a^{op}) &=
  \rho_{\alpha_{i}}(a^{op}) \stensor{\beta}{\gamma} \Id, &
  \rho_{(\beta \rt \delta_{j})}(c^{op}) &= \Id
  \stensor{\beta}{\gamma} \rho_{\delta_{j}}(c^{op}).
\end{align*}

The $C^{*}$-relative tensor product is bifunctorial: If
$\tilde \frakH=(\tilde H,\tilde \alpha_{1},\ldots,\tilde
\alpha_{m},\tilde \beta)$ is a
$C^{*}$-$(\sigma_{1},\ldots,\sigma_{m},\mu)$-module, $\tilde
\frakK=(\tilde K,\tilde \gamma,\tilde \delta_{1},\ldots,\tilde
\delta_{n})$  a
$C^{*}$-$(\mu^{op},\tau_{1},\ldots,\tau_{n})$-module, and $S \in
{\cal L}(\frakH,\tilde \frakH)$, $T \in {\cal L}(\frakK,\tilde
\frakK)$, then there exists a unique operator $S
\mutimes T \in {\cal L}(\frakH \mutimes \frakK, \tilde
\frakH \mutimes \tilde \frakK)$ such that
\begin{align*}
  (S \mutimes T)(\xi \tr \zeta \tl \eta) = S\xi  \tr \zeta
   \tl T\eta \quad \text{for all } \xi \in \beta, \, \zeta \in
  H_{\mu}, \, \eta \in \gamma.
\end{align*}

The $C^{*}$-relative tensor product is unital in the following
sense. If we consider $B,B^{op}$ embedded in ${\cal L}(H_{\mu})$ via
the GNS-representations, then the tuple
$\mathfrak{U}:=(H_{\mu},B,B^{op})$ is a
$C^{*}$-$(\mu,\mu^{op})$-module, and the maps
\begin{align*}
  H \stensor{\beta}{B^{op}} H_{\mu} \to H, \ \xi  \tr \zeta \tl 
  b^{op} \mapsto \xi b^{op}\zeta, \quad
  H_{\mu} \, \stensor{B}{\gamma} K \to K, \ b \tr \zeta \tl \eta
  \mapsto \eta b \zeta,
\end{align*}
are isomorphisms of $C^{*}$-$(\sigma_{1},\ldots,\sigma_{m},\mu)$-
and $C^{*}$-$(\mu^{op},\tau_{1},\ldots,\tau_{n})$-modules $\frakH
\mutimes \mathfrak{U} \cong \frakH$ and $\mathfrak{U} \mutimes \frakK
\cong \frakK$, natural in $\frakH$ and $\frakK$, respectively.

The $C^{*}$-relative tensor product is associative in the following
sense. Assume that $\nu,\rho_{1},\ldots,\rho_{l}$ are faithful
KMS-states on $C^{*}$-algebras $D$,$E_{1},\ldots,E_{l}$, that $\hat
\frakK =(K,\gamma,\delta_{1},\ldots,\delta_{n},\epsilon)$ is a
$C^{*}$-$(\mu^{op},\tau_{1},\ldots,\tau_{n},\nu)$-module, and
$\frakL=(L,\phi,\psi_{1},\ldots,\psi_{l})$ a
$C^{*}$-$(\nu^{op},\rho_{1},\ldots,\rho_{l})$-module. Then the
isomorphisms of Hilbert spaces
\begin{align} \label{eq:modules-rtp-associative}
  (H \stensor{\beta}{\gamma} K) \stensor{\beta \rt \epsilon}{\phi} L
  \cong  \beta \tr_{\rho_{\gamma}} K {_{\rho_{\epsilon}}\tl} \, \phi
  \cong 
  H \stensor{\beta}{\gamma \lt \phi} (K \stensor{\epsilon}{\phi} L)
\end{align}
 are isomorphisms of
$C^{*}$-$(\sigma_{1},\ldots,\sigma_{m},\tau_{1},\ldots,\tau_{n},\rho_{1},\ldots,\rho_{l})$-modules
$(\frakH \mutimes \hat \frakK) \nutimes \frakL \cong \frakH \mutimes
(\hat \frakK \nutimes \frakL)$. We shall identify the Hilbert spaces in
\eqref{eq:modules-rtp-associative} without further notice and denote
these Hilbert spaces by $H \stensor{\beta}{\gamma} K
\stensor{\epsilon}{\phi} L$.

We shall need the following simple construction not mentioned in
\cite{timmer:cpmu}:
\begin{lemma} \label{lemma:modules-antiunitary} Let $H_{\beta}$,
  $\tilde H_{\tilde \beta}$ be $C^{*}$-$\mu$-modules, $K_{\gamma}$,
  $\tilde K_{\tilde \gamma}$ $C^{*}$-$\mu^{op}$-modules, and  $I
  \colon H \to \tilde H$, $J \colon K \to \tilde K$ 
  anti-unitaries such that $I\beta J_{\mu} = \tilde \beta$ and $J
  \gamma J_{\mu} = \tilde \gamma$.
  \begin{enumerate}
  \item There exists a unique anti-unitary $I
    \rtensor{\beta}{J_{\mu}}{\gamma} J \colon H \stensor{\beta}{\gamma}
    K \to \tilde H \stensor{\tilde \beta}{\tilde \gamma} \tilde K$ such
    that
    \begin{align*}
      (I \rtensor{\beta}{J_{\mu}}{\gamma} J)(\xi \tr \zeta \tl \eta) &=
      I\xi J_{\mu} \tr J_{\mu}\zeta \tl J\eta J_{\mu} \quad \text{for
        all } \xi \in \beta, \, \zeta \in H_{\mu}, \, \eta \in \gamma.
    \end{align*}
  \item $(I \rtensor{\beta}{J_{\mu}}{\gamma} J)|\xi\rangle_{1} = |I\xi
    J_{\mu}\rangle_{1} J$ and $(I \rtensor{\beta}{J_{\mu}}{\gamma}
    J)|\eta\rangle_{2} = |J\eta J_{\mu}\rangle_{2} I$ for all $\xi \in
    \beta$, $\eta \in \gamma$.
  \item $(I \rtensor{\beta}{J_{\mu}}{\gamma} J)(S
    \stensor{\beta}{\gamma} T) =(ISI^{*} \stensor{\tilde\beta}{\tilde
      \gamma} JTJ^{*}) (I \rtensor{\beta}{J_{\mu}}{\gamma} J)$ for all
    $S \in {\cal L}(H_{\beta})$, $T \in {\cal L}(K_{\gamma})$. 
  \end{enumerate}
\end{lemma}
\begin{proof}
Straightforward.
\end{proof}

\paragraph{The fiber product of $C^{*}$-algebras}
The fiber product of $C^{*}$-algebras over KMS-states is an analogue
of the fiber product of von Neumann algebras. We
briefly summarize the definition and main properties; for details,
see \cite[Section 3]{timmer:cpmu}.  

Let $\mu$ be a faithful KMS-state on a $C^{*}$-algebra $B$, let
$A_{H}^{\beta}$ be a $C^{*}$-$\mu$-algebra, and let $C_{K}^{\gamma}$
be a $C^{*}$-$\mu^{op}$-algebra. The {\em fiber product} of  
$A_{H}^{\beta}$ and $C_{K}^{\gamma}$ is  the $C^{*}$-algebra
  \begin{align*}
    A \rfibre{\beta}{\gamma} C := \big\{ x \in {\cal L}(H
    \stensor{\beta}{\gamma} K) \,\big|\, x\kbeta{1},
    x^{*}\kbeta{1} \subseteq [ \kbeta{1} C] \text{ and }
    x\kgamma{2}, x^{*}\kgamma{2} \subseteq [\kgamma{2}A]\big\}.
  \end{align*}
  If $A$ and $C$ are unital, so is $A \rfibre{\beta}{\gamma} C$, but
  otherwise,  $A \rfibre{\beta}{\gamma} C$ may be degenerate.

 Conjugation by the flip $\Sigma \colon H
\stensor{\beta}{\gamma} K \to K \stensor{\gamma}{\beta} H$ in
\eqref{eq:modules-flip-modules} yields an isomorphism 
\begin{align} \label{eq:modules-flip-algebras}
  \Ad_{\Sigma} \colon A \rfibre{\beta}{\gamma} C \to C
  \rfibre{\gamma}{\beta} A.
\end{align}

Assume that $\frakA=(H,\alpha_{1},\ldots,\alpha_{m},\beta,A)$ is a
$C^{*}$-$(\sigma_{1},\ldots,\sigma_{m},\mu)$-algebra and 
$\frakC=(K,\gamma,\delta_{1}$, $\ldots,\delta_{n},C)$ a
$C^{*}$-$(\mu^{op},\tau_{1},\ldots,\tau_{n})$-algebra, where
$\sigma_{1},\ldots,\sigma_{m}, \tau_{1},\ldots, \tau_{n}$ are
faithful KMS-states on $C^{*}$-algebras
$A_{1},\ldots,A_{m},C_{1},\ldots,C_{n}$.  If $A
\rfibre{\beta}{\gamma} C$ is nondegenerate, then
\begin{align*} 
  \frakA \ast \frakC := ((H \stensor{\beta}{\gamma} K, \alpha_{1} \lt
  \gamma, \ldots, \alpha_{m} \lt \gamma, \beta \rt
  \delta_{1},\ldots, \beta \rt \delta_{n}), A \rfibre{\beta}{\gamma}
  C)
\end{align*}
is a $C^{*}$-$(\sigma_{1},\ldots,\sigma_{m},
\tau_{1},\ldots,\tau_{n})$-algebra, called the {\em fiber product}
of $\frakA$ and $\frakC$.

Assume furthermore that $\tilde\frakA= (\tilde H,\tilde
\alpha_{1},\ldots,\tilde \alpha_{m},\tilde \beta, \tilde A)$ is a
$C^{*}$-$(\sigma_{1},\ldots,\sigma_{m},\mu)$-algebra and $\tilde
\frakC=(\tilde K,\tilde \gamma, \tilde \delta_{1},\ldots,\tilde
\delta_{n},\tilde C)$ is a
$C^{*}$-$(\mu^{op},\tau_{1},\ldots,\tau_{n})$-algebra. Then for each
$\phi \in \Mor(\frakA,\tilde \frakA)$ and $\psi \in \Mor(\frakC,
\tilde \frakC)$, there exists a unique morphism
\begin{align*}
 \phi \ast \psi \in \Mor(\frakA \ast \frakC, \tilde \frakA \ast
 \tilde \frakC) 
\end{align*}
such that $(\phi \ast \psi)(x) (S \stensor{\beta}{\gamma} T) = (S
\stensor{\beta}{\gamma} T)x$ for all $x \in A \rfibre{\beta}{\gamma}
C$, $S \in {\cal L}(H_{\beta},\tilde H_{\tilde \beta})$, $T \in
{\cal L}(K_{\gamma},\tilde K_{\tilde \gamma})$ satisfying
$Sa=\phi(a)S$ and $Tc=\psi(c)T$ for all $a \in A$, $c \in C$.

A fundamental deficiency of the fiber product is that it need not be
associative. In our applications, however, the fiber product will only
appear as the target of a comultiplication, and the non-associativity of
the former will be compensated by the coassociativity of the latter.

We shall need the following simple construction not mentioned in
\cite{timmer:cpmu}:
\begin{lemma} \label{lemma:modules-antiunitary-morphism} Let
  $A_{H}^{\beta}$, $\tilde A_{\tilde H}^{\tilde \beta}$ be
  $C^{*}$-$\mu$-algebras,  $C^{\gamma}_{K}$, $\tilde C^{\tilde
    \gamma}_{\tilde K}$  $C^{*}$-$\mu^{op}$-algebras, and  $R
  \colon A \to \tilde A^{op}$, $S \colon C \to \tilde C^{op}$ 
  $*$-homomorphisms. Assume that $I \colon H \to \tilde H$ and $J
  \colon K \to \tilde K$ are anti-unitaries such that $I\beta
  J_{\mu} = \tilde \beta$, $R(a)=I^{*}a^{*}I$ for all $a\in A$, and
  $J \gamma J_{\mu} = \tilde \gamma$, $S(c)=J^{*}c^{*}J$ for all $c
  \in C$. Then there exists a $*$-homomorphism $R
  \rfibre{\beta}{\gamma} S \colon A \rfibre{\beta}{\gamma} C \to
  (\tilde A \rfibre{\tilde \beta}{\tilde \gamma} \tilde C)^{op}$
  such that $(R \rfibre{\beta}{\gamma} S)(x):=(I
  \rtensor{\beta}{J_{\mu}}{\gamma} J)^{*}x^{*}(I
  \rtensor{\beta}{J_{\mu}}{\gamma} J)$ for all $x \in A
  \rfibre{\beta}{\gamma} C$.  This $*$-homomorphism does not depend
  on the choice of $I$ or $J$.
\end{lemma}
\begin{proof}
  Evidently, the formula defines a $*$-homomorphism $R
  \rfibre{\beta}{\gamma} S$. The definition does not depend on the
  choice of $J$ because
$\langle \xi|_{1} (R \rfibre{\beta}{\gamma} S)(x) |\xi'\rangle_{1} =
    J^{*} \langle I \xi J_{\mu}|_{1} x^{*} |I\xi' J_{\mu}\rangle_{1} J
     = S\big(  \langle I\xi' J_{\mu}|_{1}x|I\xi J_{\mu}\rangle_{1}\big)$
  for all $x \in A \rfibre{\beta}{\gamma} C$ by Lemma
  \ref{lemma:modules-antiunitary} ii), and a similar argument shows that
  it does not depend on the choice of $I$.
\end{proof}

\section{Compact $C^{*}$-quantum groupoids}

\label{section:cqg}

In this section, we introduce the main object of study of this article
--- compact $C^{*}$-quantum groupoids. Roughly, a compact
$C^{*}$-quantum groupoid is a compact $C^{*}$-quantum graph equipped
with a coinvolution and a comultiplication subject to several relations.
Most importantly, we assume left- and right-invariance of the Haar
weights, the existence of a modular element, and a strong invariance
condition relating the coinvolution to the  Haar weights and to the
comultiplication.

We proceed as follows.  First, we discuss the appropriate notion of a
comultiplication and recall the notion of a Hopf $C^{*}$-bimodule, of
bounded invariant Haar weights, and of bounded counits.  Then, we
introduce and study the precise definition of a compact $C^{*}$-quantum
groupoid. Finally, we show that the modular element can always
be assumed to be trivial, and that the Haar weights are unique up to
scaling.

\paragraph{Hopf $C^{*}$-bimodules over KMS-states}
Throughout this paragraph, let $\mu$ be a faithful KMS-state on a
$C^{*}$-algebra $B$.
\begin{definition}[{\cite{timmer:cpmu}}] \label{definition:cqg-comultiplication}
  A {\em comultiplication} on a $C^{*}$-$(\mu,\mu^{op})$-algebra
  $A^{\alpha,\beta}_{H}$ is a morphism $\Delta \in
  \Mor(A^{\alpha,\beta}_{H}, A^{\alpha,\beta}_{H} \ast
  A^{\alpha,\beta}_{H})$  that makes
  the following diagram commute:
    \begin{align*}
      \xymatrix@C=15pt@R=10pt{ A \ar[rrr]^{\Delta} \ar[dd]^{\Delta} &&&
        {\AfibreA} \ar[d]^{\Id
          \ast \Delta} \\
        &&& {A \rfibre{\alpha}{\beta \lt \beta} (\AfibreA)}
        \ar@{^(->}[d] \\
        {\AfibreA} \ar[rr]^(0.35){\Delta\ast \Id} && {(\AfibreA)
          \rfibre{\alpha \rt \alpha}{\beta} A} \ar@{^(->}[r]
        &{\cal L}(H \fibreab H \fibreab H).  }
    \end{align*}
    A {\em Hopf $C^{*}$-bimodule over $\mu$} is a
    $C^{*}$-$(\mu,\mu^{op})$-algebra together with a comultiplication.
\end{definition}

The following important invariance conditions will be imposed on the
Haar weights of a compact $C^{*}$-quantum groupoid:
\begin{definition} \label{definition:cqg-haar-weights}
  Let $(A_{H}^{\alpha,\beta},\Delta)$ be a  Hopf
  $C^{*}$-bimodule over $\mu$.  A {\em bounded left Haar weight} for
  $(A_{H}^{\alpha,\beta},\Delta)$ is a non-zero completely
  positive contraction $\phi \colon A \to B$ satisfying
 \begin{enumerate}
 \item $\phi\big(\rho_{\beta}(b)a\rho_{\beta}(c)\big) = b\phi(a)c$ for
   all $a \in A$ and $b,c \in B$,
 \item $\phi\big(\langle \xi|_{1}\Delta(a)|\xi'\rangle_{1}\big) =
   \xi^{*}\rho_{\beta}(\phi(a))\xi'$ for all $a \in A$ and $\xi,\xi'
   \in \alpha$.
 \end{enumerate}
 A {\em bounded right Haar weight} for $(A_{H}^{\alpha,\beta},\Delta)$
 is a non-zero completely positive contraction $\psi \colon A \to
 B^{op}$ satisfying
 \begin{enumerate} 
 \item[i)'] $\psi\big(\rho_{\alpha}(b^{op})a\rho_{\alpha}(c^{op})\big) =
   b^{op}\psi(a)c^{op}$ for all $a \in A$ and $b,c \in B$,
 \item[ii)']
   $\psi\big(\langle\eta|_{2}\Delta(a)|\eta'\rangle_{2}\big) =
   \eta^{*}\rho_{\alpha}(\psi(a))\eta'$ for all $a \in A$ and
   $\eta,\eta' \in \beta$.
 \end{enumerate}
\end{definition}
\begin{remarks} \label{remarks:cqg-haar-weights}
  Let $(A_{H}^{\alpha,\beta},\Delta)$ be a Hopf
  $C^{*}$-bimodule over $\mu$.
  \begin{enumerate}
  \item If $\phi$ is a bounded left Haar weight for
    $(A_{H}^{\alpha,\beta},\Delta)$ , then $\rho_{\beta} \circ \phi
    \colon A \to \rho_{\beta}(\frakB)$ is a conditional expectation.
  \item If $\phi \colon A \to B$ satisfies condition ii) and if
    $[\balpha{1}\Delta(A)\kalpha{1}]=A$, then $\phi$ also satisfies
    condition i) because $\phi\big(\rho_{\beta}(b)\langle
    \xi|_{1}\Delta(a)|\xi'\rangle_{1}\rho_{\beta}(c)\big) =
    \phi\big(\langle \xi b|_{1}\Delta(a)|\xi'c\rangle_{1}\big) =
    b^{*}\xi^{*}\rho_{\beta} (\phi(a)) \xi'c = b^{*} \phi\big(\langle
    \xi|_{1}\Delta(a)|\xi'\rangle_{1}\big)c$ for all $a \in A$, $b,c \in
    B$, $\xi,\xi' \in \alpha$.
  \end{enumerate}
  Similar remarks apply to bounded right Haar weights.
\end{remarks}

The notion of a counit of a Hopf algebra extends to Hopf
$C^{*}$-bimodules as follows.
\begin{definition}
  Let $( A_{ H}^{ \alpha, \beta}, \Delta)$ be a Hopf
  $C^{*}$-bimodule over $\mu$.  A {\em bounded (left/right) counit}
  for $( A_{ H}^{ \alpha, \beta}, \Delta)$ is a morphism $ \epsilon
  \in \Mor\big( A^{\alpha,\beta}_{ H}, {\cal L}(H_{\mu})^{ B,
    B^{op}}_{H_{ \mu}}\big)$ satisfying (the first/second of) the
  following conditions:
  \begin{enumerate}
  \item $\epsilon\big(\langle
    \eta|_{2}\Delta(a)|\eta'\rangle_{2}\big) = \eta^{*}a\eta'$ for
    all $a \in A$ and $\eta,\eta' \in \beta$,
  \item  $\epsilon\big(\langle
  \xi|_{1}\Delta(a)|\xi'\rangle_{1}\big) = \xi^{*}a\xi'$ for all $a
  \in A$ and $\xi,\xi' \in \alpha$.
\end{enumerate}
\end{definition}
\begin{remark}
  \begin{enumerate}
  \item Condition i) and ii), respectively, hold if and only if the
    left/the right square of the following diagram commute:
    \begin{align*}
      \xymatrix@C=15pt@R=15pt{ { A \rfibre{ \alpha}{\beta} A}
        \ar[d]_{ \epsilon \ast \Id} && A \ar[ll]_{ \Delta} \ar[d]
        \ar[rr]^{ \Delta} \ar[d] && { A \rfibre{\alpha}{\beta} A}
        \ar[d]^{\Id \ast \epsilon}
        \\
        {{\cal L}(H_{\mu}) \rfibre{ B}{ \beta} A} \ar[r] & {{\cal
            L}(H_{\mu} \stensor{ B}{ \beta} H)} \ar[r]^(0.6){\cong}
        & {{\cal L}( H)} & {{\cal L}( H \stensor{ \alpha}{ B^{op}}
          H_{\mu} )} \ar[l]_(0.6){\cong} & { A \rfibre{ \alpha}{
            B^{op}} {\cal L}(H_{\mu}).} \ar[l]}
    \end{align*}
  \item A standard argument shows that if a bounded left and a
    bounded right counit exist, then they are equal and a counit.
  \end{enumerate}
\end{remark}

\paragraph{Compact $C^{*}$-quantum groupoids}
Given a compact $C^{*}$-quantum graph
$(B,\mu,A,r,\phi$, $s,\psi,\delta)$ with coinvolution $R$, we use the
notation introduced in Proposition \ref{proposition:module-graph}
and put $\nu:=\mu \circ \phi$, $\nu^{-1}:=\mu^{op} \circ \psi =
\nu_{\delta}$, $J:=J_{\nu} = J_{\nu^{-1}}$,
\begin{align} \label{eq:cqg-modules}
  H &:= H_{\nu}, & \halpha &:=[A\zeta_{\phi}], & \beta
  &:=[A^{op}\zeta_{\phi}], & \hbeta &:=[A\zeta_{\psi}], & \alpha &:=
  [A^{op}\zeta_{\psi}],
\end{align}
and define an antiunitary $I \colon H \to H$ by
$I\Lambda_{\nu^{-1}}(a)=\Lambda_{\nu}(R(a)^{*})$ for all $a \in A$.
Since $I\alpha J_{\mu} =\beta$, $I\beta J_{\mu}=\alpha$, and
$R(a)=Ia^{*}I$ for all $a \in A$, we can define a $*$-antihomomorphism
$R \rfibre{\alpha}{\beta} R \colon A \rfibre{\alpha}{\beta} A \to A
\rfibre{\beta}{\alpha} A$ by $x \mapsto (I
\rtensor{\alpha}{J_{\mu}}{\beta} I)^{*}x^{*}(I
\rtensor{\alpha}{J_{\mu}}{\alpha} I)$ (Lemma
\ref{lemma:modules-antiunitary-morphism}). 

The definition of a compact $C^{*}$-quantum groupoid involves the
following conditions that are analogues of the strong invariance
property known from quantum groups:
\begin{lemma} \label{lemma:cqg-strong-invariance} Let
  $(B,\mu,A,r,\phi,s,\psi,\delta)$ be a compact
  $C^{*}$-quantum graph with a coinvolution $R$ and a comultiplication
  $\Delta$ for $A^{\alpha,\beta}_{H}$ such that $(R
  \rfibre{\alpha}{\beta} R) \circ \Delta = \Ad_{\Sigma} \circ \Delta
  \circ R$. Then 
  \begin{align*}
    R\big(\langle \zeta_{\psi}|_{1}\Delta(a)(d^{op}
    \stensor{\alpha}{\beta} 1)|\zeta_{\psi}\rangle_{1}\big) =
    \langle\zeta_{\phi}|_{2} (1 \stensor{\alpha}{\beta}
    R(d)^{op})\Delta(R(a))|\zeta_{\phi}\rangle_{2} \quad \text{for
      all } a,d \in A.
  \end{align*}
\end{lemma}
\begin{proof}
  Let $a,d \in A$.    By Lemmas \ref{lemma:modules-antiunitary} and 
  \ref{lemma:graph-coinvolution},
  \begin{align*}
    \langle\zeta_{\phi}|_{2} (1 \stensor{\alpha}{\beta}
    R(d)^{op})\Delta(R(a))|\zeta_{\phi}\rangle_{2} &= \langle
    \zeta_{\phi}|_{2} (1 \stensor{\alpha}{\beta} I(d^{op})^{*}I)
    \Sigma (I \stensor{\alpha}{\beta} I)\Delta(a)^{*}(I
    \rtensor{\alpha}{J_{\mu}}{\beta} I)^{*} \Sigma
    |\zeta_{\phi}\rangle_{2}
    \\
    &= \langle \zeta_{\phi}|_{1} (I \rtensor{\alpha}{J_{\mu}}{\beta}
    I)((d^{op})^{*} \stensor{\alpha}{\beta} 1)\Delta(a)^{*}
    |I\zeta_{\phi}J_{\mu}\rangle_{1} I \\
    &= I \langle \zeta_{\psi}|_{1} ((d^{op})^{*}
    \stensor{\alpha}{\beta}
    1)\Delta(a)^{*} |\zeta_{\psi}\rangle_{1} I \\
    &= R\big( \langle \zeta_{\psi}|_{1} \Delta(a) (d^{op}
    \stensor{\alpha}{\beta} 1) |\zeta_{\psi}\rangle_{1}
    \big). \qedhere
  \end{align*}
  \end{proof}
As a direct consequence, we obtain the following result:
\begin{lemma}
  \label{lemma:cqg-strong-invariance-equiv} Let
  $(B,\mu,A,r,\phi,s,\psi,\delta)$ be a compact $C^{*}$-quantum graph
  with a coinvolution $R$ and a comultiplication $\Delta$ for
  $A^{\alpha,\beta}_{H}$ such that $(R \rfibre{\alpha}{\beta} R) \circ
  \Delta = \Ad_{\Sigma} \circ \Delta \circ R$. Then the following two
  conditions are equivalent:
 \begin{enumerate}
 \item $R\big(\langle\zeta_{\psi}|_{1}\Delta(a)(d^{op} \stensor{\alpha}{\beta}
   1)|\zeta_{\psi}\rangle_{1}\big) = \langle\zeta_{\psi}|_{1}(a^{op} \stensor{\alpha}{\beta}
   1)\Delta(d)|\zeta_{\psi}\rangle_{1}$ for all $a,d \in A$.
 \item $R\big(\langle\zeta_{\phi}|_{2}\Delta(a)(1
   \stensor{\alpha}{\beta} d^{op})|\zeta_{\phi}\rangle_{2}\big) =
   \langle \zeta_{\phi|_{2}}(1 \stensor{\alpha}{\beta}
   a^{op})\Delta(d)|\zeta_{\phi}\rangle_{2}$ for all $a,d \in
   A$. \qed
 \end{enumerate}
\end{lemma}
Now we come to the main definition of this article:
\begin{definition} \label{definition:cqg}
  A {\em compact $C^{*}$-quantum groupoid} is a compact
  $C^{*}$-quantum graph $(B,\mu,A$, $r,\phi,s,\psi,\delta)$
  with a coinvolution $R$ and a comultiplication $\Delta$ for
  $A^{\alpha,\beta}_{H}$ such that
    \begin{enumerate}
    \item $[\Delta(A)\kalpha{1}]=[\kalpha{1}A] =
      [\Delta(A)|\zeta_{\psi}\rangle_{1}A]$ and
      $[\Delta(A)\kbeta{2}]=[\kbeta{2}A] =
      [\Delta(A)|\zeta_{\phi}\rangle_{2}A]$;
  \item $\phi$ is a bounded left Haar weight and $\psi$ a bounded right Haar
    weight for $(A^{\alpha,\beta}_{H},\Delta)$;
  \item $R\big(\langle\zeta_{\psi}|_{1}\Delta(a)(d^{op}
    \stensor{\alpha}{\beta} 1)|\zeta_{\psi}\rangle_{1}\big) =
    \langle\zeta_{\psi}|_{1}(a^{op} \stensor{\alpha}{\beta}
    1)\Delta(d)|\zeta_{\psi}\rangle_{1}$ for all $a,d \in A$.
  \end{enumerate}
\end{definition}
Let us briefly comment on this definition.  The coinvolution $R$ is
uniquely determined by condition iii).  The Haar weights are unique
up to some scaling, as we shall see at the end of this section. At
the end of the next section, we will see that $(R
\rfibre{\alpha}{\beta} R) \circ \Delta = \Ad_{\Sigma} \circ \Delta
\circ R$; in particular,  the modified strong invariance
condition in Lemma \ref{lemma:cqg-strong-invariance-equiv} ii) holds
by Lemma \ref{lemma:cqg-strong-invariance-equiv}.

From now on, let $(B,\mu,A,r,\phi,s,\psi,\delta,R,\Delta)$ be a
compact $C^{*}$-quantum groupoid.

\begin{lemma} \label{lemma:cqg-rs} $ \{ a \in A \cap r(B)' \mid
  \Delta(a) = 1 \stensor{\alpha}{\beta} a \} = s(B^{op})$ and $ \{ a
  \in A \cap s(B^{op})' \mid \Delta(a) = a \stensor{\alpha}{\beta} 1
  \} = r(B)$.
\end{lemma}
\begin{proof}
  We only prove the first equation.  Clearly, the right hand side is
  contained in the left hand side. Conversely, if $a \in A \cap
  r(B)'$ and $\Delta(a)=1 \stensor{\alpha}{\beta} a$, then $a=
  \langle \zeta_{\psi}|_{1} \Delta(a)|\zeta_{\psi}\rangle_{1} =
  s(\psi(a))$ by right-invariance of $\psi$.
\end{proof}

\paragraph{The conditional expectation onto the $C^{*}$-algebra of
  orbits}
Let us study the maps $\tau=\psi \circ r \colon B \to Z(B^{op})
\cong Z(B)$ and $\tau^{\dag} = \phi \circ s \colon B^{op} \to
Z(B)\cong Z(B^{op})$ introduced in \eqref{eq:graph-tau}). First,
note that $\tau(b)=\tau^{\dag}(b^{op})$ for all $b \in B$ by Lemma
\ref{lemma:graph-coinvolution} ii).
\begin{proposition} \label{proposition:cqg-tau} The maps $\tau$ and
  $\tau^{\dag}$ are conditional expectations onto a $C^{*}$-subalgebra of
  $Z(B)=B \cap B^{op}$ and satisfy 
  \begin{align*}
    s \circ \tau &= r \circ \tau, &
    \sigma^{\mu}_{t} \circ \tau &=\tau \text{ for all } t \in \reals, \\
    \tau \circ \phi &= \tau^{\dag} \circ
    \psi,  &
    \tau(b\sigma_{-i/2}^{\mu}(d)) &= \tau(d\sigma^{\mu}_{-i/2}(d))
    \text{ for  all } b,d \in \Dom(\sigma^{\mu}_{-i/2}).
  \end{align*}
\end{proposition}
The proof involves the following equation:
\begin{lemma} \label{lemma:cqg-invariance-tau}
  For all $b,c,e \in B$ and $d \in \Dom(\sigma^{\mu}_{-i/2})$,
  \begin{align*}
    \langle \zeta_{\psi}|_{1}\Delta(r(b)s(c^{op}))((r(d)s(e^{op}))^{op}
    \stensor{\alpha}{\beta} 1) |\zeta_{\psi}\rangle_{1}
    = r(\tau(b\sigma^{\mu}_{-i/2}(d))) r(e) s(c^{op}) .
  \end{align*}
\end{lemma}
\begin{proof}
  Let $b,c,d,e$ as above. Then
    \begin{align*}
    \langle \zeta_{\psi}|_{1}\Delta(r(b)s(c^{op}))((r(d)s(e^{op}))^{op}
    \stensor{\alpha}{\beta} 1) |\zeta_{\psi}\rangle_{1} &= \langle
    \zeta_{\psi}|_{1}(r(b)(r(d)s(e^{op}))^{op}
    \stensor{\alpha}{\beta} s(c^{op})) |\zeta_{\psi}\rangle_{1} \\
    &=
    \rho_{\beta}\big(\zeta_{\psi}^{*}r(b)r(d)^{op}s(e^{op})^{op}\zeta_{\psi}\big) s(c^{op}) \\
    &= r\big(\zeta_{\psi}^{*}r(b)r(d)^{op}\zeta_{\psi}e\big)
    s(c^{op}),
  \end{align*}
and  by Lemma
  \ref{lemma:module} iii),
  $\zeta_{\psi}^{*}r(b)r(d)^{op}\zeta_{\psi} =
  \zeta_{\psi}^{*}r(b\sigma^{\mu}_{-i/2}(d))\zeta_{\psi} =
  \tau(b\sigma^{\mu}_{-i/2}(d))$.
\end{proof}
\begin{proof}[Proof of Proposition \ref{proposition:cqg-tau}]
  The left- and right-invariance of $\phi$ and $\psi$ imply that for all
  $a \in A$,
  \begin{align*}
    \phi(s(\psi(a)))  = \zeta_{\phi}^{*}s(\psi(a))\zeta_{\phi} 
    &= \zeta_{\psi}^{*} \langle \zeta_{\phi}|_{2} \Delta(a)
    |\zeta_{\phi}\rangle_{2} \zeta_{\psi} \\
    &= \zeta_{\phi}^{*} \langle
    \zeta_{\psi}|_{1} \Delta(a)
    |\zeta_{\psi}\rangle_{1} \zeta_{\phi}^{*} 
    =  \zeta_{\psi}^{*}r(\phi(a))\zeta_{\psi} = \psi(r(\phi(a))).
  \end{align*}
  Therefore, $\tau^{\dag} \circ \psi = \tau \circ \phi$ and $\tau
  \circ \tau = \tau^{\dag} \circ \tau = \tau^{\dag} \circ (\psi
  \circ r) = \tau \circ \phi \circ r = \tau$.  Next, $s\circ \tau =
  r \circ \tau$ because for all $b \in B$,
  \begin{align*}
    s(\psi(r(b))) =
    \langle\zeta_{\psi}|_{1}\Delta(r(b))|\zeta_{\psi}\rangle_{1} 
=    \langle\zeta_{\psi}|_{1}(r(b) \stensor{\alpha}{\beta}
    1)|\zeta_{\psi}\rangle_{1}  =
    \rho_{\beta}\big(\zeta_{\psi}^{*}r(b)\zeta_{\psi}) = r(\psi(r(b))).
  \end{align*}
  In particular, we find that for all $b,c,d \in B$,
  \begin{align*}
    \tau(b)\tau(c)\tau(d) = \tau(b)\psi(r(c))\tau(d) =
    \psi\big(s(\tau(b))r(c)s(\tau(d))\big) 
    = \psi\big(r(\tau(bcd)))  = \tau(bcd).
  \end{align*}
  Therefore, $\tau$ is a conditional expectation onto its image.

  Let $t \in \reals$. Then   $\sigma^{\mu}_{t}(\tau(B)) \subseteq
  \tau(B)$ because $\sigma^{\mu}_{t} \circ \tau =
  \sigma^{\mu^{op}}_{-t} \circ \psi \circ r = \psi \circ
  \sigma^{\nu^{-1}}_{-t} \circ r = \psi \circ \sigma^{\nu}_{t} \circ
  r = \psi \circ r \circ \sigma^{\mu}_{t}$. Since
  $\upsilon:=\mu|_{\tau(B)}$ is a trace, we can conclude from Lemma
  \ref{lemma:graph-rieffel} i) that $\sigma^{\mu}_{t} \circ \tau =
  \tau \circ \sigma^{\upsilon}_{t} = \tau$.
  
  Finally, let $b,d \in \Dom(\sigma^{\mu}_{-i/2})$. By Lemma
  \ref{lemma:cqg-invariance-tau} and condition iii) in Definition
  \ref{definition:cqg},
  \begin{align*}
    r(\tau(b\sigma^{\mu}_{-i/2}(d))) &= \langle
    \zeta_{\psi}|_{1}\Delta(r(b))(r(d)^{op} \stensor{\alpha}{\beta} 1)
    |\zeta_{\psi}\rangle_{1} \\
    &= R\big( \langle
    \zeta_{\psi}|_{1}\Delta(r(d))(r(b)^{op} \stensor{\alpha}{\beta} 1)
    |\zeta_{\psi}\rangle_{1} \big) =
    s(\tau(d\sigma^{\mu}_{-i/2}(b))).
  \end{align*}
  Since $s \circ \tau = r\circ \tau$ and $r$ is
  injective, we can conclude
  $\tau(b\sigma^{\mu}_{-i/2}(d))=\tau(d\sigma^{\mu}_{-i/2}(b))$.
\end{proof}

\paragraph{The modular element}
The modular element of a compact $C^{*}$-quantum groupoid can be
described in terms of the element
\begin{align*}
  \theta:=  \phi(\delta) = \psi(\delta^{-1}) \in B \cap B^{op}
\end{align*}
(see Lemma \ref{lemma:graph} and Lemma
\ref{lemma:graph-coinvolution} i)) as follows. 
\begin{proposition}  \label{proposition:cqg-modular}
 $\delta = r(\theta)s(\theta)^{-1}$ and $\Delta(\delta) = \delta
    \stensor{\alpha}{\beta} \delta$.
\end{proposition}
\begin{proof}
  By Lemma \ref{lemma:graph} i), the element $\tilde \delta:=
  r(\theta)s(\theta)^{-1}$ is positive, invertible, and invariant
  with respect to $\sigma^{\nu}$. Moreover, $\nu^{-1}(a) =
  \nu(\tilde \delta^{1/2}a\tilde \delta^{1/2})$ for all $a \in A$
  because 
  \begin{multline*}
    \nu^{-1}\big(s(\theta)^{1/2}as(\theta)^{1/2}\big) =
    \mu^{op}(\theta^{1/2}\psi(a)\theta^{1/2})
    = (\mu \circ \phi \circ  s \circ \psi)(a)  \\
    = (\mu^{op} \circ \psi \circ r \circ \phi)(a)
   = \mu(\theta^{1/2}\phi(a)\theta^{1/2}) =
    \nu\big(r(\theta)^{1/2}ar(\theta)^{1/2}\big)
  \end{multline*}
  for all $a \in A$ by Proposition \ref{proposition:cqg-tau} and Lemma
  \ref{lemma:graph} ii). Now,  $\delta=\tilde \delta$ by Lemma
  \ref{lemma:graph-rn} iii), and
  $ \Delta(\delta) = r(\theta) \stensor{\alpha}{\beta}
  s(\theta)^{-1} = r(\theta) \rho_{\alpha}(\theta^{-1})
  \stensor{\alpha}{\beta} \rho_{\beta}(\theta) s(\theta)^{-1} =
  \delta \stensor{\alpha}{\beta} \delta$ because $\theta \in Z(B)$.
\end{proof}
An important consequence of the preceding result is that for every
compact $C^{*}$-quantum groupoid, there exists a faithful invariant
KMS-state on the basis:
\begin{corollary} \label{corollary:cqg-mu}
    $\mu_{\theta} \circ \phi = (\mu_{\theta})^{op} \circ \psi$.
\end{corollary}
\begin{proof}
    For all $a \in A$,  we have
    $\mu\big(\theta^{1/2}\phi(a)\theta^{1/2}\big) =
   \nu\big(r(\theta)^{1/2}ar(\theta)^{1/2}\big) =
   \nu^{-1}\big(s(\theta)^{1/2}as(\theta)^{1/2}\big) =
   \mu^{op}\big(\theta^{1/2}\psi(a)\theta^{1/2}\big)$.
\end{proof}
This result implies that in principle, we could restrict to compact
$C^{*}$-quantum groupoids with trivial modular element $\delta=1_{A}$.
We shall not do so for several reasons. First, the treatment of
a nontrivial modular element does not require substantially more work.
Second, the freedom to choose the state $\mu$ might be useful in 
applications. Finally, we hope to prepare the ground for a more general
theory of locally compact quantum groupoids, where the modular element
can no longer be assumed to be trivial.

The KMS-state $\mu$ can be factorized into a state $\upsilon$ on the
commutative $C^{*}$-algebra $\tau(B) \subseteq Z(B)$ and a perturbation
of $\tau$ as follows. We define maps
\begin{align*}
  \tau_{\theta^{-1}} \colon B &\to \tau(B), \
  b \mapsto \tau(\theta^{-1/2}b\theta^{-1/2}), &
  \upsilon  = \mu_{\theta}|_{\tau(B)}\colon \tau(B) &\to \complex, \ b \mapsto
  \mu(\theta^{1/2}b\theta^{1/2}).
\end{align*}
Note that $\tau(\theta^{-1})=1$ because  $\theta =\phi(\delta) =
  \phi(r(\theta)s(\theta)^{-1}) = \theta \tau(\theta^{-1})$.
\begin{proposition} \label{proposition:cqg-mu} 
    $\mu = \upsilon \circ \tau_{\theta^{-1}}$.
\end{proposition}
\begin{proof}
  By Propositions \ref{proposition:cqg-tau} and
  \ref{proposition:cqg-modular}, we have $\mu(b) = \nu(r(b)) =
  \nu^{-1}(\delta^{-1/2}r(b)\delta^{-1/2}) =
  \mu^{op}\big(\theta^{1/2} \psi(r(\theta^{-1/2}b\theta^{-1/2}))
  \theta^{1/2}\big) = (\upsilon \circ \tau_{\theta^{-1}})(b)$ for all $b
  \in B$.
 \end{proof}

\paragraph{Uniqueness of the Haar weights}
A central result in the theory of locally compact (quantum) groups is
the uniqueness of the Haar weights up to scaling. In this
paragraph, we prove a similar uniqueness result for the Haar weights of a
compact $C^{*}$-quantum groupoid. 

The Haar weights of a compact $C^{*}$-quantum groupoid can be
rescaled by elements of $B$ as follows: For every positive element $\gamma
\in B^{op}$,  the completely positive contraction
\begin{align*}
  \phi_{s(\gamma)} \colon A \to B, \ a \mapsto
  \phi\big(s(\gamma)^{1/2}as(\gamma)^{1/2}\big),
\end{align*}
is a bounded left Haar weight for $(A^{\alpha,\beta}_{H},\Delta)$
because for all $a \in A$ and $\xi,\xi' \in \alpha$
\begin{align*}
  \phi_{s(\gamma)}\big(\langle \xi|_{1}\Delta(a)|\xi'\rangle_{1}\big) &=
  \phi\big(\langle \xi|_{1}(1 \stensor{\alpha}{\beta}
  s(\gamma)^{1/2})\Delta(a) (1 \stensor{\alpha}{\beta}
  s(\gamma)^{1/2}) |\xi'\rangle_{1}\big)  \\
  &=\phi\big(\langle
  \xi|_{1}\Delta\big(s(\gamma)^{1/2}as(\gamma^{1/2})\big)|\xi'\rangle_{1}
  \big) = \xi^{*} \phi_{s(\gamma)}(a)\xi'.
\end{align*}
Likewise, for every positive element $\gamma
\in B$,  the completely positive contraction
\begin{align*}
  \psi_{r(\gamma)} \colon A \to B^{op}, \ a \mapsto
  \psi\big(r(\gamma)^{1/2}ar(\gamma)^{1/2}\big),
\end{align*}
is a bounded right Haar weight for $(A^{\alpha,\beta}_{H},\Delta)$.
\begin{theorem} \label{theorem:cqg-unique-weight} Let $\tilde \phi$,
  $\tilde \psi$, $\tilde \delta$ be such that $(B,\mu,A,r,\tilde
  \phi,s,\tilde\psi,\tilde \delta)$ is a compact $C^{*}$-quantum
  graph.
  \begin{enumerate}
  \item If $\tilde \phi$ is a bounded left Haar weight for
    $(A^{\alpha,\beta}_{H},\Delta)$, then $\tilde \phi = \phi_{\gamma}$,
    where $\gamma=\tilde \psi(\tilde \delta^{-1}) \theta^{-1}$.
  \item If $\tilde \psi$ is a bounded right Haar weight for
    $(A^{\alpha,\beta}_{H},\Delta)$, then $\tilde \psi = \psi_{\gamma}$,
    where $\gamma=\tilde \phi(\tilde \delta) \theta^{-1}$.
  \end{enumerate}
\end{theorem}
\begin{proof}
  We only prove i), the proof of ii) is similar.  Put $\tilde \nu:=\mu
  \circ \tilde \phi$, $\tilde \nu^{-1}:=\mu^{op} \circ \tilde \psi$,
  $\tilde \theta:=\tilde \psi(\tilde\delta^{-1})$.  Let $a \in A$. Then
  \begin{align} \label{eq:cqg-unique-weight}
    \tilde \phi(s(\psi(a))) = \tilde \phi(\langle 
    \zeta_{\psi}|_{1}\Delta(a)|\zeta_{\psi}\rangle_{1}) =
    \psi(r(\tilde \phi(a))).
  \end{align}
  We apply $\mu$ to the left hand side and find, using Lemma
  \ref{lemma:graph} ii),
  \begin{align*}
    \tilde \nu(s(\psi(a))) = \mu^{op}_{ \tilde \theta}(\psi(a)) =
    \nu^{-1}\big(s(\tilde \theta)^{1/2}as(\tilde \theta)^{1/2}\big)
    = \nu\big(\delta^{1/2} s(\tilde \theta)^{1/2}as(\tilde
    \theta)^{1/2 } \delta^{1/2}\big).
  \end{align*}
  Next, we apply $\mu$ to the right hand side of equation
  \eqref{eq:cqg-unique-weight} and find
  \begin{align*}
    \nu^{-1}(r(\tilde \phi(a))) =\mu_{\theta}(\tilde \phi(a)) =
    \tilde \nu\big(r(\theta)^{1/2}ar(\theta)^{1/2}\big).
  \end{align*}
  Since the left hand side and the right hand side of equation
  \eqref{eq:cqg-unique-weight} are equal and $\delta
  =r(\theta)s(\theta)^{-1}$, we can conclude
  $\tilde\nu(d) = \nu\big(s(\gamma)^{1/2}ds(\gamma)^{1/2}\big)$ for all $d \in
  A$ and in particular
  \begin{align*}
    \mu(b^{*}\tilde \phi(a)) = \tilde \nu\big(r(b)^{*}a\big) =
    \nu\big(s(\gamma)^{1/2}r(b)^{*}as(\gamma)^{1/2}\big) =
    \mu\big(b^{*}\phi(s(\gamma)^{1/2}as(\gamma)^{1/2})\big) 
  \end{align*}
  for all $b \in B$, $a \in A$.  Since $\mu$ is faithful, we have $\tilde \phi(a) = \phi(s(\gamma)^{1/2} a s(\gamma)^{1/2})$
  for all $a \in A$.
\end{proof}

 \section{The fundamental unitary}

\label{section:pmu}

In the theory of locally compact quantum groups, a fundamental
r\^ole is played by the multiplicative unitaries of Baaj, Skandalis
\cite{baaj:2} and Woronowicz \cite{woron}: To every locally compact
quantum group, one can associate a manageable multiplicative
unitary, and to every manageable multiplicative unitary two Hopf
$C^{*}$-algebras called the ``legs'' of the unitary. One of these
legs coincides with the initial quantum group, and the other one is
its generalized Pontrjagin dual.  Moreover, the multiplicative
unitary can be used to switch between the reduced $C^{*}$-algebra
and the von Neumann algebra of the quantum group.

Similarly, we associate to every compact $C^{*}$-quantum groupoid a
generalized multiplicative unitary.  More precisely, this unitary is
a regular $C^{*}$-pseudo-multiplicative unitary in the sense of
\cite{timmer:cpmu}. The first application of this unitary will be to
prove that the coinvolution of a compact $C^{*}$-quantum groupoid
reverses the comultiplication.  The second application will be to
associate to every compact $C^{*}$-quantum groupoid a measured
quantum groupoid in the sense of Enock and Lesieur
\cite{enock:lesieur,lesieur}. The third application, given in the
next section, will be to construct a generalized Pontrjagin dual of
the compact $C^{*}$-quantum groupoid in form of a Hopf
$C^{*}$-bimodule. Finally, one can use this unitary to define
reduced crossed products for coactions of the compact
$C^{*}$-quantum groupoid as in \cite{timmer:ckac}.

\paragraph{$C^{*}$-pseudo-multiplicative unitaries}
The notion of a $C^{*}$-pseudo-multiplicative unitary extends the
notion of a multiplicative unitary \cite{baaj:2}, of a continuous
field of multiplicative unitaries \cite{blanchard}, and of a
pseudo-multiplicative unitary on $C^{*}$-modules
\cite{ouchi,timmermann:hopf}, and is closely related to
pseudo-multiplicative unitaries on Hilbert spaces \cite{vallin:2};
see \cite[Section 4.1]{timmer:cpmu}. The precise definition is as
follows.  Let $\mu$ be a faithful KMS-state on a $C^{*}$-algebra $B$.
\begin{definition}[{\cite{timmer:cpmu}}]   \label{definition:pmu}
  A
  {\em  $C^{*}$-pseudo-multiplicative unitary over $\mu$} consists
  of a $C^{*}$-$(\mu^{op},\mu,\mu^{op})$-module
  $(H,\hbeta,\alpha,\beta)$ and a unitary $V \colon \sHsource \to
  \sHrange$ such that
  \begin{gather} \label{eq:pmu-intertwine}
    \begin{aligned}
      V(\alpha \lt \alpha) &= \alpha \rt \alpha, &
      V(\hbeta \rt \beta) &= \hbeta \lt \beta, &
      V(\hbeta \rt \hbeta) &= \alpha \rt \hbeta, & 
      V(\beta \lt \alpha) &= \beta \lt \beta
    \end{aligned}
  \end{gather}
  and the following diagram commutes:
  \begin{gather} \label{eq:pmu-pentagon}
    \begin{gathered}
      \xymatrix@R=25pt@C=40pt{ {\sHone} \ar[r]^{\displaystyle V
          \stensor{\hbeta \rt \hbeta}{\alpha} \Id}
        \ar[d]_{\displaystyle \Id \stensor{\hbeta}{\alpha \lt
            \alpha} V} & {\sHtwo} \ar[r]^{\displaystyle \Id
          \stensor{\alpha}{\beta \lt \alpha} V}& { \sHthree,}
        \\
        {\sHfive} \ar[d]_{\displaystyle \Id \stensor{\hbeta}{\alpha
            \rt \alpha} \Sigma} & & {\sHfour} \ar[u]_{\displaystyle
          V \stensor{\alpha \lt \alpha}{\beta} \Id}
        \\
        {\sHfourlt} \ar[rr]^{\displaystyle V \stensor{\hbeta \rt
            \beta}{\alpha} \Id}&& {\sHfourrt} \ar[u]_{\displaystyle
          \Sigma_{\leg{23}}} }
    \end{gathered}
  \end{gather}
  where $\Sigma_{\leg{23}}$ denotes the isomorphism
  \begin{align*}
    \sHfourrt \cong (H {_{\rho_{\alpha}} \tl} \beta) {_{\rho_{\hbeta \lt
          \beta}} \tl} \alpha &\mycong (H {_{\rho_{\hbeta}} \tl} \alpha)
    {_{\rho_{\alpha \lt \alpha}} \tl} \beta \cong \sHfour, \\
    (\zeta \tl \xi) \tl \eta &\mapsto (\zeta \tl \eta) \tl \xi.
  \end{align*}
  Given a $C^{*}$-pseudo-multiplicative unitary
  $(H,\hbeta,\alpha,\beta,V)$, we adopt the following leg notation.
  We abbreviate the operators $V \stensor{\hbeta \rt \hbeta}{\alpha}
  \Id$ and $V \stensor{\alpha \lt \alpha}{\beta} \Id$ by
  $V_{\leg{12}}$, the operators $\Id \stensor{\alpha}{\beta \lt
    \alpha} V$ and $\Id \stensor{\hbeta}{\alpha \lt \alpha} V$ by
  $V_{23}$, and $(\Id \stensor{\hbeta}{\alpha \rt \alpha}
  \Sigma)V_{\leg{12}}\Sigma_{\leg{23}}$ by $V_{13}$. Thus, the
  indices indicate those positions in a relative tensor product
  where the operator acts like $V$.
\end{definition}

Let $(H,\hbeta,\alpha,\beta,V)$ be a $C^{*}$-pseudo-multiplicative
unitary.  We put
\begin{align*}
  \hA(V) &:= \lnspan \bbeta{2} V \kalpha{2}\rnspan \subseteq {\cal
    L}(H), & A(V) &:= \lnspan\balpha{1} V \khbeta{1}\rnspan
  \subseteq {\cal L}(H).
\end{align*}
These spaces satisfy $\hA(V) \subseteq {\cal L}(H_{\beta})$ and
$A(V) \subseteq {\cal L}(H_{\hbeta})$, so that we can define maps
\begin{align*}
  \widehat{\Delta}_{V} &\colon \hA \to {\cal L}\big(\Hsource\big), \
  \ha \mapsto V^{*}(1 \stensor{\alpha}{\beta} \ha)V, & \Delta_{V}
  &\colon A \to {\cal L}\big(\Hrange), \ a \mapsto V(a
  \stensor{\hbeta}{\alpha} 1)V^{*}.
\end{align*}
\begin{definition}[{\cite{timmer:cpmu}}]  
  We call a $C^{*}$-pseudo-multiplicative unitary
  $(H,\hbeta,\alpha,\beta,V)$ {\em regular} if
  $[\balpha{1}V\kalpha{2}] = [\alpha\alpha^{*}]$, and {\em
    well-behaved} if $(\hA(V)^{\hbeta,\alpha}_{H},\hDelta_{V})$ and
  $(A(V)_{H}^{\alpha,\beta},\Delta_{V})$ are Hopf $C^{*}$-bimodules
  over $\mu^{op}$ and $\mu$, respectively.
\end{definition}
\begin{theorem}[{\cite{timmer:cpmu}}] \label{theorem:pmu-legs}
 Every regular
  $C^{*}$-pseudo-multiplicative unitary is well-behaved. 
  \end{theorem}

  Let $(H,\hbeta,\alpha,\beta,V)$ be a $C^{*}$-pseudo-multiplicative
  unitary.  We put
\begin{align*}
 V^{op} :=\Sigma V^{*} \Sigma \colon H \stensor{\beta}{\alpha} H
\xrightarrow{\Sigma} \sHrange \xrightarrow{V^{*}} \sHsource
\xrightarrow{\Sigma} H \stensor{\alpha}{\hbeta} H. 
\end{align*}
Then $(H,\beta,\alpha,\hat\beta,V^{op})$ is a
$C^{*}$-pseudo-multiplicative unitary over $\mu^{op}$, called the
{\em opposite} of $(H,\hat\beta,\alpha,\beta,V)$ \cite[Remark
4.3]{timmer:cpmu}. One easily checks that $V^{op}$ is regular if $V$
is regular.

\paragraph{The fundamental unitary of a compact $C^{*}$-quantum
  groupoid}
Throughout this section, let
$(B,\mu,A,r,\phi,s,\psi,\delta,R,\Delta)$ be a compact
$C^{*}$-quantum groupoid. We use the same notation as in the
preceding section.

 The main result of this paragraph is the
following theorem.
\begin{theorem} \label{theorem:pmu} There exists a regular
  $C^{*}$-pseudo-multiplicative unitary $(H,\hbeta,\alpha,\beta,V)$
  such that $V |a\zeta_{\psi}\rangle_{1} =
  \Delta(a)|\zeta_{\psi}\rangle_{1}$ for all $a \in A$.
\end{theorem}
We prove this result in several steps. Til the end of this section,
we fix a compact $C^{*}$-quantum groupoid
$(B,\mu,A,r,\phi,s,\psi,\delta,R,\Delta)$ and use
the same notation as in Section \ref{section:cqg}. 
\begin{proposition} \label{proposition:pmu} 
  \begin{enumerate}\item 
    There exists a unique unitary $V \colon \sHsource \to \sHrange$
    such that $V |a\zeta_{\psi}\rangle_{1} =
    \Delta(a)|\zeta_{\psi}\rangle_{1}$ for all $a \in A$.
  \item $V(a\zeta_{\nu} \tl d^{op}\zeta_{\psi}) =
    \Delta(a)(\zeta_{\nu} \tl d^{op}\zeta_{\phi})$ for all $a,d \in
    A$.
  \item $V(\hbeta \rt \beta) = \hbeta \lt \beta$, $V(\hbeta \rt
    \hbeta) = \alpha \rt \hbeta$, $V(\hbeta \rt \halpha) = \alpha
    \rt \halpha$, $V(\halpha \lt \alpha) = \halpha \lt \beta$.
  \end{enumerate}
  \end{proposition}
  \begin{proof}
    i) Let $a \in A$, $\eta \in \beta$, $\zeta \in H_{\mu}$.  Since
    $\psi$ is a bounded right Haar weight for
    $(A^{\alpha,\beta}_{H},\Delta)$,
  \begin{align*}
    \big \langle \Delta(a)(\zeta_{\psi} \tr \zeta \tl \eta) \big|
    \langle \Delta(a)(\zeta_{\psi} \tr \zeta \tl \eta)
    \big\rangle_{(\sHrange)} &= \big\langle \zeta\big|\zeta_{\psi}^{*}
    \langle\eta|_{2}\Delta(a^{*}a)|\eta\rangle_{2}\zeta_{\psi}\zeta\rangle
    \\
    &= \big\langle \zeta\big|
    \eta^{*}\rho_{\alpha}(\zeta_{\psi}^{*}a^{*}a\zeta_{\psi})\eta
    \zeta\big\rangle \\ &= \big\langle a\zeta_{\psi} \tr
    \eta\zeta|a\zeta_{\psi} \tr \eta\zeta\rangle_{(\sHsource)}.
  \end{align*}
  Therefore, there exists an isometry $V \colon \sHsource \to \sHrange$ such
  that $V |a\zeta_{\psi}\rangle_{1} =
  \Delta(a)|\zeta_{\psi}\rangle_{1}$ for all $a \in A$. 
  Since $[\Delta(A)\kbeta{2}]=[\kbeta{2}A]$,
    \begin{align*}
      V(\hbeta \lt \beta) =
      [V|A\zeta_{\psi}\rangle_{1}\beta] &=
      [\Delta(A)|\zeta_{\psi}\rangle_{1}\beta] =
      [\Delta(A) \kbeta{2}\zeta_{\psi}] =
      [\kbeta{2}A\zeta_{\psi}] = \hbeta \lt \beta.
    \end{align*}
    In particular, $V$ is surjective and hence a unitary.

    ii) By Proposition \ref{proposition:cqg-modular}, we have for
    all $a,d \in A$
    \begin{align*}
      V(a\zeta_{\nu} \tl d^{op}\zeta_{\psi}) &= V(a \delta^{-1/2}
      \zeta_{\psi}\zeta_{\mu} \tl
      d^{op}\zeta_{\psi}) \\
      &= V(a \delta^{-1/2}\zeta_{\psi} \tr d^{op} \zeta_{\nu^{-1}})
      \\
      &= \Delta(a\delta^{-1/2})(\zeta_{\psi} \tr d^{op}\delta^{1/2}
      \zeta_{\nu}) \\
      &= \Delta(a)(\delta^{-1/2} \stensor{\alpha}{\beta}
      \delta^{-1/2}) (\zeta_{\nu^{-1}} \tl
      d^{op}\delta^{1/2}\zeta_{\phi}) = \Delta(a)(\zeta_{\nu} \tl
      d^{op}\zeta_{\phi}).
    \end{align*}

    iii) The first relation was already proven above. Since $[\Delta(A)
    |\zeta_{\psi}\rangle_{1}A]=[\kalpha{1}A]$,
    \begin{align*}
      V(\hbeta \rt \hbeta) = [V\khbeta{1}A\zeta_{\psi}] &=
      [V|A\zeta_{\psi}\rangle_{1}A\zeta_{\psi}] =
      [\Delta(A)|\zeta_{\psi}\rangle_{1}A\zeta_{\psi}]
      = [\kalpha{1}A\zeta_{\psi}] = \alpha \rt \hbeta
    \end{align*}
    and similarly $V(\hbeta \rt \halpha)= \alpha \rt \halpha$.
    Finally,  by ii), for all $b \in B$ and $a,d \in A$
    \begin{align*}
      V|a^{op}\zeta_{\psi}\rangle_{2} d\zeta_{\phi} b\zeta_{\mu}
      = V\big(d r(b)\zeta_{\nu} \tl
      a^{op}\zeta_{\psi}\big) 
      &= \Delta(dr(b))(\zeta_{\nu} \tl a^{op}\zeta_{\phi}) \\
      &=
      \Delta(d)\big(r(b)\zeta_{\nu} \tl
      a^{op}\zeta_{\phi}\big)  
      = \Delta(d)|a^{op}\zeta_{\phi}\rangle_{2}
      \zeta_{\phi} b\zeta_{\mu}
    \end{align*}
    and hence $[V\khalpha{2}\alpha] =
    [V|A^{op}\zeta_{\psi}\rangle_{2} A\zeta_{\phi}] =
    [\Delta(A)|A^{op}\zeta_{\phi}\rangle_{2}\zeta_{\phi}] =
    [\Delta(A)\kbeta{2}\zeta_{\phi}] = [\kbeta{2}A\zeta_{\phi}] =
    \halpha \lt \beta$.
\end{proof}

The strong invariance condition on the coinvolution yields the
following important inversion formula for the unitary $V$
constructed above.
\begin{theorem} \label{theorem:pmu-inversion} $V^{*}=(J
  \rtensor{\alpha}{J_{\mu}}{\beta} I)V(J
  \rtensor{\alpha}{J_{\mu}}{\beta} I)$.
  \end{theorem}
  \begin{proof}
    Put $\tilde V:=(J \rtensor{\alpha}{J_{\mu}}{\beta} I)V(J
    \rtensor{\alpha}{J_{\mu}}{\beta} I)$. Then for all  $a,b,c,d \in A$
    \begin{align*}
      \big\langle a\zeta_{\psi} \tr b^{op}\zeta_{\nu^{-1}} \big|
      V^{*}(c^{op}\zeta_{\nu^{-1}} \tl d^{op}\zeta_{\phi})
      \big\rangle &= \big\langle \Delta(a)(\zeta_{\psi} \tr
      b^{op}\zeta_{\nu^{-1}}) \big| c^{op}\zeta_{\psi} \tr
      d^{op}\zeta_{\nu}\big\rangle \\
      &= \big\langle \zeta_{\psi} \tr \zeta_{\nu^{-1}}\big|
      \Delta(a^{*})(c^{op} \stensor{\alpha}{\beta} 1)( \zeta_{\psi} \tr
      (b^{op})^{*}d^{op}\zeta_{\nu} )\big\rangle \\
      &= \big\langle \zeta_{\nu^{-1}}\big| \langle\zeta_{\psi}|_{1}
      \Delta(a^{*})(c^{op} \stensor{\alpha}{\beta} 1)|\zeta_{\psi}\rangle_{1}
      (db^{*})^{op}\zeta_{\nu}\big\rangle, \\
      \big\langle a\zeta_{\psi} \tr b^{op}\zeta_{\nu^{-1}}
      \big|\tilde V(c^{op}\zeta_{\nu^{-1}} \tl d^{op}\zeta_{\phi})
      \big\rangle &= \overline{\big\langle (a^{*})^{op} \zeta_{\psi} \tr
        Ib^{op}\zeta_{\nu^{-1}} \big|
        V(c^{*}\zeta_{\nu^{-1}} \tl Id^{op}\zeta_{\phi})\big\rangle} \\
      &= \overline{\big\langle  \zeta_{\psi} \tr
        Ib^{op}\zeta_{\nu^{-1}} \big| (a^{op} \stensor{\alpha}{\beta}
        1)\Delta(c^{*})(\zeta_{\psi} \tr
        Id^{op}\zeta_{\nu})\big\rangle} \\
      &= \overline{\big\langle  \zeta_{\psi} \tr
        I\zeta_{\nu^{-1}} \big| (a^{op} \stensor{\alpha}{\beta}
        1)\Delta(c^{*})(\zeta_{\psi} \tr
        I(b^{op})^{*}d^{op}\zeta_{\nu})\big\rangle} \\
      &= \overline{\big\langle I\zeta_{\nu^{-1}}\big|
        \langle\zeta_{\psi}|_{1} (a^{op} \stensor{\alpha}{\beta}
        1)\Delta(c^{*})|\zeta_{\psi}\rangle_{1} I
        (db^{*})^{op}\zeta_{\nu}\big\rangle}  \\
      &= \big\langle \zeta_{\nu^{-1}}\big| I
      \langle\zeta_{\psi}|_{1} (a^{op} \stensor{\alpha}{\beta}
      1)\Delta(c^{*})|\zeta_{\psi}\rangle_{1} I
      (db^{*})^{op}\zeta_{\nu}\big\rangle.
      \end{align*}
      Now, the claim follows from condition iii) in Definition
      \ref{definition:cqg}.
    \end{proof}

  \begin{proof}[Proof of Theorem \ref{theorem:pmu}]
    By Lemma \ref{lemma:modules-antiunitary} ii) and Propositions
    \ref{proposition:module-graph} iii) and 
    \ref{proposition:pmu} iii), left multiplication by $(J
    \rtensor{\hbeta}{J_{\mu}}{\alpha} I)V^{*}(J
    \rtensor{\hbeta}{J_{\mu}}{\alpha} I)$ acts on subspaces of
    ${\cal L}(H_{\mu}, \sHsource)$
    below as follows:
    \begin{gather*}
      [\kalpha{2}\alpha] \xrightarrow{(J
        \rtensor{\hbeta}{J_{\mu}}{\alpha} I)} [\kbeta{2} J \alpha] =
      \kbeta{2}\hbeta J_{\mu} \xrightarrow{V^{*}} [\khbeta{1}\beta
      J_{\mu}] \xrightarrow{(J \rtensor{\hbeta}{J_{\mu}}{\alpha} I)}
      [\kalpha{1} I \beta J_{\mu}]
      = [\kalpha{1}\alpha], \\
      [\kalpha{2}\beta]  \xrightarrow{(J
        \rtensor{\hbeta}{J_{\mu}}{\alpha} I)}
      [\kbeta{2} J \beta] =  [\kbeta{2} \halpha J_{\mu}]
      \xrightarrow{V^{*}}
      [\kalpha{2}\halpha J_{\mu}]  \xrightarrow{(J
        \rtensor{\hbeta}{J_{\mu}}{\alpha} I)}
      [\kbeta{2} J \halpha J_{\mu}] = [\kbeta{2}\beta], \\
      [\khbeta{1}\hbeta]  \xrightarrow{(J
        \rtensor{\hbeta}{J_{\mu}}{\alpha} I)}
      [\kalpha{1} I \hbeta] =
      [\kalpha{1} \halpha J_{\mu}] \xrightarrow{V^{*}}
      [\khbeta{1} \halpha J_{\mu}]  \xrightarrow{(J
        \rtensor{\hbeta}{J_{\mu}}{\alpha} I)}
      [\kalpha{1} I \halpha J_{\mu}] = [\kalpha{1}\hbeta].
    \end{gather*}
    Now, Theorem \ref{theorem:pmu-inversion} implies $V(\alpha \lt
    \alpha) = \alpha \rt \alpha$, $V(\beta \lt \alpha) = \beta \lt
    \beta$, $V(\hbeta \rt \hbeta) = \alpha \rt \hbeta$. These relations
    and the relations in Proposition \ref{proposition:pmu} iii) are
    precisely \eqref{eq:pmu-intertwine}.

    Let us show that diagram \eqref{eq:pmu-pentagon} commutes.  Let
    $a,d \in A$ and $\omega \in H$. Then
    \begin{gather*}
      V_{23}V_{12}(a \zeta_{\psi} \tr d \zeta_{\psi} \tr \omega) = V_{23}
      (\Delta(a) \stensor{\alpha \rt \hbeta}{\alpha} \Id) (\zeta_{\psi} \tr d\zeta_{\psi} \tr \omega) 
      = \Delta^{(2)}(a) \big(\zeta_{\psi} \tr \Delta(d)(\zeta_{\psi} \tr
      \omega)\big),
    \end{gather*}
    where $\Delta^{(2)}=(\Delta \ast \Id) \circ \Delta = (\Id \ast
    \Delta) \circ \Delta$.  On the other hand,
    \begin{align*}
      V_{12}V_{13}V_{23}(a\zeta_{\psi} \tr d\zeta_{\psi} \tr \omega) &=
      V_{12}V_{13} \big(a \zeta_{\psi} \tr \Delta(d)(\zeta_{\psi} \tr
      \omega)\big)  \\
      &= V_{12} \Delta_{13}(a)\big(\zeta_{\psi} \tr \Delta(d)(\zeta_{\psi}
      \tr \omega) \big) 
      = \Delta^{(2)}(a) \big(\zeta_{\psi} \tr \Delta(d)(\zeta_{\psi} \tr
      \omega))\big),
    \end{align*}
    where $\Delta_{13}(a)=\Sigma_{23}(\Delta(a) \stensor{\hbeta \lt
      \beta}{\alpha} \Id)\Sigma_{23}$. Since $a$ and $d$ were
    arbitrary, we can conclude $V_{23}V_{12}=V_{12}V_{13}V_{23}$.

    Finally, $V$ is regular because by Theorem
    \ref{theorem:pmu-inversion}, Lemma
    \ref{lemma:modules-antiunitary} ii) and Proposition
    \ref{proposition:module-graph} iii),
    \begin{align*}
      [\balpha{1} V \kalpha{2}] &=
      [\balpha{1}  (J
        \rtensor{\hbeta}{J_{\mu}}{\alpha} I)V^{*}(J
        \rtensor{\hbeta}{J_{\mu}}{\alpha} I) \kalpha{2}] \\ &=
        [I \bhbeta{1}V^{*}\kbeta{2}J]  \\
        &=  [I \langle \zeta_{\psi}|_{1} \Delta(A) \kbeta{2} J]
        \\
        &= [I \langle \zeta_{\psi}|_{1} \kbeta{2} AJ] 
        = [I \beta  J_{\mu} \cdot J_{\mu} \zeta_{\psi}^{*} A J]
        =[\alpha\alpha^{*}]. \qedhere
    \end{align*}
  \end{proof}

  By Theorem \ref{theorem:pmu-legs}, the regular
  $C^{*}$-pseudo-multiplicative unitary $(H,\hbeta,\alpha,\beta,V)$
  constructed above yields two Hopf $C^{*}$-bimodules
  $(A(V)^{\alpha,\beta}_{H},\Delta_{V})$ and
  $(\hA(V)^{\hbeta,\alpha}_{H},\hDelta_{V})$. 
  \begin{proposition} \label{proposition:pmu-right-leg}
    $(A(V)^{\alpha,\beta}_{H},\Delta_{V})=(A^{\alpha,\beta}_{H},\Delta)$.
  \end{proposition}
  \begin{proof}
    We have $A(V) = [\balpha{1}V\khbeta{1}] =
    [\balpha{1}\Delta(A)|\zeta_{\psi}\rangle_{1}] = [A \balpha{1}
    |\zeta_{\psi}\rangle_{1}] = [A \rho_{\alpha}(\alpha^{*}
    \zeta_{\psi})] = [A s(B^{op})] = A$ and $ \Delta_{V}(a) = V(a
    \stensor{\hbeta}{\alpha} 1) V^{*} = \Delta(a)$ for all $a \in
    A$.
  \end{proof}
  The Hopf $C^{*}$-bimodule
  $(\hA(V)^{\hbeta,\alpha}_{H},\hDelta_{V})$ will be studied in the
  next section.

Our first application of the fundamental unitary is to prove
that the coinvolution reverses the comultiplication. 
\begin{theorem} \label{theorem:pmu-coinvolution}
 $(R \rfibre{\alpha}{\beta} R)
  \circ \Delta = \Ad_{\Sigma} \circ \Delta \circ R$.
  \end{theorem}
The proof involves  the following formulas:
\begin{lemma} \label{lemma:pmu-coinvolution}
  \begin{enumerate}
  \item $\Delta(\langle \xi|_{1}V|\xi'\rangle_{1}) = \langle
    \xi|_{1}V_{12}V_{13}|\xi'\rangle_{1}$  for all $\xi \in
    \alpha, \xi' \in \hbeta$.
  \item 
    $R(\langle \xi|_{1}V|\xi'\rangle_{1})= \langle J\xi'
    J_{\mu}|_{1} V |J \xi J_{\mu}\rangle_{1}$ for all $\xi \in
    \alpha, \xi' \in \hbeta$.
  \end{enumerate}
\end{lemma}
\begin{proof}
  i) For all $\xi \in \alpha$, $\xi' \in \hbeta$,  we have that
  $\Delta(\langle \xi|_{1}V|\xi'\rangle_{1}) = V((\langle \xi|_{1} V
  |\xi\rangle_{1}) \stensor{\hbeta}{\alpha} 1)V^{*} = \langle
  \xi|_{1} V_{23}^{*}V_{12}V_{23}|\xi'\rangle_{1} = \langle
  \xi|_{1}V_{12}V_{13} |\xi'\rangle_{1}$; see also 
  \cite[Lemma 4.13]{timmer:cpmu}.
  
  ii)  By Lemma \ref{lemma:modules-antiunitary} and Theorem
  \ref{theorem:pmu-inversion}, we have that $R(\langle
  \xi|_{1}V|\xi'\rangle_{1}) = I \langle \xi'|_{1}
  V^{*}|\xi\rangle_{1} I = \langle J\xi' J_{\mu}|_{1} (J
  \rtensor{\alpha}{J_{\mu}}{\beta} I)^{*} V^{*} (J
  \rtensor{\alpha}{J_{\mu}}{\beta} I)^{*} |J \xi J_{\mu}\rangle_{1} =
  \langle J\xi' J_{\mu}|_{1} V |J \xi J_{\mu}\rangle_{1}$  for all
  $\xi \in \alpha$, $\xi' \in \hbeta$
\end{proof}
  \begin{proof}[Proof of Theorem \ref{theorem:pmu-coinvolution}]
    Let $\xi \in \alpha$ and $\xi' \in \hbeta$. By Lemma
    \ref{lemma:pmu-coinvolution} i),
    \begin{align*}
      (\Ad_{\Sigma} \circ (R \rfibre{\alpha}{\beta} R) \circ
      \Delta)\big(\langle \xi|_{1}V|\xi'\rangle_{1}\big) &=
      (\Ad_{\Sigma} \circ (R \rfibre{\alpha}{\beta}
      R))\big(\langle \xi|_{1}V_{12}V_{13}|\xi'\rangle_{1}\big)  \\
      &= \Ad_{\Sigma}\big((I \rtensor{\alpha}{J_{\mu}}{\beta} I)^{*} \langle
      \xi'|_{1} V_{13}^{*} V_{12}^{*} |\xi\rangle_{1} (I
      \rtensor{\alpha}{J_{\mu}}{\beta} I) \big).
    \end{align*}
    By Lemma \ref{lemma:modules-antiunitary} ii), we can rewrite
    this expression in the form
    \begin{align*}
      \Ad_{\Sigma}\big(\langle J\xi' J_{\mu}|_{1} 
      (J \rtensor{\hbeta}{J_{\mu}}{\alpha \rt \alpha} (I
      \rtensor{\alpha}{J_{\mu}}{\beta} I))
      V_{13}^{*}V_{12}^{*} (J \rtensor{\alpha}{J_{\mu}}{\beta} I
      \rtensor{\alpha}{J_{\mu}}{\beta} I)^{*}
      |J\xi J_{\mu}\rangle_{1}
      \big).
    \end{align*}
    Two applications of Lemma \ref{lemma:modules-antiunitary} iii)
    and Theorem \ref{theorem:pmu-inversion} and an application of
    Lemma \ref{lemma:pmu-coinvolution} ii) show that this expression
    is equal to
    \begin{align*}
      \Ad_{\Sigma}\big(\langle J\xi' J_{\mu}|_{1}
      V_{13}V_{12}|J\xi J_{\mu}\rangle_{1} \big) &= \langle J \xi'
      J_{\mu}|_{1} V_{12}V_{13} |J\xi J_{\mu}\rangle_{1} \\
      &= \Delta\big(
      \langle J \xi' J_{\mu}|_{1} V |J\xi J_{\mu}\rangle_{1} \big) 
      = \Delta(R(\langle \xi|_{1}V|\xi'\rangle_{1})). \qedhere
    \end{align*}
  \end{proof}

\paragraph{A second fundamental unitary}
Like in the theory of locally compact quantum groups, we can associate
to a given compact $C^{*}$-quantum groupoid besides
$(H,\hbeta,\alpha,\beta,V)$ a second $C^{*}$-pseudo-multiplicative
unitary $(H,\halpha,\beta,\alpha,W)$ as follows.
\begin{theorem} \label{theorem:pmu-second} There exists a regular
  $C^{*}$-pseudo-multiplicative unitary $(H,\alpha,\beta,\halpha,W)$
  such that $W^{*} |a\zeta_{\phi}\rangle_{2} =
  \Delta(a)|\zeta_{\phi}\rangle_{2}$ for all $a \in A$. Moreover,
  \begin{align*}
    W=\Sigma(I \rtensor{\hbeta}{J_{\mu}}{\alpha} I) V^{*} (I
    \rtensor{\beta}{J_{\mu}}{\alpha} I) \Sigma = (I
    \rtensor{\alpha}{J_{\mu}}{\hbeta} I)V^{op}(I
    \rtensor{\alpha}{J_{\mu}}{\beta} I).
  \end{align*}
\end{theorem}
\begin{proof}
  Let $a \in A$ and $\xi \in H$. Since
  $I\zeta_{\psi}J_{\mu}=\zeta_{\phi}$ and $\Delta(R(a)^{*})= 
  \Sigma (I \rtensor{\alpha}{J_{\mu}}{\beta} I)\Delta(a)(I
  \rtensor{\beta}{J_{\mu}}{\alpha} I)\Sigma$, 
  \begin{align*}
    \Sigma(I \rtensor{\beta}{J_{\mu}}{\alpha} I)^{*} V (I
    \rtensor{\hbeta}{J_{\mu}}{\alpha} I)^{*}\Sigma (\xi \tl
    a\zeta_{\phi}) &= \Sigma(I \rtensor{\alpha}{J_{\mu}}{\beta} I)V
    (Ia\zeta_{\phi}J_{\mu} \tr I\xi)  \\
    &= \Sigma(I \rtensor{\alpha}{J_{\mu}}{\beta} I)V
    (R(a)^{*}\zeta_{\psi} \tr I\xi)  \\
    &= \Sigma(I \rtensor{\alpha}{J_{\mu}}{\beta} I)\Delta(R(a)^{*})
    (\zeta_{\psi} \tr I\xi)  \\
    &= \Delta(a)(I \rtensor{\alpha}{J_{\mu}}{\beta} I) \Sigma
    (\zeta_{\psi} \tr I\xi) = \Delta(a)(\xi \tl \zeta_{\phi}).
  \end{align*}
  Therefore, the unitary $W=\Sigma(I
  \rtensor{\hbeta}{J_{\mu}}{\alpha} I) V^{*} (I
  \rtensor{\beta}{J_{\mu}}{\alpha} I) \Sigma$ satisfies
  $W^{*}|a\zeta_{\phi}\rangle_{2} =
  \Delta(a)|\zeta_{\phi}\rangle_{2}$ for all $a \in A$.  Since
  $(H,\beta,\alpha,\hbeta,V^{op})$ is a regular
  $C^{*}$-pseudo-multiplicative unitary, so is
  $(H,\alpha,\beta,\halpha,W)$.
\end{proof}

\paragraph{The passage to the setting of von Neumann algebras}
In this paragraph, we indicate how every compact $C^{*}$-quantum
groupoid can be completed to a measurable quantum groupoid in the
sense of Lesieur \cite{lesieur} and Enock \cite{enock:lesieur}. We
assume some familiarity with \cite{lesieur} or \cite{enock:lesieur}.

Let $\mu$ be a faithful KMS-state on a unital $C^{*}$-algebra
$B$. Then the state $\tilde \mu$ on $N:=B'' \subseteq
{\cal L}(H_{\mu})$ given by $y \mapsto
\langle\zeta_{\mu}|y\zeta_{\mu}\rangle$ is the unique normal
extension of $\mu$ and is faithful because $\zeta_{\mu}$ is cyclic
for $\pi_{\mu^{op}}(B^{op}) \subseteq N'$. Evidently, the Hilbert
space $H_{\tilde \mu}:=H_{\mu}$ and the map $\Lambda_{\tilde \mu}
\colon N \to H_{\tilde \mu}$, $y \mapsto y\zeta_{\mu}$, form a
GNS-representation for $\tilde \mu$.  
\begin{lemma}
  Let $\mu$ be a faithful KMS-state on a unital $C^{*}$-algebra $B$
  and let $(r,\phi)$ be a $\mu$-module structure on a unital
  $C^{*}$-algebra $A$. We put $N:=B'' \subseteq {\cal L}(H_{\mu})$,
  $M:=A'' \subseteq {\cal L}(H_{\nu})$,  and use the notation of
  Lemma \ref{lemma:graph-rieffel}.
  \begin{enumerate}
  \item $r$ extends uniquely to a normal embedding $\tilde r
    \colon N \to M$.
  \item $\phi$ extends uniquely to a normal completely positive map
    $\tilde \phi \colon M \to N$, and $\tilde \nu = \tilde \mu \circ
    \tilde \phi$. If $\phi$ is faithful, so is $\tilde \phi$.
  \item $\zeta y = \tilde r(y)\zeta$, $\zeta^{*}x = \tilde
    \phi(x)\zeta^{*}$, $\tilde \phi(x\tilde r(y))=\tilde
    \phi(x)\tilde r(y)$ for all $x \in M$, $y \in N$.
  \end{enumerate}
\end{lemma}
\begin{proof}
  i) Uniqueness is clear. Put $\halpha := [A\zeta]$. By Lemmas
  \ref{lemma:modules} and \ref{lemma:module}, $H_{\nu} \cong \halpha
  \tr H_{\mu}$. Hence, we can define a $*$-homomorphism $\tilde r
  \colon N \to {\cal L}(\alpha \tr H_{\mu}) \cong {\cal
    L}(H_{\nu})$ by $x \mapsto \Id_{\halpha} \tr x$. By Lemma
  \ref{lemma:module} i), $\tilde r$ extends $r$, and routine
  arguments show that $\tilde r$  is normal and injective.

  ii) $\tilde \phi$ is uniquely determined by $\tilde \phi(x)=
  \zeta^{*}x\zeta$ for all $x \in M$, and clearly $\tilde \nu(x) =
  \langle \zeta_{\nu} |x\zeta_{\nu}\rangle = \langle
  \zeta_{\mu}|\zeta^{*}x\zeta\zeta_{\mu}\rangle = (\tilde \mu \circ
  \tilde \phi)(x)$ for all $x \in M$. If $\phi$ is faithful, so is
  $\nu$ and, since $\tilde \mu$ is faithful and $\tilde \nu=\tilde
  \mu \circ \tilde \phi$, also $\tilde \phi$ is faithful.

  iii) Use Lemma \ref{lemma:graph-rieffel} iii) and
  the fact that $\tilde r$, $\tilde \phi$ are normal extensions of $r,\phi$.
\end{proof}

Let $(B,\mu,A,r,\phi,s,\psi,\delta,R,\Delta)$ be a compact
$C^{*}$-quantum groupoid. We keep the notation introduced before and
put
\begin{align*}
    N&:= B'' \subseteq {\cal L}(H_{\mu}), & N^{op} &:=(B^{op})'' = N'
    \subseteq {\cal L}(H_{\mu}), &
    M&:= A'' \subseteq {\cal L}(H).
\end{align*}
By the previous remarks, the maps $\mu$, $r$, $\phi$, $s$, $\psi$
have unique normal extensions
\begin{gather*}
  \begin{aligned}
    \tilde \mu &\colon N \to \complex, & \tilde r &\colon N \to M, &
    \tilde \phi &\colon M \to N, & \tilde s &\colon N' \to M, &
    \tilde \psi &\colon M \to N'.
  \end{aligned} 
\end{gather*}

Before we can extend the comultiplication $\Delta$ from $A$ to $M$,
we need to recall the definition of the fiber product of von Neumann
algebras \cite{sauvageot:2} and the underlying relative tensor
product of Hilbert spaces \cite{takesaki:2}; a reference is also
\cite[\S 10]{timmer:qg}.  The relative tensor product of $H$ with
itself, taken with respect to $\tilde
s,\tilde r$ and $\tilde \mu$, is defined as follows. Put
\begin{align*}
  D(H_{\tilde r};\tilde \mu) := \big\{ \eta \in H \,\big|\, \exists C>0
  \forall y \in N: \|\tilde r(y) \eta\| \leq C \| y\zeta_{\mu}\| \big\}.
\end{align*}
Evidently, an element $\eta \in H$ belongs to $D(H_{\tilde r};
\tilde \mu)$ if and only if the map $N\zeta_{\mu} \to H$ given by
$y\zeta_{\mu} \mapsto \tilde r(y)\eta$ extends to a bounded linear
map $L(\eta) \colon H_{\mu} \to H$, and $L(\eta)^{*}L(\eta') \in N'$
for all $\eta,\eta' \in D(H_{\tilde r};\tilde \mu)$. The relative tensor product $H \rtensor{\tilde s}{\tilde \mu}{\tilde
  r} H$ is the separated completion of the algebraic tensor product
$H \odot D(H_{\tilde r};\tilde \mu)$ with respect to the
sesquilinear form defined by
\begin{align*}
 \langle \omega \odot \eta|\omega' \odot  \eta'\rangle = \langle
 \omega | \tilde s(L(\eta)^{*}L(\eta'))\omega'\rangle \quad
 \text{for all } \omega,\omega' \in H, \, \eta,\eta' \in
 D(H_{\tilde r};\tilde \mu).
\end{align*}
We denote the image of an element $\omega \odot \eta$ in $H
\rtensor{\tilde s}{\tilde \mu}{\tilde r} H$ by $\omega
\stensor{\tilde s}{\tilde r} \eta$.
\begin{lemma}
  \begin{enumerate}
  \item $a^{op}\zeta_{\nu} \in D(H_{\tilde r};\tilde \mu)$ and
    $L(a^{op}\zeta_{\nu}) = a^{op}\zeta_{\phi} \in \beta$ for all $a
    \in A$.
  \item There exist inverse isomorphisms
    \begin{align*}
      \Phi_{\alpha,\beta} &\colon H \stensor{\alpha}{\beta} H \cong
      H {_{\rho_{\alpha}} \tl} \beta \to H \rtensor{\tilde s}{\tilde
        \mu}{\tilde r} H, & \Psi_{\alpha,\beta} &\colon H
      \rtensor{\tilde s}{\tilde \mu}{\tilde r} H \to \alpha
      \tr_{\rho_{\beta}} H \cong H \stensor{\alpha}{\beta} H
    \end{align*}
    such that for all $\omega \in H$, $a \in A$, $\xi \in \alpha$, $\eta \in
    D(H_{\tilde r},\tilde \mu)$, $\zeta \in H_{\mu}$
    \begin{align*}
      \Phi_{\alpha,\beta}(\omega \tl a^{op}\zeta_{\phi}) &= \omega
      \stensor{\tilde s}{\tilde r} a^{op}\zeta_{\nu}, &
      \Psi_{\alpha,\beta}(\xi\zeta \stensor{\tilde s}{\tilde r}
      \eta) &= \xi \tr L(\eta)\zeta.
    \end{align*}
  \end{enumerate}
\end{lemma}
\begin{proof}
  i) For all $a \in A$, $y \in N$,  $\tilde
  r(y)a^{op}\zeta_{\nu} = a^{op} \tilde r(y)\zeta_{\phi} \zeta_{\mu}
  = a^{op} \zeta_{\phi} y \zeta_{\mu}$. The claims follow.

  ii) The formulas for $\Phi_{\alpha,\beta}$ and
  $\Psi_{\alpha,\beta}$ define isometries because for all $\omega,
  a, \xi, \eta, \zeta$ as above,
  \begin{align*}
\|\omega \tl a^{op}\zeta_{\phi}
\|^{2} &= \big\langle \omega
    \big|\rho_{\alpha}(\zeta_{\phi}^{*}(a^{op})^{*}a^{op}\zeta_{\phi})
    \omega    \big\rangle 
     = \big\langle \omega \big|\tilde
    s(L(a^{op}\zeta_{\nu})^{*}L(a^{op}\zeta_{\nu}))\omega
    \big\rangle = \|\omega \stensor{\tilde s}{\tilde r}
    a^{op}\zeta_{\nu} \|^{2}
  \end{align*}
  and
  \begin{align*}
    \| \xi \zeta \stensor{\tilde s}{\tilde r} \eta \|^{2} = \langle
    \xi\zeta |\tilde s(L(\eta)^{*}L(\eta))\xi \zeta\rangle &=
    \langle \zeta | \xi^{*}\xi L(\eta)^{*}L(\eta) \zeta\rangle \\ &=
    \langle \zeta | L(\eta)^{*} \rho_{\beta}(\xi^{*}\xi^{*})L(\eta)
    \zeta\rangle = \| \xi \tr L(\eta)\zeta\|^{2}.
  \end{align*}
  Moreover, $\Psi_{\alpha,\beta} \circ \Phi_{\alpha,\beta}=\Id$
  because for all $a,\xi,\zeta$ as above,
  \begin{align*}
    (\Psi_{\alpha,\beta} \circ \Phi_{\alpha,\beta})(\xi \zeta \tl
    a^{op}\zeta_{\phi} ) = \xi \tr L(a^{op}\zeta_{\nu})\zeta &= \xi
    \tr a^{op}\zeta_{\phi} \zeta \equiv \xi \zeta \tl
    a^{op}\zeta_{\phi}. \qedhere
  \end{align*}
\end{proof}
 We identify $H \stensor{\alpha}{\beta} H$ with $\HsrH$
via $\Phi_{\alpha,\beta}$ and $\Psi_{\alpha,\beta}$ without further
notice. 

The fiber product $M \fibre{\tilde s}{\tilde \mu}{\tilde r} M$ is
defined as follows. One has $\tilde r(N)' D(H_{\tilde r};\tilde \mu)
\subseteq D(H_{\tilde r};\tilde \mu)$, and for each $x,x' \in M'$,
there exists a well-defined operator $x \rtensor{\tilde s}{\tilde
  \mu}{\tilde r} x' \subseteq {\cal L}(\HsrH)$ such that $(x
\rtensor{\tilde s}{\tilde \mu}{\tilde r} x')(\omega \stensor{\tilde
  s}{\tilde r} \eta) = x\omega \stensor{\tilde s}{\tilde r} x\eta$
for all $\omega \in H$, $\eta \in D(H_{\tilde r};\tilde \mu)$. Now,
\begin{align*}
  M \fibre{\tilde s}{\tilde \mu}{\tilde r} M = (M' \stensor{\tilde s}{\tilde r}
  M')' \subseteq {\cal L}(\HsrH).
\end{align*}
\begin{lemma}
  $\Delta$ extends to a normal $*$-homomorphism $\tilde \Delta
  \colon M \to   M \fibre{\tilde s}{\tilde \mu}{\tilde r} M$.
\end{lemma}
\begin{proof}
  Simply define $\tilde \Delta$ by $\tilde \Delta(x):=V(x \tl
  \Id_{\alpha})V^{*}$ for all $x \in M$.
\end{proof}
The notion of a measurable quantum groupoid was first defined in
\cite{lesieur}; later, the definition was changed in
\cite[\S 6]{enock:lesieur}.
\begin{theorem}
 $(N,M,\tilde r,\tilde s, \tilde\Delta, \tilde
  \phi,\tilde \psi, \tilde \mu)$ is a measurable quantum groupoid.
\end{theorem}
\begin{proof}
  First, one has to check that $(N,M, \tilde r,\tilde s,\tilde
  \Delta)$ is a Hopf-bimodule; this follows from the definition of
  $\tilde \Delta$ and the fact that $V$ is a
  $C^{*}$-pseudo-multiplicative unitary.

  Second, one has to check that $\tilde \phi$ and $\tilde\psi$ are
  left- and right-invariant, respectively. This follows from the
  fact that these maps are normal extensions of $\phi$ and $\psi$,
  which are left- and right-invariant, respectively.

  Finally, one has to check that the modular automorphism groups of
  $\tilde \nu=\tilde \mu \circ \tilde \phi$ and $\tilde
  \nu^{-1}=\tilde \mu^{op} \circ \tilde \psi$ commute, but this
  follows from the fact that $\tilde \nu^{-1}=\tilde
  \nu_{\delta^{1/2}}$.
\end{proof}

\section{Supplements on $C^{*}$-pseudo-multiplicative unitaries}

\label{section:add}

In this section, we interrupt our discussion of compact $C^{*}$-quantum
groupoids and study several properties $C^{*}$-pseudo-multiplicative
unitaries that shall prove useful later.  The corresponding
properties for multiplicative unitaries were introduced and studied
in \cite{baaj:2}.  Throughout this section, let $\mu$ be a faithful
KMS-state on a unital $C^{*}$-algebra $B$.

\paragraph{Fixed and cofixed elements for a
  $C^{*}$-pseudo-multiplicative unitary}
We shall study elements with the following property:
\begin{definition} \label{definition:fixed} Let $\pmu$  be a
  $C^{*}$-pseudo-multiplicative unitary over $\mu$. A {\em fixed
    element} for $V$ is an element $\eta \in \hbeta \cap \alpha$
  satisfying $V|\eta\rangle_{1} = |\eta\rangle_{1} \in {\cal L}(H,
  \sHrange)$. A {\em cofixed element} for $V$ is an element $\xi \in
  \alpha \cap \beta$ satisfying $V|\xi\rangle_{2}=|\xi\rangle_{2}
  \in {\cal L}(H,\sHrange)$. We denote the set of all fixed/cofixed
  elements for $V$ by $\Fix(V)$/$\Cofix(V)$.
\end{definition}
Til the end of this paragraph, let $\pmu$ be a $C^{*}$-pseudo-multiplicative unitary
over $\mu$.
\begin{remarks} \label{remarks:fixed} 
  \begin{enumerate}
  \item  $\Fix(V)=\Cofix(V^{op})$ and
    $\Cofix(V)=\Fix(V^{op})$.
  \item $\Fix(V)^{*} \Fix(V)$ and $\Cofix(V)^{*}\Cofix(V)$ are
    contained in $B \cap B^{op}=Z(B)$.
  \item Since $\Fix(V) \subseteq \hbeta \cap \alpha$, we have
    $\rho_{\alpha}(B^{op})\Fix(V) = \Fix(V)B^{op} \subseteq \hbeta$
    and $\rho_{\hbeta}(B)\Fix(V) =\Fix(V)B \subseteq \alpha$.
    Likewise, $\rho_{\beta}(B)\Cofix(V) \subseteq \alpha$ and
    $\rho_{\alpha}(B^{op})\Cofix(V) \subseteq \beta$.
  \end{enumerate}
\end{remarks}
\begin{lemma} \label{lemma:fixed}
  \begin{enumerate}
  \item $\langle \xi|_{2} V |\xi'\rangle_{2} =
    \rho_{\alpha}(\xi^{*}\xi') = \rho_{\hbeta}(\xi^{*}\xi')$ for all
    $\xi,\xi' \in \Cofix(V)$, and $\langle \eta|_{1} V |\eta'\rangle_{1} =
    \rho_{\beta}(\eta^{*}\eta') = \rho_{\alpha}(\eta^{*}\eta')$ for all
    $\eta,\eta' \in \Fix(V)$.
  \item $\rho_{\hbeta}(B)\Cofix(V) \subseteq \Cofix(V)$ and
    $\rho_{\beta}(B)\Fix(V) \subseteq \Fix(V)$.
    \item $[EE^{*}E]=E$ for $E\in \{\Cofix(V), \Fix(V)\}$.
    \item  $[\Cofix(V)^{*}\Cofix(V)]$ and $[\Fix(V)^{*}\Fix(V)]$ are
      $C^{*}$-subalgebras of $Z(B)$.
    \end{enumerate}
\end{lemma}
\begin{proof}
  We only prove the assertions on $\Cofix(V)$; the other
  assertions follow similarly.

  i) For all $\xi,\xi'
  \in \Cofix(V)$ and $\zeta \in H$,  $\langle \xi|_{2} V
  |\xi'\rangle_{2} \zeta = \langle \xi|_{2}|\xi'\rangle_{2} \zeta =
  \rho_{\alpha}(\xi^{*}\xi')\zeta$ and $(\langle \xi|_{2} V
  |\xi'\rangle_{2})^{*}\zeta = \langle \xi'|_{2} |\xi\rangle_{2}\zeta =
  \rho_{\hbeta}(\xi^{*}\xi')^{*} \zeta$.

  ii) Let $b \in B$ and $\xi \in \Cofix(V)$. Then
  $\rho_{\hbeta}(b)\xi \in \rho_{\hbeta}(B)\beta \cap
  \rho_{\hbeta}(B)\alpha \subseteq \beta \cap \alpha$, and   $V | \rho_{\hbeta}(b)\xi\rangle_{2} = V \rho_{(\hbeta \rt
      \hbeta)}(b)|\xi\rangle_{2} = \rho_{(\alpha \rt
      \hbeta)}(b)V|\xi\rangle_{2} = \rho_{(\alpha \rt \hbeta)}(b)
    |\xi\rangle_{2} = |\rho_{\hbeta}(b)\xi\rangle_{2}$ because
 $V(\hbeta \rt \hbeta)=\alpha \rt \hbeta$.

 iii) Let $\xi,\xi',\xi'' \in \Cofix(V)$. Then
 $\rho_{\alpha}(\xi'{}^{*}\xi'') = \rho_{\hbeta}(\xi'{}^{*}\xi'')$
 by i) and hence $V|\xi \xi'{}^{*}\xi''\rangle_{2} =
 V|\xi\rangle_{2} \rho_{\hbeta}(\xi'{}^{*}\xi'') = |\xi\rangle_{2}
 \rho_{\alpha}(\xi'{}^{*}\xi'') = |\xi\xi'{}^{*}\xi''\rangle_{2}$ in
 ${\cal L}(H,\sHrange)$.

    iv) Immediate from iii).
\end{proof}

\begin{definition} \label{definition:compact} We call $\pmu$ or
  briefly $V$ {\em \'etale} if $\eta^{*}\eta =\Id_{\frakK}$ for some
  $\eta \in \Fix(V)$, and {\em compact} if $\xi^{*}\xi =
  \Id_{\frakK}$ for some $\xi \in \Cofix(V)$.
\end{definition}
\begin{remarks} \label{remarks:compact}
  \begin{enumerate}
  \item By Remark \ref{remarks:fixed}, $V$ is \'etale/compact if and
    only if $V^{op}$ is compact/\'etale.
  \item If $V$ is compact, then $\Id_{H} \in \hA(V)$; if $V$ is
    \'etale, then $\Id_{H} \in A(V)$. This follows directly from
    Lemma \ref{lemma:fixed}.
  \end{enumerate}
\end{remarks}

The following observation supports the plausibility of the
assumptions in condition i) of Definition \ref{definition:cqg}:
\begin{remark} \label{remark:fixed-assumption}
  Let $\pmu$ be a regular $C^{*}$-pseudo-multiplicative unitary over
  $\mu$. If $\xi_{0} \in \Fix(V)$ and $\hbeta=[A(V)\xi_{0}]$, then
by \cite[Lemma 5.8]{timmer:ckac},  \begin{align*}
    [\Delta_{V}(A(V))|\xi_{0}\rangle_{1}A(V)] &= [V(A(V)
    \stensor{\hbeta}{\alpha} 1)V^{*}|\xi_{0}\rangle_{1}A(V)] \\ &=
    [V|A(V)\xi_{0}\rangle_{1}A(V)] = [V\khbeta{1}A(V)] =
    [\kalpha{1}A(V)].
  \end{align*}
  Likewise, if
  $\eta_{0} \in \Cofix(V)$ and $\alpha=[\hA(V)\eta_{0}]$, then
  $[\hDelta_{V}(\hA(V))|\eta_{0}\rangle_{2} \hA(V)] = [\kbeta{2}\hA(V)]$.
\end{remark}

The (co)fixed vectors of the $C^{*}$-pseudo-multiplicative unitaries
introduced in Theorems \ref{theorem:pmu} and
\ref{theorem:pmu-second} are easily determined:
\begin{proposition} \label{proposition:pmu-cqg-fixed} Let
  $(B,\mu,A,r,\phi,s,\psi,\delta,R,\Delta)$ be a compact
  $C^{*}$-quantum groupoid.  
  \begin{enumerate}
  \item The associated $C^{*}$-pseudo-multiplicative
    unitary $(H,\hbeta,\alpha,\beta,V)$ is compact and
    $\Fix(V)=[r(B)\zeta_{\psi}]$.
  \item The associated  $C^{*}$-pseudo-multiplicative
    unitary $(H,\alpha,\beta,\halpha,W)$ is \'etale and
    $\Cofix(W)=[s(B^{op})\zeta_{\phi}]$.
  \end{enumerate}
\end{proposition}
\begin{proof}
  i) Evidently, $\zeta_{\psi} \in \Fix(V)$, and by Lemma
  \ref{lemma:fixed} ii),
  $[r(B)\zeta_{\psi}]=[\rho_{\beta}(B)\zeta_{\psi}] \subseteq
  \Fix(V)$. Conversely, if $\eta_{0} \in \Fix(V)$, then $\eta_{0}
  \in \hbeta = [A\zeta_{\psi}]$ and therefore
  \begin{align*}
    \eta_{0} = \rho_{\alpha}(\zeta_{\phi}^{*}\zeta_{\phi})\eta_{0}
    =
    \langle\zeta_{\phi}|_{2}|\eta_{0}\rangle_{1} \zeta_{\phi} 
    &=    \langle\zeta_{\phi}|_{2}V|\eta_{0}\rangle_{1} \zeta_{\phi}  \\
    & \in
    [    \langle\zeta_{\phi}|_{2}
    \Delta(A)|\zeta_{\psi}\rangle_{1}   \zeta_{\phi} ] \\ 
    &=
    [\langle\zeta_{\phi}|_{2} \Delta(A)|\zeta_{\phi}\rangle_{2}
    \zeta_{\psi}] = [r(\phi(A))\zeta_{\psi}] =
    [r(B)\zeta_{\psi}]. 
  \end{align*}

ii) This follows easily from i) and the relation $W=(I
\rtensor{\alpha}{J_{\mu}}{\hbeta} I)V^{op}(I
\rtensor{\alpha}{J_{\mu}}{\beta} I)$.
\end{proof}

\paragraph{Haar weights and counits obtained from (co)fixed
elements}
Fixed and cofixed elements for a $C^{*}$-pseudo-multiplicative
unitary yield bounded Haar weights and bounded counits on the legs
as follows:
\begin{theorem} \label{theorem:fixed-cofixed}
  Let $\pmu$ be a well-behaved $C^{*}$-pseudo-multiplicative unitary
  over $\mu$.
  \begin{enumerate}
  \item Assume that $\pmu$ is \'etale and that $\eta_{0} \in \fix(V)$ satisfies
    $\eta_{0}^{*}\eta_{0} = \Id_{H_{\mu}}$.
    \begin{enumerate}
    \item A bounded left counit $\hepsilon$ for
      $(\hA(V)^{\hbeta,\alpha}_{H},\hDelta_{V})$ is given by
      $\hepsilon(\ha):=\eta_{0}^{*} \ha \eta_{0}$. For all $\eta\in
      \beta,\xi \in\alpha$, we have $\hepsilon(\langle
      \eta|_{2}V|\xi\rangle_{2})=\eta^{*}\xi$. In particular,
      $\hepsilon$ does not depend on the choice of $\eta_{0}$,
      $\hepsilon(\hA(V)) = [\beta^{*}\alpha]$, and
      $[\beta^{*}\alpha]$ is a $C^{*}$-algebra.  If $V$ is regular,
      then $\hepsilon$  is a bounded counit.
    \item A bounded right Haar weight $\psi$ for
      $(A(V)^{\alpha,\beta}_{H},\Delta_{V})$ is given by
      $\psi(a):=\eta_{0}^{*}a\eta_{0}$.
    \end{enumerate}
  \item Assume that $\pmu$ is compact and that $\xi_{0} \in
    \cofix(V)$ satisfies $\xi_{0}^{*}\xi_{0} = \Id_{H_{\mu}}$.
    \begin{enumerate}
    \item A bounded right counit $\epsilon$ for
      $(A(V)^{\alpha,\beta}_{H},\Delta)$ is given by
      $\epsilon(a):=\xi_{0}^{*}a\xi_{0}$. For all $\eta\in
      \alpha,\xi \in\hbeta$, we have $\epsilon(\langle
      \eta|_{1}V|\xi\rangle_{1})=\eta^{*}\xi$. In particular,
      $\epsilon$ does not depend on the choice of $\eta_{0}$,
      $\epsilon(A(V)) = [\alpha^{*}\hbeta]$, and
      $[\alpha^{*}\hbeta]$ is a $C^{*}$-algebra.  If $V$ is regular,
      then $\hepsilon$  is a bounded counit.  \item A
      bounded left Haar weight $\widehat\phi$ for
      $(\hA(V)^{\hbeta,\alpha}_{H},\hDelta_{V})$ is given by
      $\widehat\phi(\ha):=\xi_{0}^{*}\ha\xi_{0}$.
    \end{enumerate}
  \end{enumerate}
\end{theorem}
\begin{proof}
  We only prove the assertions concerning
  $(\hA(V)^{\hbeta,\alpha}_{H},\hDelta_{V})$, the corresponding
  assertions for $(A(V)^{\alpha,\beta}_{H},\Delta_{V})$ follow by
  replacing $V$ by $V^{op}$.

  i) (a) Evidently, $\hepsilon$ is a completely positive
  contraction. Let $\eta,\eta' \in \beta$ and $\xi,\xi' \in
  \alpha$. Then
  \begin{align} \label{eq:counit} \langle
    \eta|_{2}V|\xi\rangle_{2}\eta_{0} = \langle \eta|_{2}
    |\eta_{0}\rangle_{1}\xi = \eta_{0} \eta^{*}\xi =
 \eta_{0}\eta_{0}^{*} \langle
    \eta|_{2}V|\xi\rangle_{2}\eta_{0} =
    \eta_{0}\hepsilon(\eta|_{2}V|\xi\rangle_{2}).
  \end{align}
  Now, $\hepsilon$ is a $*$-homomorphism and
  $\hepsilon(\hA(V))=[\beta^{*}\alpha]$ because
  \begin{align*}
    \eta_{0}^{*}\langle \eta|_{2}V|\xi\rangle_{2} \eta_{0}
    \eta_{0}^{*} \langle \eta'|_{2}V|\xi'\rangle_{2}\eta_{0} &=
    \eta_{0}^{*}\langle \eta|_{2}V|\xi\rangle_{2} \langle
    \eta'|_{2}V|\xi'\rangle_{2}\eta_{0}, & \hepsilon(\langle
    \eta|_{2}V|\xi\rangle_{2}) &= \eta_{0}^{*}\eta_{0} \eta^{*}\xi =
    \eta^{*}\xi.
  \end{align*}
  Since $[\eta_{0}^{*}\alpha]=B$ and $[\eta_{0}^{*}\hbeta]=B^{op}$,
  the map $\hepsilon$ is morphism of
  $C^{*}$-$(\mu,\mu^{op})$-algebras $\hA(V)^{\hbeta,\alpha}_{H}$ and
  $[\beta^{*}\alpha]^{B^{op},B}_{H_{\mu}}$. It is a left
  counit because $(\hepsilon \ast \Id)\big(\hDelta(\ha)\big) =
  \langle \eta_{0}|_{1}V^{*}(1 \stensor{\alpha}{\beta}
  \ha)V|\eta_{0}\rangle_{1} = \langle \eta_{0}|_{1}(1
  \stensor{\alpha}{\beta} \ha)|\eta_{0}\rangle_{1} = \ha$ for all
  $\ha \in \hA(V)$.

  Assume that $V$ is regular, and
  consider the following diagram:
  \begin{gather*}\smalldiagram 
    \xymatrix@C=17pt@R=25pt{ {H} \ar[d]^{|\eta_{0}\rangle_{2}}
      \ar[r]^{|\xi\rangle_{2}} & {\sHsource}
      \ar[d]^{|\eta_{0}\rangle_{2}} \ar[r]^{V}
      \ar@{}[dr]|{\mathrm{(*)}} & {\sHrange}
      \ar[rd]+<-2pt,7pt>_(0.4){|\eta_{0}\rangle_{2}} \ar[rr]^{\Id} &&
      {\sHrange} \ar@{<-}[ld]+<2pt,7pt>^(0.4){\langle\eta_{0}|_{2}}
      \ar[r]^{\langle \eta|_{2}} & H \\
      {\sHsource} \ar[r]^(0.4){|\xi\rangle_{3}} \ar `d/0pt[d]
      `r[drrrrr]^{\hDelta_{V}(\langle \eta|_{2}V|\xi\rangle_{2})}
      [rrrrr] & {\sHone} \ar[r]^(0.7){V_{13}V_{23}} &&
      {\hspace{-2cm}\sHfour\hspace{-2cm}} &
      \ar[r]^(0.3){\langle\eta|_{3}} & {\sHsource}
      \ar[u]^{\langle\eta_{0}|_{2}}  \\
      &&&&& }
  \end{gather*}
  The lower cell commutes  by the proof of 
  \cite[Lemma 4.13]{timmer:cpmu}, cell (*) commutes because
  $V_{23}|\eta_{0}\rangle_{2} = |\eta_{0}\rangle_{2}$, and the other
  cells commute as well. Since $\eta \in \beta$ and $\xi \in
  \alpha$  were arbitrary, $\hepsilon$ is a bounded right counit.

  ii) (b) By Remark \ref{remarks:fixed} i), $[\xi_{0}^{*} \hA(V)
\xi_{0}] = [\xi_{0}^{*} \rho_{\alpha}(B^{op})\hA(V)
\rho_{\alpha}(B^{op})\xi_{0}] \subseteq [\beta^{*}\hA(V)\beta]
\subseteq B^{op}$. Hence, the given formula defines a completely
positive contraction $\widehat{\phi} \colon \hA(V) \to B^{op}$. Since
$\rho_{\alpha}(b^{op})\xi_{0}=\xi_{0}b^{op}$ for all $b^{op}
\in B^{op}$, condition i) of Definition \ref{definition:cqg-haar-weights}
holds.  Condition ii)  holds
because  for all $\ha \in \hA(V)$ and $\eta,\eta' \in \hbeta$,
  \begin{align*}
    \xi_{0}^{*} \langle \eta|_{1} \hDelta_{V}(\ha)|\eta'\rangle_{1}\xi_{0} &=
    \eta^{*}\langle\xi_{0}|_{2}V^{*}(\Id \stensor{\alpha}{\beta}
    \ha)V|\xi_{0}\rangle_{2} \eta' \\ &=
    \eta^{*}\langle\xi_{0}|_{2}(\Id \stensor{\alpha}{\beta}
    \ha)|\xi_{0}\rangle_{2} \eta' =
    \eta^{*} \rho_{\alpha}\big(\xi_{0}^{*}\ha\xi_{0}\big)\eta'. \qedhere
  \end{align*}
\end{proof}

\paragraph{Balanced $C^{*}$-pseudo-multiplicative unitaries and
  $C^{*}$-pseudo-Kac systems}

Weak $C^{*}$-pseudo-Kac systems were introduced in
\cite{timmer:ckac} as a framework to construct reduced crossed
products for coactions of Hopf $C^{*}$-bimodules. Let us briefly
recall the definition.
\begin{definition}[{\cite{timmer:ckac}}] \label{definition:add-kac}
  A {\em balanced  $C^{*}$-pseudo-multiplicative unitary over $\mu$}
  is a tuple   $(H,\alpha,\halpha,\beta, \hbeta,U,V)$, where
  $(H,\alpha,\halpha,\beta,\hbeta)$ is a
  $C^{*}$-$(\mu,\mu,\mu^{op},\mu^{op})$-module, $V \colon \sHsource
  \to \sHrange$ is a unitary and $U  \colon H \to H$ is a symmetry
  satisfying the following conditions:
  \begin{enumerate}
  \item $U\alpha=\halpha$ and $U\beta=\hbeta$;
  \item  $(H,\hbeta,\alpha,\beta,V)$,
    $(H,\halpha,\hbeta,\alpha,\widecheck{V})$,
    $(H,\alpha,\beta,\halpha,\widehat{V})$ are well-behaved
    $C^{*}$-pseudo-multiplicative unitaries, where $\checkV$ and
    $\hatV$ are defined by
    \begin{align*}
      \widecheck{V} &:=\Sigma (1 \stensor{\alpha}{\beta} U)V(1
      \stensor{\hbeta}{\halpha} U)\Sigma \colon \checkHsource \to
      \checkHrange, \\ \widehat{V} &:= \Sigma (U
      \stensor{\alpha}{\beta} 1)V(U \stensor{\beta}{\alpha} 1)\Sigma \colon
      \hatHsource \to \hatHrange.
    \end{align*}
  \end{enumerate}

  A {\em weak $C^{*}$-pseudo-Kac system over $\mu$}  is a balanced
  $C^{*}$-pseudo-multiplicative unitary
  $(H,\alpha,\halpha,\beta,\hbeta,U,V)$ such that
  $(H,\beta,\alpha,\hbeta,V)$ is well-behaved and
  $[\hA(V),U\hA(V)U]=0=[A(V), UA(V)U]$.
  A weak $C^{*}$-pseudo-Kac system
  $(H,\alpha,\beta,\halpha,\hbeta,V,U)$ is a {\em $C^{*}$-pseudo-Kac
    system} if 
  $(H,\hbeta,\alpha,\beta,V)$,
  $(H,\halpha,\hbeta,\alpha,\widecheck{V})$,
  $(H,\alpha,\beta,\halpha,\widehat{V})$ are regular and 
  $\big(\Sigma(1 \stensor{\alpha}{\beta} U)V\big)^{3} = \Id\in {\cal
    L}(\sHsource)$.
\end{definition}
Let    $(H,\alpha,\halpha,\beta, \hbeta,U,V)$ be a balanced
$C^{*}$-pseudo-multiplicative unitary over $\mu$. Then by
\cite[Proposition 3.3]{timmer:ckac},
\begin{gather} \label{eq:add-legs}
  \begin{gathered}
    \hA(\checkV) = \Ad_{U}(A(V)), \quad \hDelta_{\checkV} = \Ad_{(U
      \rtensorh U)} \circ \Delta_{V} \circ \Ad_{U}, \quad
    A(\checkV) = \hA(V), \quad
    \Delta_{\checkV} = \hDelta_{V}, \\
    A(\hatV) = \Ad_{U}(\hA(V)), \quad \Delta_{\hatV} = \Ad_{(U
      \rtensorh U)} \circ \hDelta_{V} \circ \Ad_{U}, \quad
    \hA(\hatV) = A(V), \quad \hDelta_{\hatV} = \Delta_{V}.
  \end{gathered}
  \end{gather}
  In particular, $\checkV$ and $\hatV$ are well-behaved if $V$ is
  well-behaved.  
\begin{lemma} \label{lemma:add-fixed-kac}
  If $\ckac$ is a balanced $C^{*}$-pseudo-multiplicative unitary,
  then $\Fix(V)=U\Cofix(\hatV)=\Cofix(\checkV)$ and
  $\Cofix(V)=\Fix(\hatV)=U\Fix(\checkV)$. \qed
\end{lemma}
\begin{corollary} \label{corollary:fixed-cofixed}
  Let $\ckac$ be a balanced $C^{*}$-pseudo-multiplicative
  unitary over $\mu$, where $V$ is well-behaved.
  \begin{enumerate}
  \item Assume that $\pmu$ is \'etale and that $\eta_{0} \in
    \fix(V)$ satisfies $\eta_{0}^{*}\eta_{0} = \Id_{H_{\mu}}$.
    \begin{enumerate}
    \item A  bounded counit for
      $(\hA(V)^{\hbeta,\alpha}_{H},\hDelta_{V})$ is given by
      $\ha\mapsto\eta_{0}^{*} \ha \eta_{0}$.
    \item A bounded left Haar weight for
      $(A(V)^{\alpha,\beta}_{H},\Delta_{V})$ is given by
      $a \mapsto \eta_{0}^{*}U^{*}a U\eta_{0}$.
    \end{enumerate}
   \item  Assume that $\pmu$ is proper and $\xi_{0} \in \cofix(V)$ satisfies
    $\xi_{0}^{*}\xi_{0} = \Id_{H_{\mu}}$.
    \begin{enumerate}
    \item A bounded counit for
      $(A(V)^{\alpha,\beta}_{H},\Delta_{V})$ is given by
      $a \mapsto \xi_{0}^{*}a\xi_{0}$.
    \item A bounded right Haar weight for
      $(\hA(V)^{\hbeta,\alpha}_{H},\hDelta_{V})$ is given by
      $\ha \mapsto\xi_{0}^{*}U^{*}\ha U\xi_{0}$.
    \end{enumerate}
  \end{enumerate}
\end{corollary}
\begin{proof}
  Apply Theorem \ref{theorem:fixed-cofixed} to $\checkV$ or $\hatV$,
  respectively, and use Remark \ref{remarks:fixed} ii) and
\eqref{eq:add-legs}.
\end{proof}

\section{The dual Hopf $C^{*}$-bimodule}

\label{section:dual}

In the preceding section, we saw that the fundamental unitary
associated to a compact $C^{*}$-quantum groupoid gives rise to two
Hopf $C^{*}$-bimodules and that one of these two coincides with the
underlying Hopf $C^{*}$-bimodule of the initial $C^{*}$-quantum
groupoid. In this short section, we study the other Hopf
$C^{*}$-bimodule, which can be considered as (the underlying Hopf
$C^{*}$-bimodule of) the generalized Pontrjagin dual of the initial
$C^{*}$-quantum groupoid.  

In principle, the dual Hopf $C^{*}$-bimodules of compact
$C^{*}$-quantum groupoids should precisely exhaust the class of
\'etale $C^{*}$-quantum groupoids with compact base, but a precise
definition of \'etale $C^{*}$-quantum groupoids is not yet
available. However, we can describe some important ingredients like
the underlying Hopf $C^{*}$-bimodule, the unitary antipode, and the
counits of the dual of a compact $C^{*}$-quantum groupoid.

Throughout this section, let
$(B,\mu,A,r,\phi,s,\psi,\delta,R,\Delta)$ be a compact
$C^{*}$-quantum groupoid. We use the notation introduced in the
preceding sections.

\paragraph{The dual Hopf $C^{*}$-bimodule}
In Theorem \ref{theorem:pmu} and Proposition
\ref{proposition:pmu-right-leg}, we associated to the compact
$C^{*}$-quantum groupoid a regular $C^{*}$-pseudo-multiplicative
unitary $(H,\hbeta,\alpha,\beta,V)$. Now, we determine the
$C^{*}$-algebra of the associated Hopf $C^{*}$-bimodule
$(\hA(V)^{\hbeta,\alpha}_{H},\hDelta_{V})$.
\begin{proposition} \label{proposition:dual-algebra}
  \begin{enumerate}\item 
    For each $a \in A$, there exists an operator $\lambda(a) \in
    {\cal L}(H)$ such that $\lambda(a) \Lambda_{\nu}(d) =
    \Lambda_{\nu}\big(\langle\zeta_{\phi}|_{2}\Delta(d)|
    a^{op}\zeta_{\phi}\rangle_{2}\big)$ for all $d \in A$, and
    $\lambda(a)^{*} = J\lambda(R(a))J$.
  \item  $\langle
    x^{op}\zeta_{\phi}|_{2} V|y^{op}\zeta_{\psi}\rangle_{2} =
    \lambda(yx^{*})$ for all $x,y \in A$.
  \item  $\hA(V)=[\lambda(A)]$.
  \end{enumerate}
\end{proposition} 
  \begin{proof}
By definition,     $\hA(V)$ is the closed linear span of all
    operators of the form $\langle
    x^{op}\zeta_{\phi}|_{2}V|y^{op}\zeta_{\psi}\rangle_{2}$, where $x,y\in
    A$. But for all $x,y,d \in A$,
    \begin{align*}
      \langle x^{op}\zeta_{\phi}|_{2}
      V|y^{op}\zeta_{\psi}\rangle_{2} d\zeta_{\nu} &= \langle
      x^{op}\zeta_{\phi}|_{2} V(d\zeta_{\nu} \tl y^{op}\zeta_{\psi})
      \\ &= \langle x^{op}\zeta_{\phi}|_{2} \Delta(d)(\zeta_{\nu}
      \tl y^{op} \zeta_{\phi}) = \Lambda_{\nu}\big( \langle
      \zeta_{\phi}|_{2}\Delta(d)|
      (x^{op})^{*}y^{op})\zeta_{\phi}\rangle_{2}\big).
    \end{align*}
    This calculation proves the existence of the operators
    $\lambda(a)$ for all $a \in A$ and that
    $\hA(V)=[\lambda(A)]$. Finally, by Theorem
    \ref{theorem:pmu-inversion}, Lemma
    \ref{lemma:modules-antiunitary} and Proposition
    \ref{proposition:module-graph},
    \begin{align*}
      \lambda(yx^{*})^{*} &=
      \big(\langle x^{op}\zeta_{\phi}|_{2} V
      |y^{op}\zeta_{\psi}\rangle_{2}\big)^{*}  \\
      &=\langle y^{op}\zeta_{\psi}|_{2} (J
      \rtensor{\alpha}{J_{\mu}}{\beta} I)V(J
      \rtensor{\alpha}{J_{\mu}}{\beta} I)
      |x^{op}\zeta_{\phi}\rangle_{2} \\
      &=  J \langle I y^{op} \zeta_{\psi} J_{\mu}|_{2} V|
      Ix^{op}\zeta_{\phi} J_{\mu}\rangle_{2} J \\
      &= J \langle R(y^{*})^{op} \zeta_{\phi}|_{2} V | R(x^{*})^{op}
      \zeta_{\psi}\rangle_{2} J =
      J \lambda(R(x)^{*}R(y))J = J \lambda(R(yx^{*}))J. \qedhere
    \end{align*}
  \end{proof}

\paragraph{The associated weak $C^{*}$-pseudo-Kac system} 
 Put $U:=IJ=JI$.
\begin{theorem} \label{theorem:dual-kac}
  $(H,\alpha,\halpha,\beta,\hbeta,U,V) $ is a weak
  $C^{*}$-pseudo-Kac system.
\end{theorem}
The proof involves the following formula:
\begin{lemma}
  $I\lambda(a)I \Lambda_{\nu^{-1}}(d) =
  \Lambda_{\nu^{-1}}\big(\langle \zeta_{\psi}|_{1} \Delta(d)
  |R(a^{*})^{op} \zeta_{\psi}\rangle_{1}\big)$ for all $a,d \in A$.
  \end{lemma}
  \begin{proof}
    By Lemma \ref{lemma:cqg-strong-invariance}, we have for all $a,d
    \in A$
    \begin{align*}
      I\lambda(a)I\Lambda_{\nu^{-1}}(d)
      &= I\lambda(a) \Lambda_{\nu}(R(d)^{*}) \\
      &= I \langle \zeta_{\phi}|_{2} \Delta(R(d)^{*})
      |a^{op}\zeta_{\phi}\rangle_{2}
      I      \zeta_{\nu^{-1}} 
      = \big(\langle \zeta_{\psi}|_{1} \Delta(d) |R(a^{*})^{op} \zeta_{\psi}\rangle_{1}\big)
      \zeta_{\nu^{-1}}. \qedhere
    \end{align*}
 \end{proof}
 \begin{lemma} \label{lemma:dual-w-vhat}
   $\hatV=W$ and $\checkV = (J \rtensor{\alpha}{J_{\mu}}{\hbeta} J)V^{op}(J
   \rtensor{\alpha}{J_{\mu}}{\hbeta} J)$.
 \end{lemma}
 \begin{proof}
 Theorems \ref{theorem:pmu-inversion} and
   \ref{theorem:pmu-second} imply  $\checkV = (U
   \stensor{\beta}{\halpha} U)W(U \stensor{\halpha}{\hbeta} U) =
   (J \rtensor{\alpha}{J_{\mu}}{\hbeta} J)V^{op}(J
   \rtensor{\alpha}{J_{\mu}}{\hbeta} J)$ and   $\hatV = \Sigma (U
   \stensor{\alpha}{\beta} 1) (J \rtensor{\hbeta}{J_{\mu}}{\alpha}
   I) V^{*} (J \rtensor{\hbeta}{J_{\mu}}{\alpha} I) (U
   \stensor{\beta}{\alpha} 1)\Sigma = \Sigma(I
   \rtensor{\hbeta}{J_{\mu}}{\alpha} I) V^{*} (I
   \rtensor{\beta}{J_{\mu}}{\alpha} I) \Sigma = W$.
 \end{proof}
\begin{proof}[Proof of Theorem \ref{theorem:dual-kac}]
  By Lemma \ref{lemma:dual-w-vhat}, $(H,\alpha,\beta,\halpha,\hatV)$
  and  $(H,\halpha,\hbeta,\alpha,\checkV)$ are regular
  $C^{*}$-pseudo-multiplicative unitaries.  Clearly, we have $[A(V), UA(V)U]
  = [A(V),JIA(V)IJ] = [A(V),JA(V)J] = 0$.  It remains to show that
  $[\hA(V),U\hA(V)U] = 0$.  But for all $x,y,d \in A$,
  \begin{align*}
    I\lambda(x)I\lambda(y) d \zeta_{\nu} &= I\lambda(x) I \langle
    \zeta_{\phi}|_{2}\Delta(d)|y^{op}\zeta_{\phi}\rangle_{2}
    \delta^{-1/2}\zeta_{\nu^{-1}} \\
    &= \langle \zeta_{\psi}|_{1}\Delta\big(\langle
    \zeta_{\phi}|_{2}\Delta(d)|y^{op}\zeta_{\phi}\rangle_{2}
    \delta^{-1/2}\big) |R(x^{*})^{op}\zeta_{\psi}\rangle_{1} \zeta_{\nu^{-1}} \\
    &= \langle \zeta_{\psi}|_{1} \langle \zeta_{\phi}|_{3}
    \Delta^{(2)}(d)|y^{op}\zeta_{\phi}\rangle_{3}
    |\delta^{-1/2}R(x^{*})^{op}\zeta_{\psi}\rangle_{1} \delta^{-1/2}
    \zeta_{\nu^{-1}} \\
    &= \langle \zeta_{\phi}|_{2} \langle \zeta_{\psi}|_{1}
    \Delta^{(2)}(d) |\delta^{-1/2}R(x^{*})^{op}\zeta_{\psi}\rangle_{1}
    |y^{op}\zeta_{\phi}\rangle_{2} \zeta_{\nu} \\
    &= \lambda(y) \langle \zeta_{\psi}|_{1} \Delta(d)
    |\delta^{-1/2}R(x^{*})^{op}\zeta_{\psi}\rangle_{1} \zeta_{\nu} \\
    &= \lambda(y) \langle \zeta_{\psi}|_{1} \Delta(d\delta^{-1/2})
    |R(x^{*})^{op}\zeta_{\psi}\rangle_{1}\delta^{1/2} \zeta_{\nu} \\
    &= \lambda(y)I\lambda(x)I d \delta^{-1/2} \zeta_{\nu^{-1}}
    =\lambda(y)I\lambda(x)I d \zeta_{\nu}.
  \end{align*}
  Therefore, $[\hA(V),U\hA(V)U] = [\hA(V),IJ\hA(V)JI] =
  [\hA(V),I\hA(V)I] = 0$.
\end{proof}

\paragraph{Coinvolution and counit on the dual Hopf $C^{*}$-bimodule}
Proposition \ref{proposition:dual-algebra} immediately implies:
\begin{corollary} \label{corollary:dual-antipode} There exists a
  $*$-antiautomorphism $\hR \colon \hA(V) \to \hA(V)$, $\ha \mapsto
  J\ha^{*}J$. \qed
\end{corollary}
This $*$-antiautomorphism is a coinvolution of the Hopf
$C^{*}$-bimodule $\big(\hA(V)^{\hbeta,\alpha}_{H},\hDelta_{V}\big)$
in the sense that it reverses the comultiplication:
\begin{proposition} \label{proposition:dual-antipode-comult}
  $\hDelta \circ \hR = \Ad_{\Sigma} \circ (\hR
  \rfibre{\hbeta}{\alpha} \hR) \circ \hDelta$.
  \end{proposition}
  \begin{proof}
  By \eqref{eq:add-legs} and Lemma
  \ref{lemma:dual-w-vhat}, we have for all $\ha \in \hA(V)$
    \begin{align*}
      \hDelta_{V}(\ha) = \checkV(\ha \, \stensor{\halpha}{\hbeta}
      1)\checkV^{*} &= (J \rtensor{\alpha}{J_{\mu}}{\hbeta} J)
      \Sigma V^{*}\Sigma (J \rtensor{\halpha}{J_{\mu}}{\hbeta} J)
      (\ha \stensor{\halpha}{\hbeta} 1) (J
      \rtensor{\alpha}{J_{\mu}}{\hbeta} J)^{*} \Sigma V\Sigma (J
      \rtensor{\alpha}{J_{\mu}}{\hbeta} J)^{*} \\ &=  (\Ad_{\Sigma} \circ
      (\hR \rfibre{\hbeta}{\alpha}
      \hR) \circ \hDelta_{V}) (\hR(\ha)). \qedhere
    \end{align*}
  \end{proof}
 The constructions in Section \ref{section:add} yield a counit on $\hA(V)$:
\begin{proposition}
    \begin{enumerate}
    \item The Hopf $C^{*}$-bimodule
      $(\hA(V)^{\hbeta,\alpha}_{H},\hDelta_{V})$ has a bounded counit
      $\hepsilon$, given by $\hepsilon(\lambda(y^{*}x)) =
      \zeta_{\psi}^{*} \lambda(y^{*}x) \zeta_{\psi} = J_{\mu}
      \zeta_{\phi}^{*}x^{*}y\zeta_{\psi} J_{\mu}$ for all $x,y \in
      A$.
      \item $\hepsilon(\hR(\ha)) = J_{\mu}\hepsilon(\ha)^{*}J_{\mu}$ for
        all $\ha \in \hA(V)$.
    \end{enumerate}
  \end{proposition}
  \begin{proof}
    i) By Proposition \ref{proposition:pmu-cqg-fixed}, Theorem
    \ref{theorem:fixed-cofixed} i), and Corollary
    \ref{corollary:fixed-cofixed} i), the map $\hepsilon \colon
    \hA(V) \to {\cal L}(H_{\mu})$, $\ha \mapsto \zeta_{\psi}^{*}\ha
    \zeta_{\psi}$, is a bounded counit, and by Theorem
    \ref{theorem:fixed-cofixed} i) and Proposition
    \ref{proposition:dual-algebra}, $\hepsilon(\lambda(y^{*}x)) =
    \hepsilon\big(\langle J xJ\zeta_{\phi}|_{2}
    V|JyJ\zeta_{\psi}\rangle_{2}\big) = J_{\mu}
    \zeta_{\phi}^{*}x^{*}y\zeta_{\psi} J_{\mu}$ for all $x,y \in A$.

    ii) For all $\ha \in \hA(V)$, we have $\hepsilon(\hR(\ha)) =
    \zeta_{\psi}^{*} J \ha^{*} J \zeta_{\psi} =
    J_{\mu}(\zeta_{\psi}^{*}\ha \zeta_{\psi})^{*}J_{\mu} =
    J_{\mu}\hepsilon(\ha)^{*}J_{\mu}$.
  \end{proof}

\section{Principal compact $C^{*}$-quantum groupoids}

\label{section:principal}

In this section, we study compact $C^{*}$-quantum groupoids that are
principal.  Most importantly, we show that a principal compact
$C^{*}$-quantum groupoid is essentially determined by the
conditional expectation $\tau \colon B \to \tau(B) \subseteq Z(B)$
and the state $\mu|_{\tau(B)}$, and that the dual of a principal
compact $C^{*}$-quantum groupoid is the $C^{*}$-algebra of compact
operators on a certain $C^{*}$-module.

\paragraph{Principal compact $C^{*}$-quantum groupoids}
Recall that a compact groupoid $G$ is principal if the map $G \to
G^{0} \times G^{0}$ given by $x \mapsto (r(x),s(x))$ is injective
or, equivalently, if $C(G)=[r^{*}(C(G^{0}))s^{*}(C(G^{0}))]$. The
second condition suggests the following definition:
\begin{definition} \label{definition:principal} A compact
  $C^{*}$-quantum graph $(B,\mu,A,r,\phi,s,\psi,\delta)$ is {\em
    principal} if $A=[r(B)s(B^{op})]$, and a compact $C^{*}$-quantum
  groupoid $(B,\mu,A,r,\phi,s,\psi, \delta,R,\Delta)$ is {\em
    principal} if $A=[r(B)s(B^{op})]$.
\end{definition}
To simplify the following discussion, we only consider the case
where $\delta=1_{A}$.  Corollary \ref{corollary:cqg-mu} shows that
this is not a serious restriction.

Let $(B,\mu,A,r,\phi,s,\psi,1_{A})$ be a principal compact
$C^{*}$-quantum graph. Then there exist at most one coinvolution $R$
for $(B,\mu,A,r,\phi,s,\psi,1_{A})$ and at most one comultiplication
$\Delta$ for $A^{\alpha,\beta}_{H}$ because necessarily
$R(r(b)s(c^{op})) = s(b^{op}) r(c)$ and $ \Delta(r(b)s(c^{op})) =
r(b) \stensor{\alpha}{\beta} s(c^{op})$ for all $b,c \in B$.  We
shall give conditions for the existence of such a coinvolution and a
comultiplication, and determine when $(B,\mu,A,r,\phi,s,\psi,
1_{A},R,\Delta)$ is a principal compact $C^{*}$-quantum
groupoid. These conditions involve the completely positive
contractions $\tau =\psi \circ r \colon B \to Z(B^{op})\cong Z(B)$
and $\tau^{\dag} =\phi \circ s \colon B^{op} \to Z(B)\cong
Z(B^{op})$ introduced in \eqref{eq:graph-tau}.  
\begin{theorem} \label{theorem:principal-graph} Let
  $(B,\mu,A,r,\phi,s,\psi,1_{A})$ be a principal compact
  $C^{*}$-quantum graph.  Then the following two conditions
  are equivalent:
  \begin{enumerate}
  \item There exist $R$, $\Delta$ such that $(B,\mu,A,r,\phi,s,\psi,
    1_{A},R,\Delta)$ is a compact $C^{*}$-quantum groupoid.
  \item $\tau(b)=\tau^{\dag}(b^{op})$ for all
    $b \in B$,  $\tau\colon B \to \tau(B)$ is a conditional
    expectation, $\mu \circ \tau = \mu$, $r \circ \tau = s \circ
    \tau$, and
    $\tau(b\sigma_{-i/2}^{\mu}(d))=\tau(d\sigma^{\mu}_{-i/2}(b))$
    for all $b,d \in \Dom(\sigma^{\mu}_{-i/2})$.
  \end{enumerate}
\end{theorem}

Before we prove this result, let us give an application: every
compact $C^{*}$-quantum groupoid has an underlying principal compact 
$C^{*}$-quantum groupoid. The nontrivial part of this assertion is
that the comultiplication restricts to a morphism of
$C^{*}$-$(\mu,\mu^{op})$-algebras.
\begin{corollary}
  Let $(B,\mu,A,r,\phi,s,\psi,1_{A},R,\Delta)$ be a compact
  $C^{*}$-quantum groupoid. Put $\tilde A:=[r(B)s(B^{op})] \subseteq
  A$, $\tilde \phi  := \phi|_{\tilde A}$, $\tilde \psi
  :=\psi|_{\tilde A}$, $\tilde R:=R|_{\tilde A}$. Then there exists
  a unique $\tilde \Delta$ such that $(B,\mu,\tilde
  A,r,\tilde \phi,s,\tilde \psi, 1_{\tilde A},\tilde R,\tilde
  \Delta)$ is a compact $C^{*}$-quantum groupoid. \qed
\end{corollary}

The proof of Theorem \ref{theorem:principal-graph} is divided into
several  steps. First, note that for
all $b,c \in B$,
\begin{align} \label{eq:principal-tau}
  \begin{aligned}
    \phi(s(b^{op})r(c)) &=\tau^{\dag}(b^{op})c, &
    \psi(r(b)s(c^{op})) &=
    \tau(b)c^{op}, \\
    \nu(s(b^{op})r(c)) &= \mu(\tau^{\dag}(b^{op})c), &
    \nu^{-1}(r(b)s(c^{op})) &= \mu^{op}(\tau(b)c^{op}).
  \end{aligned}
\end{align}

\begin{lemma}
  Let $(B,\mu,A,r,\phi,s,\psi,1_{A})$ be a principal compact
  $C^{*}$-quantum graph.  There exists a coinvolution $R$ for
  $(B,\mu,A,r,\phi,s,\psi,1_{A})$ if and only if
  $\tau(b)=\tau^{\dag}(b^{op})$ for all $b \in B$.
  \end{lemma}
\begin{proof}
  The only if part is Lemma \ref{lemma:graph-coinvolution} ii). So,
  assume that $\tau(b)=\tau^{\dag}(b^{op})$ for all $b \in B$. Then
  there exists an antiunitary $I \colon H \to H$ such that $I
  r(b)s(c^{op})\zeta_{\nu^{-1}} = s(b^{op})^{*}r(c)^{*}\zeta_{\nu}$
  for all $b,c \in B$ because by \eqref{eq:principal-tau},
  \begin{align*}
    \big\| s(b^{op})^{*}r(c)^{*}\zeta_{\nu}
    \big\|^{2}
    &= \nu\big(s((b^{*}b)^{op})r(cc^{*})\big) 
=    \nu^{-1}\big(r(b^{*}b)s((cc^{*})^{op})\big) 
    = \big\| r(b)s(c^{op})\zeta_{\nu^{-1}}
    \big\|^{2}.
  \end{align*}
  A short calculation shows that $Ir(b)^{*}s(c^{op})^{*}I
  =s(b^{op})r(c)$ for all $b,c \in B$. Therefore, we can define a
  $*$-homomorphism $R \colon A \to A$ by $a \mapsto Ia^{*}I$, and
  $R(r(b)s(c^{op}))=s(c^{op})r(b)$ for all $b,c \in B$. Finally,
  \eqref{eq:principal-tau} implies that
  $(\phi \circ R)(a)= \psi(a)^{op}$ for all $a \in A$.
\end{proof}
\begin{lemma} \label{lemma:principal-delta} Let
  $(B,\mu,A,r,\phi,s,\psi,1_{A})$ be a principal compact
  $C^{*}$-quantum graph such  that $\tau(b)=\tau^{\dag}(b^{op})$
  for all $b \in B$,  $\tau\colon B \to \tau(B)$ is a
  conditional expectation, and  $\mu \circ \tau = \mu$.
  \begin{enumerate}
  \item For all $d,e \in B$, where $e$ is analytic for
    $\sigma^{\mu}$, there exists an operator $T_{d,e} \in {\cal
      L}(H,\sHrange)$ such that for all $x\in r(B) \cup r(B)^{op}$ and $y \in
    s(B^{op}) \cup s(B^{op})^{op}$,
    \begin{align} \label{eq:principal-delta-1} 
      T_{d,e}xy \zeta_{\nu}
      &=x\zeta_{\psi} \tr de^{op}\zeta_{\mu} \tl y \zeta_{\phi},
  \end{align}
  and for all $b,c,b',c',d',e' \in B$, where $e'$ is analytic for
  $\sigma^{\mu}$,
  \begin{align}\label{eq:principal-delta-3}
    T_{d,e}^{*} \big(r(b')\zeta_{\psi} \tr  d'e'{}^{op}\zeta_{\mu} \tl s(c'{}^{op})
    \zeta_{\phi} \big)= r(\tau(d^{*}d'
    \sigma^{\mu}_{-i/2}(e'e^{*}))r(b') s(c'{}^{op})\zeta_{\nu}.
  \end{align}
\item Put ${\cal T} :=\{ T_{d,e} \mid d,e \in B, \, e$ analytic for
  $\sigma^{\mu}\}$. Then $[{\cal T}\alpha] = \alpha \rt \alpha$,
  $[{\cal T}^{*}(\alpha \rt \alpha)]=\alpha$ and $[{\cal
    T}\beta]=\beta \lt \beta$, $[{\cal T}^{*}(\beta \lt
  \beta)]=\beta$.
\item There exists a comultiplication $\Delta$ for
  $A^{\alpha,\beta}_{H}$.
  \end{enumerate}
\end{lemma}
\begin{proof}
i)  Let $d,e$ be as in i). Then there exists a $T_{d,e} \in
  {\cal L}(H,\sHrange)$ such that equations
  \eqref{eq:principal-delta-1} and \eqref{eq:principal-delta-3} hold
for all $x \in r(B)$, $y \in s(B^{op})$  because
  \begin{align*}
    \big\langle r(b)\zeta_{\psi} \tr de^{op}\zeta_{\mu} \tl s(c^{op})
    \zeta_{\phi} &\big|r(b')\zeta_{\psi}d' \tr \zeta_{\mu} \tl
    s(c'{}^{op}) \zeta_{\phi}e'{}^{op}\big\rangle_{\sHrange} \\
    &=\big\langle \zeta_{\mu} \big|d^{*}    (e^{op})^{*} \psi(r(b^{*}b'))
\phi(s((c'c^{*})^{op}))d'e'{}^{op}\zeta_{\mu}\big\rangle \\
    &=\big\langle \zeta_{\mu} \big|\tau(b^{*}b')
    \tau(c'c^{*}) d^{*}d'   (e'e^{*})^{op }\zeta_{\mu}\big\rangle \\
    &= \mu\big(\tau(b^{*}b')
    \tau(c'c^{*})d^{*}d'\sigma^{\mu}_{-i/2}(e'e^{*})\big) \\
    &= \mu\big(\tau(b^{*}b') \tau(c'c^{*})
    \tau(d^{*}d'\sigma^{\mu}_{-i/2}(e'e^{*}))\big) \\
    &= \mu\big(\tau(b^{*}b'\tau(d^{*}d'\sigma^{\mu}_{-i/2}(e'e^{*})))
    \tau(c'c^{*})\big) \\
    &= \big\langle r(b)s(c^{op}) \zeta_{\nu} \big| r(\tau(d^{*}d'
    \sigma^{\mu}_{-i/2}(e'e^{*}))r(b')
    s(c'{}^{op})\zeta_{\nu}\big\rangle
  \end{align*}
  for all $b,c,b',c',d',e' \in B$, where $e'$ is analytic for
  $\sigma^{\mu}$. Using Lemma \ref{lemma:module} iii), one easily
  concludes that $T_{d,e}xy \zeta_{\nu} =x\zeta_{\psi} \tr de^{op}
  \zeta_{\mu} \tl y \zeta_{\phi}$ for all $x \in 
  r(B)^{op}$ and $y \in s(B^{op})^{op}$.

  ii) Let $b,c \in B$ and $d,e$ as in i). Then for all $f \in B$,
  \begin{align*}
    T_{d,e}r(b)^{op}s(c^{op})^{op}\zeta_{\phi} f^{op}\zeta_{\mu} &=
    T_{d,e} r(fb)^{op}s(c^{op})^{op}\zeta_{\nu} \\
    &=r(fb)^{op}\zeta_{\psi} \tr de^{op}\zeta_{\mu} \tl s(c^{op})^{op}
    \zeta_{\phi} \\
   &= |s(c^{op})^{op}r(e)^{op} \zeta_{\phi}\rangle_{2}
    r(b)^{op}s(d^{op})^{op}  \zeta_{\phi} f^{op} \zeta_{\mu}.
  \end{align*}
   This relation implies $[{\cal T}\beta]=\beta \lt
  \beta$, and the remaining assertions follow  similarly.

  iii) For all $b,c,d,e \in B$, where $e$ is analytic for
  $\sigma^{\mu}$, we have $T_{d,e}r(b)s(c^{op}) = r(b)\stensor{\alpha}{\beta}
  s(c^{op})$.  Now, the claim follows from  ii).
\end{proof}

\begin{proof}[Proof of Theorem \ref{theorem:principal-graph}]
 i) implies ii) by Proposition
  \ref{proposition:cqg-tau}. Conversely, assume that ii) holds.  Then the
  preceding lemmas imply that there exist a coinvolution $R$ for
  $(B,\mu,A,r,\phi$, $s,\psi,1_{A})$ and a comultiplication $\Delta$
  for $A^{\alpha,\beta}_{H}$. We  show that the conditions in
 Definition \ref{definition:cqg} hold.

 First, we check condition \ref{definition:cqg} i). Since
 $\rho_{\beta}=r$ and $\rho_{\alpha}=s$,
\begin{gather*}
  [\Delta(A)|\alpha\rangle_{1}] =
  [|\rho_{\beta}(B)\alpha\rangle_{1} \rho_{\alpha}(B^{op})] =
  [|\alpha B\rangle_{1} \rho_{\alpha}(B^{op})] =
  [|\alpha\rangle_{1}\rho_{\beta}(B)\rho_{\alpha}(B^{op})] =
  [|\alpha\rangle_{1}A], \\
  [\Delta(A)|\zeta_{\psi}\rangle_{1}A] = [
  |\rho_{\beta}(B)\zeta_{\psi}\rangle_{1} 
  \rho_{\beta}(B)A] =
  [|\rho_{\beta}(B)\zeta_{\phi}B\rangle_{1} A]
  =[|\rho_{\beta}(B)s(B^{op})^{op}\zeta_{\psi}\rangle_{1}A] =
  [|\alpha\rangle_{1}A].
\end{gather*}
Similar calculations show that
  $[\Delta(A)|\beta\rangle_{2}]=[\kbeta{2}A]$ and
  $[\Delta(A)|\zeta_{\phi}\rangle_{2}A] = [\kbeta{2}A]$.

  Next, $\phi$ is a bounded left  Haar weight for
  $(A^{\alpha,\beta}_{H},\Delta)$ because for all $b,c \in B$,
  \begin{align*}
    \langle
    \zeta_{\phi}|_{2}\Delta(r(b)s(c^{op}))|\zeta_{\phi}\rangle_{2} &=
    \langle \zeta_{\phi}|_{2}(r(b) \stensor{\alpha}{\beta}
    s(c^{op}))|\zeta_{\phi}\rangle_{2}  \\
    &= r(b)\rho_{\alpha}(\zeta_{\phi}^{*}s(c^{op})\zeta_{\phi}) \\
    &=
    r(b)s(\phi(s(c^{op})))  
    = r(b) s(\tau(c)) = r(b\tau(c)) = r(\phi(r(b)s(c^{op}))).
  \end{align*}
  A similar calculation shows that $\psi$ is a bounded right Haar weight
  for $(A^{\alpha,\beta}_{H},\Delta)$. 

  Finally, we prove that $\phi,\psi$ and $R$ satisfy the strong
  invariance condition \ref{definition:cqg} iii). By Lemma
  \ref{lemma:cqg-invariance-tau}, we have for all $b,c,d,e \in
  \Dom(\sigma^{\mu}_{-i/2})$
  \begin{align*}
    \langle \zeta_{\psi}|_{1}\Delta(r(b)s(c^{op}))((r(d)s(e^{op}))^{op}
    &\,\stensor{\alpha}{\beta} 1) |\zeta_{\psi}\rangle_{1} = r(e)
    s(c^{op}) r(\tau(b\sigma^{\mu}_{-i/2}(d))) \\
    &= R\big(s(e^{op})r(c) s(\tau(d\sigma^{\mu}_{-i/2}(b)\big) \\
    &= R\big( \langle
    \zeta_{\psi}|_{1}\Delta(r(d)s(e^{op}))((r(b)s(c^{op}))^{op}
    \stensor{\alpha}{\beta} 1) |\zeta_{\psi}\rangle_{1} \big).
  \end{align*}
  Since $\Dom(\sigma^{\mu}_{-i/2}) \subseteq B$ is dense, condition
  \ref{definition:cqg} iii)
    holds.
\end{proof}

\paragraph{The reconstruction of a principal compact $C^{*}$-quantum
groupoid}
A principal compact $C^{*}$-quantum groupoid is completely
determined by the conditional expectation $\tau \colon B \to \tau(B)
\subseteq Z(B)$ and
can be reconstructed from $\tau$ as follows.  Assume that
\begin{itemize}
\item $C$ is a commutative unital $C^{*}$-algebra with a faithful
  state $\upsilon$,
\item $B$ is a unital $C^{*}$-algebra with a $\upsilon$-module
  structure $(\iota,\tau)$ such that $\iota(C) \subseteq Z(B)$.
\end{itemize}
We put $\mu:=\upsilon \circ \tau$ and identify $C$ with $\iota(C)$
via $\iota$.
\begin{lemma}
  $\tau(b\sigma^{\mu}_{-i/2}(d))=\tau(d\sigma^{\mu}_{-i/2}(b))$ for
  all $b,d \in \Dom(\sigma^{\mu}_{-i/2})$. 
\end{lemma}
\begin{proof}
  For all $c \in C$, we have $\sigma^{\mu}_{t}(c) =
  \sigma^{\upsilon}_{t}(c)=c$ for all $t \in \reals$ by Lemma
  \ref{lemma:graph-rieffel} and hence
  \begin{align*}
    \upsilon(c^{*}\tau(b\sigma^{\mu}_{-i/2}(d)) =
    \mu(c^{*}b\sigma^{\mu}_{-i/2}(d))
    &=\langle\Lambda_{\mu}(b^{*}c)|J\Lambda_{\mu}(d^{*})\rangle \\
    &=
    \langle \Lambda_{\mu}(d^{*})| J\Lambda_{\mu}(b^{*}c)\rangle \\
    &=
    \mu(d\sigma^{\mu}_{-i/2}(c^{*}b))  = 
    \mu(dc^{*}\sigma^{\mu}_{-i/2}(b))
    =     \upsilon(c^{*}\tau(d\sigma^{\mu}_{-i/2}(b)).
  \end{align*}
  Since $c \in C$ was arbitrary and $\upsilon$ faithful, the claim follows.
\end{proof}
As in Proposition \ref{lemma:module}, we define an isometry
$\zeta_{\tau} \colon H_{\upsilon} \to H_{\mu}$ by $
\Lambda_{\upsilon}(c) \mapsto \Lambda_{\mu}(c)$, identify $B,B^{op}$
with $C^{*}$-subalgebras of ${\cal L}(H_{\mu})$ via the
GNS-representations, and put 
\begin{align*}
  \gamma &:=[B\zeta_{\tau}] \subseteq {\cal
    L}(H_{\upsilon},H_{\mu}), & \gamma^{op} &:=[B^{op}\zeta_{\tau}]
  \subseteq {\cal L}(H_{\upsilon},H_{\mu}).
\end{align*}
\begin{proposition} \label{proposition:principal-construction} There
  exists a unique principal compact $C^{*}$-quantum groupoid
  $(B,\mu,A$, $r,\phi,s,\psi, 1_{A},R,\Delta)$ such that $A = B
  \stensor{\gamma^{op}}{\gamma} B^{op} \subseteq {\cal L}\big((H_{\mu})
  \stensor{\gamma^{op}}{\gamma} (H_{\mu})\big)$ and for all $b,c \in
  B$,
    \begin{align*}
      r(b) &= b \gtensor 1^{op}, & \phi(b \gtensor c^{op}) &=
      b\tau(c), & s(c^{op}) &= 1 \gtensor c^{op}, & \psi(b
      \gtensor c^{op}) &= \tau(b)c^{op}.
    \end{align*}
\end{proposition}
\begin{proof}
  Routine arguments show that there exists a unique principal
  compact $C^{*}$-quantum graph $(B,\mu,A, r,\phi,s,\psi, 1_{A})$
  with $A$, $r$, $s$, $\phi$, $\psi$ as above; let us only note that
  the completely positive contractions $\phi\colon A \to B$ and
  $\psi \to B^{op}$ are well-defined because they are given by $x
  \mapsto \langle \zeta_{\tau}|_{2}x |\zeta_{\tau}\rangle_{2}$ and
  $x \mapsto \langle \zeta_{\tau}|_{1}x|\zeta_{\tau}\rangle_{1}$,
  respectively. Now, the assertion follows from Theorem
  \ref{theorem:principal-graph}. 
\end{proof}
Every principal compact $C^{*}$-quantum groupoid is of the form
constructed above:
\begin{proposition} \label{proposition:principal-description}
  Let $(B,\mu,A,r,\phi,s,\psi,1_{A},R,\Delta)$ be a principal
  compact $C^{*}$-quantum group\-oid and put $\tau=\psi \circ r$.
  \begin{enumerate}
  \item $C:=\tau(B)$ is a commutative unital $C^{*}$-algebra,
    $\upsilon:=\mu|_{C}$ is a faithful state on $C$, 
    $(\Id,\tau)$ is a
    $\upsilon$-module structure on $B$, and $\mu=\upsilon \circ \tau$.
\end{enumerate}
Denote by $\zeta_{\tau} \colon H_{\upsilon} \to H_{\tau}$ the
isometry $c \zeta_{\upsilon} \mapsto c\zeta_{\tau}$  and put
$\gamma:=[B\zeta_{\tau}]$, $\gamma^{op}:=[B^{op}\zeta_{\tau}]$.
 \begin{enumerate}\setcounter{enumi}{1}
  \item There exists a unitary $\Xi\colon H_{\nu} \to (H_{\mu})
    \gtensor (H_{\mu})$ such that for all $b,c \in  B$,
    \begin{align*}
      \Xi \big( r(b)^{op}s(c^{op})^{op}\zeta_{\nu}\big)&=
      b^{op}\zeta_{\tau} \tr \zeta_{\upsilon} \tl c\zeta_{\tau} &&
      \text{and} & \Xi\big(r(b)s(c^{op})\zeta_{\nu}\big) &=
      b\zeta_{\tau} \tr \zeta_{\upsilon} \tl c^{op}\zeta_{\tau}.
    \end{align*}
    Moreover, $\Xi  \hbeta=[|\gamma\rangle_{1} B^{op}]$ and
    $\Xi\alpha=[|\gamma^{op}\rangle_{1}B]$.
  \item $\Ad_{\Xi}$ restricts to an isomorphism $A \to B
    \stensor{\gamma^{op}}{\gamma} B^{op}$, $r(b)s(c^{op}) \mapsto b
    \gtensor c^{op}$.
  \end{enumerate}
\end{proposition}
\begin{proof}
  i) This follows directly from Proposition
  \ref{proposition:cqg-tau} and Proposition \ref{proposition:cqg-mu}.

  ii) There exists an isomorphism $\Xi \colon H_{\nu} \to (H_{\mu})
  \gtensor (H_{\mu})$ satisfying the first equation in ii) because
  by Proposition \ref{proposition:cqg-mu}, \eqref{eq:principal-tau},
  and i) $\| r(b)^{op}s(c^{op})^{op}\zeta_{\nu}\|^{2} =
  \nu(r(bb^{*})s((c^{*}c)^{op})) =
  \upsilon(\tau(bb^{*})\tau(c^{*}c)) = \| b^{op}\zeta_{\tau} \tr
  \zeta_{\upsilon} \tl c\zeta_{\tau}\|^{2}$ for all $b,c \in B$.
  From Lemma \eqref{lemma:modules} iii), one easily deduces
  $Tr(b)s(c^{op})\zeta_{\nu} = b\zeta_{\tau} \tr \zeta_{\upsilon}
  \tl c^{op}\zeta_{\tau}$ for all $b,c \in B$.
  Finally, $\Xi \hbeta = [|\gamma\rangle_{1}B^{op}]$ and
  $\Xi\alpha = [|\gamma^{op}\rangle_{1}B]$ because for all
  $b,c,d \in B$,
  \begin{gather*}
    \Xi r(b)s(c^{op})\zeta_{\psi} d^{op}\zeta_{\mu} =
\Xi    r(b)s(c^{op}d^{op}) \zeta_{\nu^{-1}} =
    b\zeta_{\tau} \tr \zeta_{\upsilon} \tl c^{op}d^{op}\zeta_{\tau}
    =
    |b\zeta_{\tau}\rangle_{1}c^{op} d^{op}\zeta_{\mu}, \\
    \Xi r(b)^{op}s(c^{op})^{op}\zeta_{\psi}d^{op}\zeta_{\mu}
    = r(b)^{op}s(c^{op})^{op}s(d^{op})\zeta_{\nu} =
    b^{op}\zeta_{\tau} \tr \zeta_{\upsilon} \tl cd^{op}\zeta_{\tau} =
    |b^{op}\zeta_{\tau}\rangle_{1}c d^{op}\zeta_{\mu}.
  \end{gather*}

  iii) Straightforward.
\end{proof}

\paragraph{The dual  Hopf $C^{*}$-bimodule}
Let $(B,\mu,A,r,\phi,s,\psi,\delta,R,\Delta)$ be a principal compact
$C^{*}$-quantum groupoid and $(H,\hbeta,\alpha$, $\beta,V)$ the
associated $C^{*}$-pseudo-multiplicative unitary (see Theorem
\ref{theorem:pmu}). We show that the dual Hopf $C^{*}$-bimodule
$(\hA(V)^{\hbeta,\alpha}_{H},\hDelta_{V})$ studied in Section
\ref{section:dual} can be identified with the $C^{*}$-algebra of
compact operators on a Hilbert $C^{*}$-module over $\tau(B)$. This
result is a (reduced) analogue of the result that for every
principal compact groupoid $G$, the irreducible representations of
$C^{*}(G)$ are labelled by the orbits $G^{0}/G$ and that each such
 representation is by all compact operators
\cite{renault}.

We use the notation of Proposition
\ref{proposition:principal-description} and denote by ${\cal
  K}_{\tau} \subseteq {\cal L}(H_{\mu})$ the $C^{*}$-algebra
corresponding to ${\cal K}_{C}(\gamma) \tr \Id \subseteq {\cal
  L}(\gamma \tr H_{\upsilon})$ with respect to the natural
isomorphism $\gamma \tr H_{\upsilon} \cong H_{\mu}$, $\xi \tr \zeta
\mapsto \xi\zeta$. Thus, ${\cal K}_{\tau} = [ \{ k_{b,c} \mid b,c
\in B\}]$, where $k_{b,c}\colon H_{\mu} \to H_{\mu}$ is given by $
d\zeta_{\mu} \mapsto b\tau(c^{*}d)\zeta_{\mu}$ for all $b,c \in B$.
Note that ${\cal K}_{\tau} \subseteq {\cal L}((H_{\mu})_{\gamma})$.
\begin{lemma}
    $({\cal K}_{\tau})_{H_{\mu}}^{B^{op},B}$ is a
    $C^{*}$-$(\mu^{op},\mu)$-algebra.
\end{lemma}
\begin{proof}
  Clearly, $(H_{\mu},B^{op},B)$ is a
  $C^{*}$-$(\mu^{op},\mu)$-module. We have $[\rho_{B^{op}}(B){\cal
    K}_{\tau}]={\cal K}_{\tau}=[\rho_{B}(B^{op}){\cal K}_{\tau}]$
  because for all $a,b,c,d \in B$, $a' \in
  \Dom(\sigma^{\mu}_{-i/2})$, 
  $\rho_{B^{op}}(a)k_{b,c} d\zeta_{\mu} =
    ab\tau(c^{*}d)\zeta_{\mu} = k_{ab,c}d\zeta_{\mu}$ and
    $\rho_{B}(a'{}^{op})k_{b,c} d\zeta_{\mu} = a'{}^{op}
    b\tau(c^{*}d)\zeta_{\mu} =
    b\sigma^{\mu}_{-i/2}(a')\tau(c^{*}d)\zeta_{\nu} =
    k_{b\sigma^{\mu}_{-i/2}(a'),c} d\zeta_{\mu}$.
\end{proof}
The comultiplication $\hDelta_{V}$ can be described in terms of the
isomorphism
\begin{align*}
  \Upsilon \colon (H_{\mu}) \stensor{B^{op}}{B} (H_{\mu}) = B^{op}
  \tr H_{\mu} \tl B \xrightarrow{\cong} H_{\mu}, \quad b^{op} \tr
  \zeta \tl c \mapsto b^{op}c\zeta.
\end{align*}
Note that $\Upsilon^{*}{\cal K}_{\tau}\Upsilon \subseteq ({\cal
  K}_{\tau}) \rfibre{B^{op}}{B} ({\cal K}_{\tau})$ because
$[\Upsilon^{*} {\cal K}_{\tau}\Upsilon|B^{op}\rangle_{1}] =
[\Upsilon^{*} {\cal K}_{\tau}B^{op}] = [\Upsilon^{*} B^{op} {\cal
  K}_{\tau}] = [|B^{op}\rangle_{1}{\cal K}_{\tau}]$ and similarly
$[\Upsilon^{*} {\cal K}_{\tau}\Upsilon|B\rangle_{2}] =
[|B\rangle_{1}{\cal K}_{\tau}]$.
\begin{theorem} \label{theorem:principal-dual} Let
  $(B,\mu,A,r,\phi,s,\psi,1_{A},R,\Delta)$ be a principal compact
  $C^{*}$-quantum groupoid and
  $(\hA(V)^{\hbeta,\alpha}_{H},\hDelta_{V})$ the dual Hopf
  $C^{*}$-bimodule.  
  \begin{enumerate}
  \item There exists an isomorphism of
    $C^{*}$-$(\mu^{op},\mu)$-algebras $j \colon ({\cal
      K}_{\tau})^{B^{op},B}_{H_{\mu}} \to
    \hA(V)^{\hbeta,\alpha}_{H}$, given by $k \mapsto \Xi^{*}(\Id
    \gtensor k)\Xi$.
  \item $\hDelta_{V}\circ j = (j \ast  j) \circ \Ad_{\Upsilon}^{-1}$.
  \item $\hR(j(k_{b,c})) = j(k_{c',b'})$, where $c' =
    \sigma^{\mu}_{i/2}(c)^{*}$ and $b' = \sigma^{\mu}_{i/2}(b)^{*}$
    for all $b,c \in \Dom(\sigma^{\mu}_{i/2})$.
  \item $\hepsilon \circ j = \Id_{{\cal K}_{\tau}}$.
  \end{enumerate}
  \end{theorem}
\begin{proof}
  i) Let $b,c \in B$ be analytic for $\sigma^{\mu}$ and put
  $a:=r(b)s(c^{op})$. Then the operator $\lambda(a)$
  defined in Proposition \ref{proposition:dual-algebra} acts as
  follows. For all $d,e \in B$,
  \begin{align*}
    \lambda(a) r(d)s(e^{op})\zeta_{\nu} &=
    \Lambda_{\nu}\big( \langle \zeta_{\phi}|_{2}(r(d)
    \stensor{\alpha}{\beta} s(e^{op})r(b)^{op}s(c^{op})^{op})
    |\zeta_{\phi}\rangle_{2}\big)  \\
    &=\rho_{\alpha}
    \big(b^{op}\zeta_{\phi}^{*}s(e^{op})s(c^{op})^{op}\zeta_{\phi}
    \big) r(d)\zeta_{\nu} 
    = s\big(b^{op}\phi(s((\sigma^{\mu}_{i/2}(c)e)^{op}))\big)r(d)\zeta_{\nu},
  \end{align*}
  and hence $\Xi \lambda(a) \Xi^{*}(d\zeta_{\tau} \tr
  \zeta_{\upsilon} \tl e^{op}\zeta_{\tau}) = d\zeta_{\tau} \tr
  \zeta_{\upsilon} \tl
  b^{op}\tau(\sigma^{\mu}_{i/2}(c)e)\zeta_{\tau}$.  Assume that $e
  \in \Dom(\sigma^{\mu}_{-i/2})$.  Then by Proposition
  \ref{proposition:cqg-tau}, Lemma \ref{lemma:module} iii), and
  $\sigma^{\mu}$-invariance of $\tau$,
\begin{align*}
    \Xi \lambda(a)\Xi^{*}(d\zeta_{\tau} \tr
    \zeta_{\upsilon} \tl \sigma^{\mu}_{-i/2}(e)\zeta_{\tau}) &= 
     d\zeta_{\tau} \tr \zeta_{\upsilon} \tl
    \sigma^{\mu}_{-i/2}(b)\tau(\sigma^{\mu}_{i/2}(c)e)\zeta_{\tau}
    \\
    &= d\zeta_{\tau} \tr \zeta_{\upsilon} \tl \sigma^{\mu}_{-i/2}(b)
    \tau(c\sigma^{\mu}_{-i/2}(e))\zeta_{\tau}.
  \end{align*}
  Therefore, $\Xi \lambda(a)\Xi^{*} = (\Id
  \gtensor k_{b',c^{*}})$, where $b' = \sigma^{\mu}_{-i/2}(b)$, and 
  $\hA(V) = \Xi^{*}(\Id \gtensor {\cal K}_{\tau})\Xi$. Since
  $\upsilon$ is faithful and $\Xi$  unitary, the map $j\colon {\cal
    K}_{\tau} \to \hA(V)$ given by $k \mapsto \Xi^{*}(\Id \gtensor
  k)\Xi$ is an isomorphism of $C^{*}$-algebras.

  It remains to show that $j$ is a morphism of
  $C^{*}$-$(\mu^{op},\mu)$-algebras.  Evidently, $tk=j(k)t$ for all
  $k \in {\cal K}_{\tau}$ and all $t \in
  [\Xi^{*}|\gamma^{(op)}\rangle_{1}]$.  By Proposition
  \ref{proposition:principal-description} ii),
  $[\Xi^{*}|\gamma\rangle_{1}B^{op}]= \hbeta$,
  $[\Xi^{*}|\gamma^{op}\rangle_{1}B] = \alpha$,
  $[\langle\gamma|_{1}\Xi \hbeta] = [\langle
  \gamma|_{1}|\gamma\rangle_{1}B^{op}] = [CB^{op}] = B^{op}$, and $
  [\langle \gamma^{op}|_{1}\Xi \alpha] =
  [\langle\gamma^{op}|_{1}|\gamma^{op}\rangle_{1}B] = [CB] = B$. The
  claim follows.

    ii) By definition of $\hDelta_{V}$ and $j$, we have for all $x \in
    \alpha$, $y \in \gamma$, $k \in {\cal K}_{\tau}$
    \begin{align*}
      \hDelta_{V}(j(k)) V^{*}|x\rangle_{1}\Xi^{*}|y\rangle_{1} &=
      V^{*}(1 \stensor{\alpha}{\beta}
      j(k))|x\rangle_{1}\Xi^{*}|y\rangle_{1} =
      V^{*}|x\rangle_{1} j(k)\Xi^{*}|y\rangle_{1} =
      V^{*}|x\rangle_{1} \Xi^{*}|y\rangle_{1} k.
    \end{align*}
    Likewise, by definition of $j\ast j$ and $\Upsilon$, we have for
    all $u\in \gamma^{op}$, $v \in \gamma$, $k \in {\cal K}_{\tau}$
    \begin{align*}
      ((j \ast j)(\Upsilon^{*} k \Upsilon)) (\Xi^{*} |u\rangle_{1}
      \stensor{B^{op}}{B} \Xi^{*}|v\rangle_{1})\Upsilon^{*} &=
      (\Xi^{*}|u\rangle_{1} \stensor{B^{op}}{B}
      \Xi^{*}|v\rangle_{1}) \Upsilon^{*}k\Upsilon \Upsilon^{*} \\ &=
      (\Xi^{*}|u\rangle_{1} \stensor{B^{op}}{B}
      \Xi^{*}|v\rangle_{1}) \Upsilon^{*}k.
    \end{align*}
    Now,  $[V^{*}|\alpha\rangle_{1} \Xi^{*}|\gamma\rangle_{1}] =
      [(\Xi^{*}|\gamma^{op}\rangle_{1} \stensor{B^{op}}{B}
      \Xi^{*}|\gamma\rangle_{1}) \Upsilon^{*}]$ because for
    all $b,c,d,e \in B$,
    \begin{align*}
      V^{*}|r(b)^{op}s(c^{op})^{op}\zeta_{\psi}\rangle_{1}
      \Xi^{*}|d\zeta_{\tau}\rangle_{1} e^{op} \zeta_{\mu} &=
      V^{*}\big(r(b)^{op}s(c^{op})^{op}\zeta_{\psi} \tr
      r(d)s(e^{op})\zeta_{\nu}\big)  \\
      &=V^{*}\big( r(b)^{op}\zeta_{\psi} \tr r(cd)s(e^{op}) \zeta_{\nu}\big) \\
      &= r(b)^{op}\zeta_{\psi} \tr r(cd)s(e^{op}) \zeta_{\nu} \\
      &= (\Xi^{*}|b^{op}\zeta_{\tau}\rangle_{1} \stensor{B^{op}}{B}
      \Xi^{*}|cd\zeta_{\tau}\rangle_{1}) (1 \tr \zeta_{\mu} \tl
      e^{op}) \\
      &= (\Xi^{*}|b^{op}\zeta_{\tau}\rangle_{1} \stensor{B^{op}}{B}
      \Xi^{*}|cd\zeta_{\tau}\rangle_{1}) \Upsilon^{*}
      e^{op}\zeta_{\mu}. 
    \end{align*}
    Since $[V^{*}|\alpha\rangle_{1} \Xi^{*}|\gamma\rangle_{1} H_{\mu}]
    = H$, we can conclude $ \hDelta_{V}(j(k)) = (j \ast j)(\Upsilon^{*}
    k \Upsilon) $ for all $k \in {\cal K}_{\tau}$.

    iii) Let $e \in \Dom(\sigma^{\mu}_{i/2})$ and $b,c \in
    \Dom(\sigma^{\mu}_{i/2})$. Since $Js(f^{op})^{op}\zeta_{\nu} =
    \sigma^{\nu^{op}}_{i/2}(s(f^{op})^{op})^{*} \zeta_{\nu} =
    s(\sigma^{\mu}_{-i/2}(f^{*})^{op})^{op}\zeta_{\nu}$ for all $f
    \in \Dom(\sigma^{\mu}_{i/2})$ and
    $\tau(b^{*}\sigma^{\mu}_{-i/2}(e^{*})) =
    \tau(e^{*}\sigma^{\mu}_{-i/2}(b^{*}))$, 
    \begin{align*}
      \hR(j(k_{b,c})) s(e^{op})^{op}\zeta_{\nu} &= Jj(k_{b,c})^{*}J
      s(e^{op})^{op}\zeta_{\nu} \\
      &= J j(k_{c,b}) s(\sigma^{\mu}_{-i/2}(e^{*})^{op})^{op}  \zeta_{\nu}\\
      &= J
      s((c\tau(b^{*}\sigma^{\mu}_{-i/2}(e^{*})))^{op})^{op}\zeta_{\nu}
\\
&= s(\sigma^{\mu}_{-i/2}
      (c^{*}\tau(\sigma^{\mu}_{-i/2}(b^{*})^{*}e))^{op})^{op}\zeta_{\nu}
      = s((c'\tau(b'^{*}e))^{op})^{op}\zeta_{\nu},
    \end{align*}
    where $b'=\sigma^{\mu}_{i/2}(b)^{*}$ and
    $c'=\sigma^{\mu}_{i/2}(c)^{*}$. The claim follows.

    iv) For all $b,c,d \in B$, 
    \begin{align*}
      \hepsilon(j(k_{b,c})) d\zeta_{\mu} = \zeta_{\psi}^{*}
      j(k_{b,c}) \zeta_{\psi} d\zeta_{\mu} &= \zeta_{\psi}^{*}
      j(k_{b,c}) s(d^{op})^{op}\zeta_{\nu}  \\ &=\zeta_{\psi}^{*}
      s((b\tau(c^{*}d))^{op})^{op} \zeta_{\nu} = b\tau(c^{*}d)
      \zeta_{\mu} = k_{b,c}d\zeta_{\mu}. \qedhere
    \end{align*}
\end{proof}

\paragraph{The  $C^{*}$-pseudo-Kac system}
Recall that in Theorem \ref{theorem:dual-kac}, we associated to
every compact $C^{*}$-quantum group\-oid a weak $C^{*}$-pseudo-Kac
system.
\begin{theorem} \label{theorem:principal-kac}
  Let $(B,\mu,A,r,\phi,s,\psi,1_{A},R,\Delta)$ be a principal
  compact $C^{*}$-quantum group\-oid. Then the weak
  $C^{*}$-pseudo-Kac system $(H,\alpha,\halpha,\beta,\hbeta,U,V)$ is
  a $C^{*}$-pseudo-Kac system.
\end{theorem}
\begin{proof}
  The $C^{*}$-pseudo-multiplicative unitaries
  $(H,\hbeta,\alpha,\beta,V)$,
  $(H,\halpha,\hbeta,\alpha,\widecheck{V})$,
  $(H,\alpha,\beta,\halpha,\widehat{V})$ are regular by Theorems
  \ref{theorem:pmu}, \ref{theorem:pmu-second} and Lemma
  \ref{lemma:dual-w-vhat}, and the operator $X:=\Sigma(1
  \stensor{\alpha}{\beta} U)V \in {\cal L}(\sHsource)$ satisfies $X^{3}
  = \Id$ because for all
  $b,c,d,e \in B$,
  \begin{align*}
    X^{3} \big(r(b)s(c^{op})\zeta_{\psi} &\tr
    r(d)s(e^{op})\zeta_{\nu}\big) = X^{2}\Sigma(1 \tr U)
    \big(r(b)\zeta_{\psi} \tr
    r(d)s(c^{op})s(e^{op})\zeta_{\nu}\big) \\
    &=X^{2} \big(s(d^{op})r(ec)\zeta_{\nu} \tl r(b)
    \zeta_{\psi }\big) \\
    &= X\Sigma(1 \tr U) \big(
    r(ec)\zeta_{\nu} \tl s(d^{op})r(b)\zeta_{\phi}\big) \\
    &= X\big( r(d)s(b^{op}) \zeta_{\psi} \tr
    r(ec)\zeta_{\nu}
    \big) \\
    &=\Sigma(1 \tr U) \big(r(d) \zeta_{\psi} \tr
    s(b^{op})r(ec)\zeta_{\nu}\big)
    \\
    &= r(b)s((ec)^{op})\zeta_{\nu} \tl r(d)
    \zeta_{\psi} \\ &= r(b)s(c^{op})\zeta_{\psi}e^{op} \tr  r(d)
    \zeta_{\nu} = r(b)s(c^{op})\zeta_{\psi} \tr
    r(d)s(e^{op})\zeta_{\nu}. \qedhere
    \end{align*}
\end{proof}

\section{Compact and \'etale groupoids}

Prototypical examples of compact $C^{*}$-quantum groupoids are the
function algebra of a compact groupoid and the reduced groupoid
$C^{*}$-algebra of an \'etale groupoid with compact space of units.
In this section, we construct these examples, determine the
associated dual Hopf $C^{*}$-bimodules, and show that the associated
weak $C^{*}$-pseudo-Kac systems are $C^{*}$-pseudo Kac systems.  We
shall use some results from \cite{timmer:cpmu} and \cite{timmer:ckac}
which we recall first.

\paragraph{Preliminaries on locally compact groupoids}
Let us fix some notation and terminology related to locally compact
groupoids; for details, see \cite{renault} or \cite{paterson}.  

Throughout this section, let $G$ be a locally compact, Hausdorff,
second countable groupoid. We denote its unit space by $G^{0}$, its
range map by $r_{G}$, its source map by $s_{G}$, and put
$G^{u}:=r_{G}^{-1}(\{u\})$, $G_{u}:=s_{G}^{-1}(u)$ for each $u \in
G^{0}$.

Let $\lambda$ be a left Haar system on $G$ and denote by
$\lambda^{-1}$ the associated right Haar system. Let $\mu_{G}$ be a
probability measure on $G^{0}$ with full support and define measures
$\nu_{G},\nu_{G}^{-1}$ on $G$ by
 \begin{align*}
   \int_{G} f \,d \nu_{G} &:= \int_{G^{0}} \int_{G^{u}} f(x) \, d\lambda^{u}(x)
\,   d\mu_{G}(u),
&   \int_{G} f d\nu_{G}^{-1} &= \int_{G^{0}}  \int_{G_{u}} f(x)
   d\lambda^{-1}_{u}(x) \, d\mu_{G}(u)
 \end{align*}
 for all $f \in C_{c}(G)$. Thus, $\nu_{G}^{-1}=i_{*}\nu_{G}$, where $i\colon G
 \to G$ is given by $x \mapsto x^{-1}$.  We assume that 
 $\mu$ is quasi-invariant in the sense that $\nu_{G}$ and $\nu_{G}^{-1}$ are
 equivalent, and denote by
 $D:=d\nu_{G}/d\nu_{G}^{-1}$ the Radon-Nikodym derivative. 

{\em
 In the following applications, we shall always assume that the unit space
 $G^{0}$ is compact and that the Radon-Nikodym derivative $D$ is
 continuous.}

\paragraph{The $C^{*}$-pseudo-Kac system of a locally
  compact groupoid}
In \cite{timmer:cpmu} and \cite{timmer:ckac}, we associated to $G$ a
$C^{*}$-pseudo-multiplicative unitary and a $C^{*}$-pseudo-Kac
system as follows.

Denote by $\mu$ the trace on $C(G^{0})$ given by $f \mapsto
\int_{G^{0}} f d\mu_{G}$.  Put $K:=L^{2}(G,\nu_{G})$ and define
representations $r,s\colon C(G^{0}) \to {\cal L}(K)$ such that
 for all $x \in G$, $\xi \in C_{c}(G)$, and $f \in C(G^{0})$,
 \begin{align*}
  \big(r(f)\xi\big)(x) &:=
 f\big(r_{G}(x)\big)\xi(x), & 
\big(s(f)\xi\big)(x) &:=
 f\big(s_{G}(x)\big) \xi(x).
 \end{align*}
 We define Hilbert $C^{*}$-$C(G^{0})$-modules $L^{2}(G,\lambda)$ and
 $L^{2}(G,\lambda^{-1})$ as the respective completions of the
 pre-$C^{*}$-module $C_{c}(G)$, where for all $\xi,\xi' \in
 C_{c}(G)$, $u \in G^{0}$, $f \in C(G^{0})$, $x \in G$,
 \begin{align*}
   \langle \xi'|\xi\rangle(u)&= \int_{G^{u}}
   \overline{\xi'(x)}\xi(x) d\lambda^{u}(x), & (\xi f)(x) &=
   \xi(x)f(r_{G}(x))  &&\text{in case of } L^{2}(G,\lambda), \\
   \langle \xi'|\xi\rangle(u)&= \int_{G_{u}}
   \overline{\xi'(x)}\xi(x) d\lambda^{-1}_{u}(x), & (\xi f)(x) &=
   \xi(x)f(s_{G}(x)) &&\text{in case of } L^{2}(G,\lambda^{-1}).
\end{align*}

There exist isometric embeddings $j\colon L^{2}(G,\lambda) \to {\cal
  L}(H_{\mu},K)$ and $\hat j \colon L^{2}(G,\lambda^{-1}) \to {\cal
  L}\big(H_{\mu},K\big)$ such that for all $\xi \in C_{c}(G)$,
$\zeta \in H_{\mu}$, $x \in G$,
\begin{align} \label{eq:groupoid-j}
  \big(j(\xi) \zeta\big)(x) &= \xi(x)\zeta(r_{G}(x)), &
  \big(\hat j(\xi)\zeta\big)(x) &= \xi(x) D(x)^{-1/2} \zeta(s_{G}(x)).
\end{align}
Put $\rho:=j(L^{2}(G,\lambda))$ and $ \sigma:= \hat
j(L^{2}(G,\lambda^{-1}))$. Then $(K,\sigma,\rho,\rho)$ is a
$C^{*}$-$(\mu^{op},\mu,\mu^{op})$-module, and $j$ and
$\hat j$ are unitary maps of Hilbert $C^{*}$-modules over $C(G^{0})$.

 Define measures $\nu^{2}_{s,r}$ on $G^{2}_{s,r}:=\{(x,y)\in
G \times G \mid s_{G}(x)=r_{G}(y)\}$ and $\nu^{2}_{r,r}$ on $G^{2}_{r,r}:=\{
(x,y) \in G^{2}\mid r_{G}(x)=r_{G}(y)\}$ by
 \begin{align*}
   \int_{G^{2}_{s,r}} \!\! f\, d\nu^{2}_{s,r} &:= \int_{G^{0}}
   \int_{G^{u}} \int_{G^{s_{G}(x)}} f(x,y) \, d\lambda^{s_{G}(x)}(y)
   \, d\lambda^{u}(x) \, d\mu_{G}(u), \\
    \int_{G^{2}_{r,r}} \!\! g\, d\nu^{2}_{r,r} &:= \int_{G^{0}}
   \int_{G^{u}}\int_{G^{u}} g(x,y)\, d\lambda^{u}(y)\, d\lambda^{u}(x)\,
   d\mu_{G}(u)
 \end{align*}
 for all $f \in C_{c}(G^{2}_{s,r})$, $g\in C_{c}(G^{2}_{r,r})$. Then
 there exist isomorphisms
  \begin{align*}
    \Phi_{\sigma,\rho} \colon \gHsource &\to
    L^{2}(G^{2}_{s,r},\nu^{2}_{s,r}), &
    \Phi_{\rho,\rho} \colon \gHrange &\to L^{2} (G^{2}_{r,r},\nu^{2}_{r,r})
  \end{align*}
such that for all $\eta,\xi \in C_{c}(G)$, $\zeta \in
C_{c}(G^{0})$, $(x,y) \in G^{2}_{s,r}$, $(x',y') \in G^{2}_{r,r}$,
\begin{align*}
  \Phi_{\sigma,\rho}\big(\hat j(\eta) \tr \zeta \tl
  j(\xi)\big)(x,y) &= \eta(x)D(x)^{-1/2} \zeta(s_{G}(x)) \xi(y),  \\
  \Phi_{\rho,\rho}\big(j(\eta) \tr \zeta \tl j(\xi)\big)(x',y') &=
  \eta(x') \zeta(r_{G}(x')) \xi(y').    
\end{align*}
From now on, we identify $\gHsource$ with
$L^{2}(G^{2}_{s,r},\nu^{2}_{s,r})$ via $\Phi_{\sigma,\rho}$ and
$\gHrange$ with $L^{2}(G^{2}_{r,r},\nu^{2}_{r,r})$ via
$\Phi_{\rho,\rho}$ without further notice.

\begin{theorem}[{\cite{timmer:cpmu},\cite{timmer:ckac}}] \label{theorem:groupoid-vg}
  There exists a $C^{*}$-pseudo-Kac system
  $(K,\rho,\sigma,\rho,\sigma,U_{G},V_{G})$ such that for all $\omega
  \in C_{c}(G^{2}_{s,r})$, $(x,y) \in G^{2}_{r,r}$, $\xi \in
  C_{c}(G)$, $z \in G$,
  \begin{align*}
    (V_{G}\omega)(x,y) &= \omega(x,x^{-1}y), & (U_{G}\xi)(z)=\xi(z^{-1})D(z)^{-1/2}.
  \end{align*}
\end{theorem}
\begin{proposition} \label{proposition:groupoid-fixed}
  \begin{enumerate}
  \item If $G$ is $r$-discrete and $1_{G^{0}} \in C(G)$ denotes the
    characteristic function of the unit space, then
    $j(1_{G^{0}})=\hat j(1_{G^{0}}) \in \Fix(V_{G})$ and
    $(K,\sigma,\rho,\rho,V_{G})$ is compact.
  \item If $G$ is compact, then $j(1_{G}) \in \Cofix(V_{G})$ and
    $(K,\sigma,\rho,\rho,V_{G})$ is \'etale.
  \end{enumerate}
\end{proposition}
\begin{proof}
  Straightforward.
\end{proof}
The concrete Hopf $C^{*}$-bimodules
$\big(\hA(V_{G})^{\sigma,\rho}_{K},\hDelta_{V_{G}}\big)$ and
$\big(A(V_{G})_{K}^{\rho,\rho},\Delta_{V_{G}}\big)$ can be described
as follows.  Denote by $m \colon C_{0}(G) \to {\cal
  L}(L^{2}(G,\nu_{G}))$ the representation given by multiplication
operators.  Recall that for each $g \in C_{c}(G)$, there exists an operator $L(g) \in {\cal L}(K)$ such that
\begin{align*}
  \big(L(g)f\big)(y) = \int_{G^{r_{G}(y)}} g(x)D(x)^{-1/2}
  f(x^{-1}y) d\lambda^{r_{G}(y)}(x) \quad \text{for all } f \in
  C_{C}(G), y \in G,
\end{align*}
and that $L(g)^{*}=L(g^{*})$, where $g^{*}(x) = \overline{g(x^{-1})}$ for
all $x \in G$.  The reduced groupoid $C^{*}$-algebra $C^{*}_{r}(G)$
is the closed linear span of all operators $L(g)$, where $g \in
C_{c}(G)$ \cite{renault}. 
\begin{theorem}[{\cite{timmer:cpmu}}] \label{theorem:groupoid-vg-legs}
\begin{enumerate}
\item $\hA(V_{G})=m(C_{0}(G))$ and
  $\big(\hDelta_{V_{G}}(m(f))\omega\big)(x,y) = f(xy) \omega(x,y)$
  for each $f \in C_{0}(G)$, $\omega \in
  L^{2}(G^{2}_{s,r},\nu^{2}_{s,r})$, $(x,y) \in G^{2}_{s,r}$.
  \item   $A(V_{G})=C^{*}_{r}(G)$, and  for each $g \in
    C_{c}(G)$, $\omega \in L^{2}(G^{2}_{r,r},\nu^{2}_{r,r})$, $(x,y) \in G^{2}_{r,r}$,
    \begin{align*} \hspace{-1ex}
      \big(\Delta_{V_{G}}(L(g))\omega\big)(x,y) =
      \int_{G^{r_{G}(x)}} g(z)D(z)^{-1/2}
      \omega(z^{-1}x,z^{-1}y)
      d\lambda^{r_{G}(x)}(z).
    \end{align*}
  \end{enumerate}
\end{theorem}

\paragraph{The  function algebra of a compact
  groupoid}

Let $G$ be a locally compact groupoid as before  but assume that $G$
is compact. Then the following assertions are evident: 
\begin{lemma} \label{lemma:groupoid-compact}
  \begin{enumerate}
  \item There exists a compact $C^{*}$-quantum graph
    $(C(G^{0}),\mu,C(G),r,\phi,s,\psi,D^{-1})$ with coinvolution $R$ such
    that
    \begin{align*}
      (r(f))(x) &= f(r_{G}(x)), & (s(f))(x) &= f(s_{G}(x))
      &&\text{ for all } f \in C(G^{0}), \, x \in G, \\
      (\phi(g))(u) &= \int_{G^{u}} g(y) d\lambda^{u}(y), &
      (\psi(g))(u) &= \int_{G_{u}} g(y) d\lambda^{-1}_{u}(y) &&
      \text{ for all } g\in C(G), \, u \in G^{0},
    \end{align*}
    and $(R(g))(x)=g(x^{-1})$ for all $g \in C(G)$, $x \in G$.  The
    states $\nu=\mu \circ \phi$ and $\nu^{-1}=\mu \circ \psi$ are
    given by $\nu(g) = \int_{G} g d\nu_{G}$ and $\nu^{-1}(g) =
    \int_{G} gd\nu^{-1}_{G}$ for all $g \in C(G)$.
  \item If we identify $H=H_{\nu}$ with $K=L^{2}(G,\nu)$ via
    $f\zeta_{\nu} \equiv f$ for all $f\in C(G)$, then
    $j(f)=f\zeta_{\phi}$ and $\hat j(f)=f\zeta_{\psi}$ for all $f\in
    C(G)$, and 
    $(H,\halpha,\beta,\hbeta,\alpha)=(K,\rho,\rho,\sigma,\sigma)$.
    \qed
  \end{enumerate}
\end{lemma}
With respect to the canonical identification $H=H_{\nu}\cong
L^{2}(G,\nu)=K$, the representation $m \colon C(G) \to {\cal
  L}(K)\cong {\cal L}(H)$ corresponds to the GNS-representation for
$\nu$. We identify $C(G)$ with $m(C(G)) \subseteq {\cal
  L}(K)={\cal L}(H)$ via $m$. 
\begin{theorem} \label{theorem:groupoid-compact}
  \begin{enumerate}
  \item $(C(G^{0}),\mu,C(G),r,\phi,s,\psi,D^{-1},\hDelta_{V_{G}},R)$ is a
    compact $C^{*}$-quantum group\-oid.
  \item The associated $C^{*}$-pseudo-multiplicative unitary
    $(H,\alpha,\beta,\halpha,W)$ equals
    $(K,\sigma,\rho,\rho,V_{G})$.
  \item The associated weak $C^{*}$-pseudo-Kac system
    $(H,\alpha,\halpha,\beta,\hbeta,U,V)$ is a $C^{*}$-pseudo-Kac
    system.
  \item $\Ad_{U}$ defines an isomorphism between the dual Hopf
    $C^{*}$-bimodule $(\hA(V)^{\hbeta,\alpha}_{H},\hDelta_{V})$ and
    $(C^{*}_{r}(G)_{K}^{\rho,\rho}, \Delta_{V_{G}})$.
  \end{enumerate}
\end{theorem}
\begin{proof}
  i) Put $\Delta:=\hDelta_{V_{G}}$.  By Theorems
  \ref{theorem:groupoid-vg}, \ref{theorem:groupoid-vg-legs},
  \cite[Theorem 4.14]{timmer:cpmu}, and Lemma
  \ref{lemma:groupoid-compact}, $(C(G)^{\alpha,\beta}_{H},\Delta) =
  (\hA(V_{G})^{\sigma,\rho}_{K},\hDelta_{V_{G}})$ is a Hopf
  $C^{*}$-bimodule and $[\Delta(C(G))\kalpha{1} ] =
  [\kalpha{1}C(G)]$,
  $[\Delta(C(G))\kbeta{2}]=[\kbeta{2}C(G)]$.

  By Lemma \ref{lemma:groupoid-compact} iii) and Proposition
  \ref{proposition:groupoid-fixed} ii), $\zeta_{\phi} = j(1_{G}) \in
  \Cofix(V_{G})$. Remark \ref{remark:fixed-assumption} shows that
  $[\Delta(C(G))|\zeta_{\phi}\rangle_{2}C(G)]=[\kbeta{2}C(G)]$. Moreover,
  $\zeta_{\psi}=U\zeta_{\phi} \in \Fix(\checkV_{G})$ by Lemma
  \ref{lemma:add-fixed-kac}, $C(G)=\hA(V_{G})=A(\checkV_{G})$ by
  equation \eqref{eq:add-legs}, and now a second application of Remark
  \ref{remark:fixed-assumption} shows that
  $[\Delta(C(G))|\zeta_{\psi}\rangle_{1}C(G)]= [\kalpha{1}C(G)]$.

  By Theorem \ref{theorem:fixed-cofixed} and Corollary
  \ref{corollary:fixed-cofixed}, $\phi \colon a \mapsto
  \zeta_{\phi}^{*}a\zeta_{\phi}$ and $\psi \colon a \mapsto
  \zeta_{\phi}^{*}UaU\zeta_{\phi}$ are a bounded left and a bounded
  right Haar weight for $(C(G)^{\alpha,\beta}_{H},\Delta)$,
  respectively.

  Finally, we show that  the strong invariance condition iii) in Definition
  \ref{definition:cqg} holds. For all $f,g  \in C(G)$,  the operator
  $\langle \zeta_{\psi}|_{1} \Delta(f)(g^{op}
    \stensor{\alpha}{\beta} 1)|\zeta_{\psi}\rangle_{1}$
    is given by pointwise multiplication by the function
    \begin{align*}
      H_{f,g} \colon G \to \complex, \
      y \mapsto \int_{G_{r_{G}(y)}} f(xy)g(x) d\lambda^{-1}_{r_{G}(y)}(x),
    \end{align*}
    and by right-invariance of $\lambda^{-1}$, 
    \begin{align*}
      \big(R(H_{g,f})\big)(y)  =
      H_{g,f}(y^{-1}) &=
      \int_{G_{r_{G}(y^{-1})}} g(xy^{-1})f(x)
      d\lambda^{-1}_{r_{G}(y^{-1})}(x) \\
      &= \int_{G_{r_{G}(y)}} g(x') f(x'y) d\lambda^{-1}_{r_{G}(y)}
      (x') = H_{f,g}(y) \quad \text{for all } y \in G.
    \end{align*}

    ii) With respect to the identifications $H
    \stensor{\alpha}{\beta} H=K \stensor{\sigma}{\rho} K \cong
    L^{2}(G^{2}_{s,r},\nu^{2}_{s,r})$ and $H
    \stensor{\beta}{\halpha} H = K \stensor{\rho}{\rho} K \cong
    L^{2}(G^{2}_{r,r},\nu^{2}_{r,r})$,
    \begin{align*}
      \big(W^{*}|j(g)\rangle_{2} f\big)(x,y) =
      \big(\Delta(g)|\zeta_{\phi}\rangle_{2}f\big)(x,y) =f(x) g(xy)
    \end{align*}
    for all $(x,y) \in G^{2}_{s,r}$ and $f,g \in C(G)$ and hence
    $(W^{*}\omega)(x,y) = \omega(x,xy) = (V_{G}^{*}\omega)(x,y)$ for
    all $\omega \in L^{2}(G^{2}_{s,r},\nu^{2}_{s,r})$.

    iii) Since $C(G)$ is commutative, $J\xi = \overline{\xi}$,
    and  $(U\xi)(x)=
    (IJ\xi)(x)=\xi(x^{-1})D(x)^{-1/2}=(U_{G}\xi)(x)$ for all $\xi
    \in L^{2}(G,\nu_{G})$ and $x\in G$.  By Theorem
    \ref{theorem:groupoid-vg},
    $(K,\rho,\sigma,\rho,\sigma,U_{G},V_{G})$ is a
    $C^{*}$-pseudo-Kac system, and since $V_{G}=W=\hatV$ by ii) and
    Lemma \ref{lemma:dual-w-vhat}, we can use \cite[Proposition
    5.5]{timmer:ckac} to conclude that
    $(H,\alpha,\halpha,\beta,\hbeta,U,V)$ is a $C^{*}$-pseudo-Kac
    system.

    iv) This assertion follows from the relation $\hatV=V_{G}$ (see
    above), equation \eqref{eq:add-legs}, and Theorem
    \ref{theorem:groupoid-vg-legs}.
\end{proof}

\paragraph{The groupoid $C^{*}$-algebra of an \'etale groupoid} 
Let $G$ be a locally compact groupoid as above but assume that $G$
is $r$-discrete and that $\lambda$ is the family of counting
measures. Since $G^{0} \subseteq G$ is closed and open,  we can embed $C(G^{0})$ in $C(G)$ by extending each function
by $0$ outside of $G^{0}$.  Thus, $1_{G^{0}}$ gets identified with
the characteristic function of $G^{0}$.   Denote by $r,s \colon
C(G^{0}) \to C(G)$ the transpose of the range
and the source map $r_{G}$ and $s_{G}$, respectively.
\begin{lemma} \label{lemma:groupoid-etale}
  \begin{enumerate}
  \item There exists a compact $C^{*}$-quantum graph
    $(C(G^{0}),\mu,C^{*}_{r}(G),\iota, \phi, \iota, \phi, 1)$ with
    coinvolution $R$ such that
    \begin{gather*}
      \begin{aligned}
        \iota(f) &=L(f) \text{ for each } f \in C(G^{0}), &
        \big(\phi(L(g))\big)(u) &= g(u) \text{ for each } g \in
        C_{c}(G), \, u \in G^{0}, 
      \end{aligned}  \\
      RL(f) = L(f^{\dag}), \text{ where } f^{\dag}(x)=f(x^{-1})
      \text{ for all } x \in G, \, f \in C_{c}(G).
    \end{gather*}
    The state $\nu=\mu \circ \phi$ is given by $\nu(a)=\langle
    1_{G^{0}}|a1_{G^{0}}\rangle$ for all $a \in C^{*}_{r}(G)$, and
    its modular automorphism group is given by
    $\sigma^{\nu}_{t}(L(f)) = L(D^{it}f)$ for all $f \in C_{c}(G)$, $t
    \in \reals$.
  \item There exists an isomorphism $\Xi \colon H_{\nu} \to K$,
    $L(f)\zeta_{\nu} \mapsto fD^{-1/2}$, and $\Xi
    L(f)^{op}\zeta_{\nu} = f$,  $\Xi L(f)\zeta_{\phi} =
    \hat j(f)$,
    $\Xi L(f)^{op}\zeta_{\phi} =  j(f)$ for all $f\in C_{c}(G)$. In
    particular, $\Xi \halpha=\Xi \hbeta = \sigma$ and $\Xi \alpha =
    \Xi \beta = \rho$.
  \end{enumerate}
\end{lemma}
\begin{proof}
  i) The $*$-homomorphism $\iota$ clearly is well-defined. Denote by
  $\zeta \colon L^{2}(G^{0},\mu_{G}) \to L^{2}(G,\nu_{G})$ the
  embedding that extends each function outside of $G^{0}$ by
  $0$. Then for each $g \in C_{c}(G)$, the operator $\zeta^{*}L(g)\zeta
  \subseteq {\cal L}(L^{2}(G^{0},\mu_{G}))$ is given by pointwise
  multiplication by the function $g|_{G^{0}}$, and we can define
  $\phi \colon C^{*}_{r}(G) \to C(G^{0}) \subseteq {\cal
    L}(L^{2}(G^{0},\mu_{G}))$ by $a \mapsto
  \zeta^{*}a\zeta$. Clearly, $\iota \circ \phi \colon C^{*}_{r}(G)
  \to \iota(C(G^{0}))$ is a conditional expectation.  Since
  $\zeta\zeta_{\mu} =1_{G^{0}}$, the state $\mu \circ \phi$ is given
  by $\nu(a) = \mu\big(\zeta^{*}a\zeta\big) = \langle
  1_{G^{0}}|a1_{G^{0}}\rangle$ for all $a \in C^{*}_{r}(G)$.  By
  \cite[\S II.5]{renault}, this is a KMS-state with modular
  automorphism group as stated above --- indeed, for all $f \in
  C_{c}(G)$,
  \begin{align*}
    \nu(L(f)^{*}L(f)) &= \int_{G^{0}} \int_{G^{u}} \overline{f(x^{-1})} f(x^{-1})
    d\lambda^{u}(x) d\mu_{G}(u) \\
    &=
    \int_{G}  \overline{f(x)}f(x) D(x)^{-1} d\nu_{G}(x) 
    = \nu(L(D^{-1/2}f)L(D^{-1/2}f)^{*}).
  \end{align*}
  Finally, $\sigma^{\nu}_{t} \circ \iota=\iota$ and $\phi \circ
  \sigma^{\nu}_{t} = \phi$ for all $t \in \reals$ because $D|_{G^{0}} =
  1_{G^{0}}$.

  ii) First, observe that for all $f \in C_{c}(G)$ and $x\in G$
\begin{align} \label{eq:groupoit-etale-embed}
  \big(L(f)1_{G^{0}}\big)(x) = f(x)D^{-1/2}(x)
\end{align}
and hence $\| L(f)\zeta_{\nu}\|_{H}^{2} = \nu(L(f)^{*}L(f)) =
\langle L(f) 1_{G^{0}}|L(f)1_{G^{0}}\rangle = \|fD^{-1/2}\|^{2}_{K}
$.  The existence of $\Xi$ follows. The remaining assertions hold
because for all $f \in C_{c}(G)$, $g \in C(G^{0})$,
  \begin{gather*}
    \Xi L(f)^{op}\zeta_{\nu} = \Xi
    \sigma^{\nu}_{-i/2}(L(f))\zeta_{\nu} = \Xi
    L(D^{1/2}f)\zeta_{\nu} = f, \\
    \Xi L(f)\zeta_{\phi} g\zeta_{\mu} = \Xi L(f)L(g)\zeta_{\nu} =
    \Xi L(f s(g)) \zeta_{\nu} = fs(g)D^{-1/2} = \hat j(f)g, \\
    \Xi L(f)^{op} \zeta_{\phi} g\zeta_{\mu} = \Xi
    L(g)L(f)^{op}\zeta_{\nu} =\Xi r(g) L(D^{1/2}f)\zeta_{\nu} =
    r(g)f = j(f)g. \qedhere
  \end{gather*}
\end{proof}
From now on, we identify $H=H_{\nu}$ with $K$ via $\Xi$ without
further notice.
\begin{lemma}
  $\big(L(f)^{op}g\big)(x)=\int_{G^{s_{G}(x)}} g(xy)f(y^{-1})
  d\lambda^{s_{G}(x)}(y)$ for all $f,g \in C_{c}(G)$, $x \in G$.
\end{lemma}
\begin{proof}
  Let $f,g \in C_{c}(G)$. Then
  \begin{align*}
    L(f)^{op}g = \Xi L(f)^{op}  L(g D^{1/2}) \zeta_{\nu} &= \Xi
    L(D^{1/2}g )\sigma^{\nu}_{-i/2}(L(f)) \zeta_{\nu}  
    \\ &= \Xi L(D^{1/2}g ) L(D^{1/2}f)
    \zeta_{\nu} = \Xi L(h)\zeta_{\nu} = h D^{-1/2},
  \end{align*}
where for all $x \in G$,
\begin{align*}
  h(x) &= \int_{G^{r_{G}(x)}} D^{1/2}(z) g (z) D^{1/2}(z^{-1}x)f(z^{-1}x)
  d\lambda^{r_{G}(x)}(z) \\ &= D^{1/2}(x) \int_{G^{s_{G}(x)}}
  g (xy)f(y^{-1})d\lambda^{s_{G}(x)}(y). \qedhere
\end{align*}
\end{proof}
\begin{theorem} \label{theorem:groupoid-etale}
  \begin{enumerate}
  \item $(C(G^{0}),\mu,C^{*}_{r}(G), \iota, \phi, \iota, \phi ,
    1,R,\Delta_{V_{G}})$ is a compact $C^{*}$-quantum groupoid.
  \item The associated $C^{*}$-pseudo-multiplicative unitary is $V_{G}$.
  \item The associated weak $C^{*}$-pseudo-Kac system is a
    $C^{*}$-pseudo-Kac system.
  \item The dual Hopf $C^{*}$-bimodule is
    $(C(G)_{K}^{\sigma,\rho},\hDelta_{V_{G}})$. 
  \end{enumerate}
\end{theorem}
\begin{proof}
  i) Put $\Delta:=\Delta_{V_{G}}$. By Theorems
  \ref{theorem:groupoid-vg}, \ref{theorem:groupoid-vg-legs},
  \cite[Theorem 4.14]{timmer:cpmu}, and Lemma
  \ref{lemma:groupoid-etale},
  $(C^{*}_{r}(G)^{\alpha,\beta}_{H},\Delta) =
  (C^{*}_{r}(G)^{\rho,\rho}_{K},\Delta_{V_{G}})$ is a Hopf
  $C^{*}$-bimodule and
  $[\Delta(C^{*}_{r}(G))\kalpha{1}]=[\kalpha{1}C^{*}_{r}(G)]$,
  $[\Delta(C^{*}_{r}(G))\kbeta{2}]=[\kbeta{2}C^{*}_{r}(G)]$. 

  By Lemma \ref{lemma:groupoid-etale} and Proposition
  \ref{proposition:groupoid-fixed}, $1_{C^{*}_{r}(G)}\zeta_{\phi} =
  L(1_{G^{0}})\zeta_{\phi} = j(1_{G^{0}}) = \hat j(1_{G^{0}}) \in
  \Fix(V_{G})$.  Remark \ref{remark:fixed-assumption} shows that
  $[\Delta(C^{*}_{r}(G))|\zeta_{\phi}\rangle_{1}C^{*}_{r}(G)]=[\kalpha{1}C^{*}_{r}(G)]$. Moreover,
  $\zeta_{\phi} \in \Cofix(\hatV_{G})$ by Lemma
  \ref{lemma:add-fixed-kac}, $C^{*}_{r}(G)=A(V_{G})=\hA(\hatV_{G})$
  by equation \eqref{eq:add-legs}, and now a second application of
  Remark \ref{remark:fixed-assumption} shows that
  $[\Delta(C^{*}_{r}(G))|\zeta_{\phi}\rangle_{2}C^{*}_{r}(G)]=
  [\kbeta{2}C^{*}_{r}(G)]$.

  The map $\phi \colon a \mapsto \zeta_{\phi}^{*}a\zeta_{\phi}$ is a
  bounded left and a bounded right Haar weight for
  $(C^{*}_{r}(G)^{\alpha,\beta}_{H},\Delta)$ by Theorem
  \ref{theorem:fixed-cofixed} and Corollary
  \ref{corollary:fixed-cofixed}.

  Finally, let us prove that the strong invariance condition iii) in
  Definition \ref{definition:cqg} holds. Let $f,g,\xi \in
  C_{c}(G)$. By the previous Lemma,
  \begin{align*}
    \big((L(g^{op})
    \stensor{\alpha}{\beta} 1)|\zeta_{\psi}\rangle_{1}\xi\big)(x,y)
    = \xi(y)g(x) \quad \text{for all } (x,y) \in G^{2}_{r,r}
  \end{align*}
  and hence
  \begin{align*}
    \big( \langle \zeta_{\psi}|_{1} \Delta(L(f))(L(g^{op})
    \stensor{\alpha}{\beta} 1)|\zeta_{\psi}\rangle_{1}\xi\big)(y)
    &= \int_{G^{r_{G}(y)}} f(z)D^{-1/2}(z)g(z^{-1})\xi(z^{-1}y)
    d^{r_{G}(y)}(z)   \\
    &= \big(L(h)\xi\big)(y),
  \end{align*}
  where $h(z)=f(z)g(z^{-1})$ for all $z \in G$. Switching $f,g$, we find
  \begin{align*}
    \big( \langle \zeta_{\psi}|_{1} \Delta(L(g))(L(f^{op})
    \stensor{\alpha}{\beta} 1)|\zeta_{\psi}\rangle_{1}\big)  =
    L(h^{\dag}) = R(L(h)).
  \end{align*}
  Since $f,g$ were arbitrary,  condition  iii) in
  Definition \ref{definition:cqg} holds.

  ii) By Lemma \ref{lemma:groupoid-etale}, we have for all $f,g \in
  C_{c}(G)$, $(x,y) \in G^{2}_{r,r}$,
\begin{align*}
  \big(V(L(f)\zeta_{\psi} \tr L(g)\zeta_{\nu})\big)(x,y) &
  = \big( \Delta(L(f))(\zeta_{\psi} \tr L(g)\zeta_{\nu})\big)(x,y) \\
  &=\int_{G^{r_{G}(x)}} f(z)D^{-1/2}(z) \big(\zeta_{\psi} \tr
  L(g)\zeta_{\nu}\big)(z^{-1}x,z^{-1}y)d\lambda^{r_{G}(x)}(z)
  \\
  &= \int_{G^{r_{G}(x)}} f(z)D^{-1/2}(z) 1_{G^{0}}(z^{-1}x)
  g(z^{-1}y)D^{-1/2}(z^{-1}y)d\lambda^{r_{G}(x)}(z) \\
  &= f(x)D^{-1/2}(x)
  g(x^{-1}y)D^{-1/2}(x^{-1}y)  
\\  &= (\hat j(f) \tr gD^{-1/2})(x,x^{-1}y) = \big(L(f)\zeta_{\psi} \tr
  L(g)\zeta_{\nu}\big)(x,x^{-1}y). 
\end{align*}

  iii), iv) The proof is similar to the proof of  Theorem
  \ref{theorem:groupoid-compact} iii), iv).
\end{proof}

\def\cprime{$'$}


\begin{thebibliography}{10}

\bibitem{baaj:2}
S.~Baaj and G.~Skandalis.
\newblock Unitaires multiplicatifs et dualit\'e pour les produits crois\'es de
  {$C\sp *$}-alg\`ebres.
\newblock {\em Ann. Sci. \'Ecole Norm. Sup. (4)}, 26(4):425--488, 1993.

\bibitem{blanchard}
{\'E}.~Blanchard.
\newblock D\'eformations de {$C\sp *$}-alg\`ebres de {H}opf.
\newblock {\em Bull. Soc. Math. France}, 124(1):141--215, 1996.

\bibitem{bratteli:2}
O.~Bratteli and D.~W. Robinson.
\newblock {\em Operator algebras and quantum statistical mechanics. 2}.
\newblock Texts and Monographs in Physics. Springer-Verlag, Berlin, second
  edition, 1997.
\newblock Equilibrium states. Models in quantum statistical mechanics.

\bibitem{enock:9}
M.~Enock.
\newblock Inclusions of von {N}eumann algebras and quantum groupo\"\i ds.
  {III}.
\newblock {\em J. Funct. Anal.}, 223(2):311--364, 2005.

\bibitem{enock:lesieur}
M.~Enock.
\newblock On {L}esieur's quantum groupoids, June 2007.

\bibitem{enock:1}
M.~Enock and J.-M. Vallin.
\newblock Inclusions of von {N}eumann algebras, and quantum groupoids.
\newblock {\em J. Funct. Anal.}, 172(2):249--300, 2000.

\bibitem{enock:2}
M.~Enock and J.-M. Vallin.
\newblock Inclusions of von {N}eumann algebras and quantum groupoids. {II}.
\newblock {\em J. Funct. Anal.}, 178(1):156--225, 2000.

\bibitem{lance}
E.~C. Lance.
\newblock {\em Hilbert {$C\sp *$}-modules. A toolkit for operator algebraists},
  volume 210 of {\em London Mathematical Society Lecture Note Series}.
\newblock Cambridge University Press, Cambridge, 1995.

\bibitem{lesieur}
F.~Lesieur.
\newblock {\em Groupo\"ides quantiques mesur\'es: axiomatique, \'etude,
  dualit\'e, exemples.}
\newblock PhD thesis, Universite d'Orl\'eans, 2003.
\newblock Available at
  {\texttt{www.univ-orleans.fr/mapmo/publications/lesieur/these.php}}.

\bibitem{ouchi}
M.~O'uchi.
\newblock Pseudo-multiplicative unitaries on {H}ilbert {$C\sp \ast$}-modules.
\newblock {\em Far East J. Math. Sci. (FJMS)}, Special Volume, Part
  II:229--249, 2001.

\bibitem{paterson}
A.~L.~T. Paterson.
\newblock {\em Groupoids, inverse semigroups, and their operator algebras},
  volume 170 of {\em Progress in Mathematics}.
\newblock Birkh\"auser Boston Inc., Boston, MA, 1999.

\bibitem{pedersen}
G.~K. Pedersen.
\newblock {\em {$C\sp{\ast} $}-algebras and their automorphism groups},
  volume~14 of {\em London Mathematical Society Monographs}.
\newblock Academic Press Inc. [Harcourt Brace Jovanovich Publishers], London,
  1979.

\bibitem{renault}
J.~Renault.
\newblock {\em A groupoid approach to {$C\sp{\ast} $}-algebras}, volume 793 of
  {\em Lecture Notes in Mathematics}.
\newblock Springer, Berlin, 1980.

\bibitem{sauvageot:2}
J.-L. Sauvageot.
\newblock Produits tensoriels de {${Z}$}-modules et applications.
\newblock In {\em Operator algebras and their connections with topology and
  ergodic theory (Bu\c steni, 1983)}, volume 1132 of {\em Lecture Notes in
  Mathematics}, pages 468--485. Springer, Berlin, 1985.

\bibitem{takesaki:2}
M.~Takesaki.
\newblock {\em Theory of operator algebras. {II}}, volume 125 of {\em
  Encyclopaedia of Mathematical Sciences}.
\newblock Springer-Verlag, Berlin, 2003.

\bibitem{timmer:ckac}
T.~Timmermann.
\newblock {$C^*$}-pseudo-{K}ac systems and duality for coactions of concrete
  {H}opf {$C^*$}-bimodules.
\newblock Technical Report 484, SFB 478, WW-University Muenster, October 2007.
\newblock arXiv:0709.4617.

\bibitem{timmer:cpmu}
T.~Timmermann.
\newblock {$C^*$}-pseudo-multiplicative unitaries.
\newblock Technical Report 481, SFB 478, WW-University Muenster, September
  2007.
\newblock arXiv:0709.2995.

\bibitem{timmermann:hopf}
T.~Timmermann.
\newblock Pseudo-multiplicative unitaries on {$C^*$}-modules and {H}opf
  {$C^*$}-families.
\newblock {\em J. {N}oncomm. {G}eom.}, 1:497--542, 2007.

\bibitem{timmer:qg}
T.~Timmermann.
\newblock {\em An invitation to quantum groups and duality}.
\newblock EMS Publishing House, 2008.

\bibitem{vallin:1}
J.-M. Vallin.
\newblock Bimodules de {H}opf et poids op\'eratoriels de {H}aar.
\newblock {\em J. Operator Theory}, 35(1):39--65, 1996.

\bibitem{vallin:2}
J.-M. Vallin.
\newblock Unitaire pseudo-multiplicatif associ\'e \`a un groupo\"\i de.
  {A}pplications \`a la moyennabilit\'e.
\newblock {\em J. Operator Theory}, 44(2):347--368, 2000.

\bibitem{woron}
S.~L. Woronowicz.
\newblock From multiplicative unitaries to quantum groups.
\newblock {\em Internat. J. Math.}, 7(1):127--149, 1996.

\end{thebibliography}
\end{document}